\documentclass[10pt]{iopart}
\usepackage{graphicx}
\usepackage{iopams}
\usepackage{amsopn}
\usepackage{algorithm}
\usepackage{algpseudocode}

\expandafter\let\csname equation*\endcsname\relax
\expandafter\let\csname endequation*\endcsname\relax
\usepackage{amsmath}

\newcommand{\diag}[1]{{\rm diag}\left( #1 \right)}
\newcommand{\stack}[2]{\left( \begin{array}{c} #1 \\ #2 \end{array} \right)}
\DeclareMathOperator*{\argmax}{argmax\,}
\DeclareMathOperator*{\argmin}{argmin\,}

\DeclareMathOperator{\proj}{proj}
\DeclareMathOperator{\prox}{prox}

\DeclareMathOperator{\sh}{shrink}
\DeclareMathOperator*{\solve}{solve}

\begin{document}

\title[CPPD]{Notes on the primal-dual algorithm for convex optimization applied to X-ray tomographic image reconstruction}

\author{Emil Y. Sidky and X. Pan}

\address{Department of Radiology, University of Chicago, MC-2026, 5841 S. Maryland Ave., Chicago, Illinois 60637}
\ead{sidky@uchicago.edu, xpan@uchicago.edu}
\vspace{10pt}
\begin{indented}
\item[]\today
\end{indented}

\begin{abstract}
The purpose of these notes is to provide background on understanding the primal-dual
algorithm of Chambolle and Pock \cite{chambolle2011first} for imaging scientists.
The presentation focuses on providing intuition and an algorithmic system that is
amenable to pre-conditioning. The document aims to be self-contained, providing
background on the essential facts of non-smooth convex analysis.\cite{rockafellar1970convex}
\end{abstract}

%
%
%
%
%

\section{Introduction}
In Ref. \cite{SidkyCP:12}, we demonstrated
how to implement the primal-dual (CPPD) algorithm of Chambolle and Pock \cite{chambolle2011first} for solving
non-smooth convex optimization problems of interest to X-ray tomographic image reconstruction.
The CPPD algorithm has proven useful in that it allows a great deal of flexibility in formulating convex optimization
problems that allow for a variety of data models, regularization terms, and image constraints.
These notes focus on a simplified CPPD system, i.e. no $G(x)$ function, which allows for a simple preconditioning technique.

Future updates of this document will focus on including the $G(x)$ function, how to effectively deal with system matrices composed of multiple blocks,
and the connection with the alternating direction method of multipliers (ADMM) algorithm \cite{boyd2011distributed}. The
connection to ADMM is useful for developing solvers for non-convex optimization problems \cite{barber2024convergence}.

The main body of these notes presents a fairly direct narrative toward the goal of writing a CPPD algorithm for total variation (TV)
regularized least-squaresimage reconstruction, which includes preconditioning. The appendices have supporting discussions
on some of the convex analysis concepts, additional derivations, and numerical experiments. Please contact us if you have any questions,
comments, or spotted errors.

\section{Methods}

\subsection{CPPD background}
\label{sec:CPPDback}

The CPPD algorithm is designed to solve the generic convex optimization problem
\begin{equation}
\label{mainmin}
x^\star = \argmin_x F(A x),
\end{equation}
where $x$ is a $n$-dimensional vector; $A$ is a $m \times n$; $F(\cdot)$ is a simple,
convex function that is possibly non-smooth. This problem is actually a special
case of what is considered in Ref. \cite{chambolle2011first}, which also considered
an additional convex term $G(x)$. The minimization in Eq. (\ref{mainmin}) encompasses
many optimization problems of interest for X-ray CT.

The difficulty in solving Eq. (\ref{mainmin}) results from the composition of $F$ with
a large linear transform $A$, where large means that it is only computationally feasible
to calculate matrix-vector products such as $Ax$ and matrix operations on $A$, such
as diagonalization, is not feasible.  The CPPD algorithm results from
converting Eq. (\ref{mainmin}) to a saddle point problem, where the convex function $F$
appears in a different than the linear transform $A$. The derivation of the saddle point
problem goes as follows:

The minimization in Eq. (\ref{mainmin}) is equivalent to the equality-constrained minimization
\begin{equation*}
\min_{x,y} F(y) \text{ such that } y=Ax,
\end{equation*}
where the splitting variable, $y$, is a $m$-dimensional vector.
This in turn can be converted to an unconstrained saddle point problem by forming
the Lagrangian $L$ 
\begin{equation*}
L(x,y,\lambda) = F(y) + \lambda^\top(Ax-y),
\end{equation*}
where $\lambda$, also a $m$-dimensional vector, is the Lagrange multiplier or dual variable.
The saddle point problem optimizing $L$ is
\begin{equation}
\label{saddle1}
\min_{x,y} \max_\lambda L(x,y,\lambda),
\end{equation}
and the solution can be identified formally by setting the gradient of $L$ to zero,
see \ref{app:critical},
\begin{align}
\partial_x L(x,y,\lambda) &= A^\top \lambda = 0, \label{cond1} \\
\partial_y L(x,y,\lambda) &= \partial F(y) - \lambda  = 0 \label{cond2}, \\
\partial_\lambda L(x,y,\lambda) &= Ax - y = 0, \label{cond3}.
\end{align}
We will use
the first and third of these equations to provide convergence checks for the CPPD algorithm.
The widely used alternating direction method of multipliers (ADMM) algorithm solves this system
of equations
with update steps derived from a modified Lagrangian, where a term quadratic in $\|Ax-y\|$
is added to the Lagrangian $L(x,y,\lambda)$ \cite{boyd2011distributed}.

For the CPPD algorithm, the size of the saddle point problem in Eq. (\ref{saddle1}) is reduced
by carrying out the minimization over $y$.
The second condition for the solution of this problem, Eq. (\ref{cond2}),
can be solved directly if $F$ is a simple
function.
The expression in Eq. (\ref{cond2}) originated from the gradient of the terms in $L(x,y,\lambda)$
that involved $y$. Isolating the minimization over $y$ from Eq. (\ref{saddle1}) yields
\begin{equation*}
\min_y \left\{ F(y) - \lambda^\top y \right\},
\end{equation*}
which is essentially the Legendre-Fenchel transform of $F(y)$, see \ref{app:legendre}:
\begin{equation*}
\min_y \left\{ F(y) - \lambda^\top y \right\} =
- \max_y \left\{ \lambda^\top y -F(y) \right\} \equiv - F^*(\lambda).
\end{equation*}
Recall that for convex $F(y)$
\begin{align*}
F^*(\lambda) &\equiv \max_y  \left\{ \lambda^\top y -F(y) \right\}, \\
F(y) &\equiv \max_\lambda  \left\{ y^\top \lambda -F^*(\lambda) \right\}.
\end{align*}
If $F^*(\lambda)$ is easily computed, the saddle point problem in Eq. (\ref{saddle1})
can be reduced to
\begin{equation}
\label{saddle2}
\min_x \max_\lambda \left\{ \lambda^\top A x  -F^*(\lambda) \right\}.
\end{equation}
For many optimization problems of interest for CT image reconstruction, $F^*$ can
be derived analytically.

\subsection{Heuristics for the CPPD algorithm and preconditioning}
\label{sec:heuristics}

The CPPD algorithm solves the saddle point problem in Eq. (\ref{saddle2}), but
the form of the algorithm is difficult to understand at an intuitive level.
In addition to the development of this section, please refer also to \ref{app:intuition}.
Intuition on the algorithm can be acquired by focusing only on the saddle
term 
\begin{equation*}
s(x,\lambda) = \lambda^\top A x.
\end{equation*}
The critical point of this potential is found by setting its gradient to zero. If
$A$ is full-rank and invertible, there is only one critical point
at $x_\text{crit} = 0, \lambda_\text{crit} = 0$. More generally, the critical
points satisfy
\begin{eqnarray*}
A x_\text{crit} = 0, \\
A^\top \lambda_\text{crit}=0.
\end{eqnarray*}

That the critical point
is a saddle point is seen by computing the Hessian
\begin{equation*}
H_s = 
\left(
\begin{array}{cc}
\frac{\partial^2}{\partial x^2} s(x,\lambda) & \frac{\partial^2}{\partial \lambda \partial x} s(x,\lambda)   \\
   & \\
\frac{\partial^2}{\partial x \partial \lambda} s(x,\lambda) &  \frac{\partial^2}{\partial \lambda^2} s(x,\lambda)  
\end{array} \right) = 
\left(
\begin{array}{cc}
0 & A^\top \\
A & 0
\end{array} \right)  .
\end{equation*}
The trace of the Hessian is zero; therefore the sum of eigenvalues is zero.
As a result, there must be negative and positive eigenvalues,
meaning that $s(x,\lambda)$ has directions of negative and positive curvature.
Thus, $(x_\text{crit},\lambda_\text{crit})$ is a saddle point of $s(x,\lambda)$.

\subsubsection{Forward Euler iteration}
\label{sec:fei}
A saddle point solver should be able to find the critical saddle point(s) from arbitrary initialization
$(x_0,\lambda_0)$. Because Eq. (\ref{saddle2}) calls for minimization over $x$ and maximization
over $\lambda$, a first attempt at an algorithm might involve taking a step in the direction
of $-\partial_x s$ and $\partial_\lambda s$. Forming the forward Euler iteration yields
\begin{eqnarray}
\label{forwardEuler1}
x_{k+1} = x_k - \alpha A^\top \lambda_k, \\
\label{forwardEuler2}
\lambda_{k+1} = \lambda_k + \alpha A x_k, 
\end{eqnarray}
where $\alpha$ is a step-size parameter.
This iteration, however, is unstable for all step-sizes $\alpha$; this can be shown
by computing the magnitude of the solution estimate, $\|(x_{k+1},\lambda_{k+1})^\top\|$,
in terms of $\|(x_{k},\lambda_{k})^\top\|$. The ratio of the former to the latter is greater
than one for any $\alpha=0$. Thus, the iteration in Eqs. (\ref{forwardEuler1}) and (\ref{forwardEuler2})
spirals away from the critical saddle point(s) unless it is initialized with the solution, i.e.
a critical saddle point.

\subsubsection{Backward Euler iteration}
\label{sec:bei}
The second attempt at an algorithm is to use backward Euler iteration, where the step direction
is computed based on the gradients at the updated point instead. For the $s(x,\lambda)$ 
saddle-point solver this update is written
\begin{eqnarray}
\label{backwardEuler1}
x_{k+1} = x_k - \alpha A^\top \lambda_{k+1}, \\
\label{backwardEuler2}
\lambda_{k+1} = \lambda_k + \alpha A x_{k+1}, 
\end{eqnarray}
where the difference with the forward Euler iteration is in the last terms of both
update equations. Those terms use $\lambda_{k+1}$ and $x_{k+1}$ instead of $\lambda_k$ and $x_k$.
Bringing all terms involving $k+1$ to the left-side, the backward Euler update equations become
\begin{equation}
\label{backwardEulerMatrix}
\left(
\begin{array}{cc}
1  & \alpha A^\top \\
- \alpha A & 1
\end{array}
\right)
\left(
\begin{array}{c}
 x_{k+1}\\
\lambda_{k+1} 
\end{array}
\right) =
\left(
\begin{array}{c}
 x_{k}\\
\lambda_{k} 
\end{array}
\right).
\end{equation}
This iteration can be shown to converge to a critical point of $s(x,\lambda)$ for
any positive step-size, $\alpha$. The problem
with this algorithm, however, is that each update computation involves solving the linear
system in Eq. (\ref{backwardEulerMatrix}). If $A$ is a large matrix, this computation may have
to be addressed with an iterative linear systems solver, and thus the whole algorithm would
involve a "loop inside of a loop", which could be computationally infeasible.

\subsubsection{Approximate backward Euler iteration}
\label{sec:abei}
One of the key insights of the CPPD algorithm \cite{chambolle2011first}
is to use the backward Euler
step in an approximate fashion by extrapolating $x$ from the previous two iterations
\begin{eqnarray}
\label{cps1A}
\bar{x}_{k+1} =& x_k  + \theta (x_k - x_{k-1}), \\
\label{cps1B}
\lambda_{k+1} =& \lambda_k + \alpha A \bar{x}_{k+1}, \\
\label{cps1C}
x_{k+1} =& x_k - \alpha A^\top \lambda_{k+1},
\end{eqnarray}
where $\theta \in [0,1]$ is the extrapolation parameter; for the following sections
of the article we employ $\theta=1$, which extrapolates $x_{k+1}$ exactly if
$x_k$ is a linear function of $k$.
By using $\bar{x}_{k+1}$ instead of $x_{k+1}$ in Eq. (\ref{cps1B}), it is no longer
necessary to solve a linear system for each iteration; all update steps can be
computed directly with matrix vector products.

Even though this scheme is reasoned heuristically, its convergence is
rigorously proven in Ref. \cite{chambolle2011first}, but unlike the backward iteration
case where any $\alpha$ yields convergence, the condition on $\alpha$ here is
\begin{equation*}
\alpha < 1/\|A\|_2,
\end{equation*}
where the matrix norm $\|A\|_2$ is the largest singular value of $A$;
in practice, setting $\alpha = 1/\|A\|_2$ usually results in a convergent
iteration. The dual update, Eq. (\ref{cps1B}), and the primal update,
Eq. (\ref{cps1C}), can also employ different stepsizes, $\sigma$ and $\tau$,
respectively, as long as they satisfy
\begin{equation*}
\sigma \tau < 1/\|A\|^2_2,
\end{equation*}
where again this condition, in practice, can include the equality.
Implementing the different step sizes and reordering the steps yields
\begin{eqnarray}
\label{cps2A}
x_{k+1} =& x_k - \tau A^\top \lambda_{k},\\
\label{cps2B}
\bar{x}_{k+1} =& x_{k+1} + \theta (x_{k+1}  - x_{k}), \\
\label{cps2C}
\lambda_{k+1} =& \lambda_k + \sigma A \bar{x}_{k+1}.
\end{eqnarray}
This form is slightly more convenient in that it only requires $x_k$ and $\lambda_k$
to compute values for the $k+1$-iteration.

\subsubsection{Generalization to matrix-mapping steps}
\label{sec:matrix-mapping}

The CPPD iteration is also shown to be a generalization of
the proximal point algorithm \cite{he2012convergence,boyd2011distributed}.
This generalizaition allows for matrix-mapping steps, where $\sigma$ and $\tau$
are replaced with symmetric positive matrices $\Sigma$ and $T$
\begin{eqnarray}
\label{cps3A}
x_{k+1} =& x_k - T A^\top \lambda_{k},\\
\label{cps3B}
\bar{x}_{k+1} =& 2 x_{k+1} - x_{k}, \\
\label{cps3C}
\lambda_{k+1} =& \lambda_k + \Sigma A \bar{x}_{k+1},
\end{eqnarray}
where we only consider $\theta=1$.
For this algorithm to converge, the condition on the matrices $\Sigma$ and $T$ is
arrived at through the matrix
\begin{equation*}
B = \left(
\begin{array}{cc}
T^{-1} & -A^\top \\
-A    & \Sigma^{-1}
\end{array}
\right).
\end{equation*}
For convergence, $B$ should be a positive definite matrix, see \ref{app:fixedpoint}.
Using the Schur complement, this condition on $B$ reduces to either
\begin{eqnarray}
\label{tauineq}
& T^{-1} - A^\top \Sigma A > 0 \\
 \text{ or } & \notag \\
\label{sigineq}
& \Sigma^{-1} - A T A^\top > 0,
\end{eqnarray}
where the inequalities are shorthand for indicating that the matrix on the left
is positive definite. Again, in practice, we can use the equality on either of these
conditions. In particular, take
\begin{equation}
\Sigma^{-1} = A T A^\top \label{matrixCond}.
\end{equation}
This relation between the dual and primal matrix step-mappings proves
useful for designing non-diagonal step-preconditioners.

\subsubsection{Preconditioning based on inverse of $(A^\top A)$}
\label{sec:ndpc}

To demonstrate perfect preconditioning, take the
case
\begin{equation}
\label{taubestdef}
 T= (\rho A^\top A)^{-1},
\end{equation}
 and 
\begin{equation}
\label{sigdef}
\Sigma=\rho \, I,
\end{equation}
where $I$ is the identity matrix; the step-size ratio parameter
$\rho$ is a positive scalar; and 
$A^\top A$ is assumed invertible.
This choice of $T$ and $\Sigma$ satisfies Eq. (\ref{matrixCond}).
Substituiting this $T$ expression into Eq.(\ref{cps3A}) yields
\begin{equation*}
x_{k+1} = x_k - (1/\rho) A^{-1} \lambda_{k};
\end{equation*}
and multiplying through by $A$ and setting $u=Ax$, we have
\begin{equation*}
u_{k+1} = u_k - (1/\rho) \lambda_{k}.
\end{equation*}
Combining Eqs. (\ref{cps3B}), (\ref{cps3C}), and (\ref{sigdef}) yields
\begin{equation*}
\lambda_{k+1} = \lambda_k + \rho ( 2 u_{k+1} -u_k),
\end{equation*}
where again $u=Ax$ is used.
Simplifying the expression for $\lambda_{k+1}$ reduces
the update equations to
\begin{eqnarray*}
u_{k+1} = u_k - (1/\rho) \lambda_{k}, \\
\lambda_{k+1} = \rho (u_k - (1/\rho) \lambda_{k}).
\end{eqnarray*}
Direct computation shows that this iteration will terminate in two
steps for any $u_0=Ax_0$, $\lambda_0$, and $\rho>0$ with the result $u_2=Ax_2=0$
and $\lambda_2=0$. The parameter $\rho$ has no effect on the number of iterations, here,
but it proves useful for the general CPPD iteration.

For large-scale $A$, it may not be feasible or possible to compute $(A^\top A)^{-1}$,
and the convergence in two steps is lost once the more general saddle point problem
in Eq. (\ref{saddle2}) is considered. But the presented argument is a heuristic
for use of an approximate inverse for $T$
\begin{equation*}
T \approx (A^\top A)^{-1}
\end{equation*}
as a non-diagonal step-preconditioner.
If $T$ is a matrix that is only an approximate inverse of $(A^\top A)^{-1}$
it will also be infeasible to compute $\Sigma$ using Eq. (\ref{matrixCond}).
Instead it is more practical to assume $\Sigma$ is diagonal and choose it so
it satisfies the inequality in Eq. (\ref{sigineq}). This is accomplished
by setting the diagonal elements of $\Sigma$ to a value less than the largest
singular value of $ATA^\top$; the largest singular value can be obtained by
the power method. In some cases it is more convenient to do power iteration
with the matrix $T A^\top A$, which has the same largest eigenvalue as the
matrix $ATA^\top$. The step-size ratio parameter $\rho$ can still be used 
when $T$ is an approximate inverse of $(A^\top A)^{-1}$ because it results from
multiplying both sides of the inequality in Eq. (\ref{sigineq}) by $\rho$.



\subsection{Heurstic derivation of CPPD}
\label{sec:cppd}

The update steps in Eqs. (\ref{cps3A}) - (\ref{cps3C}) combined with strategies
for setting the step matrices $\Sigma$ and $T$ only address solution of the bilinear
saddle point optimization $\lambda^\top A x$. The extension to the main problem of interest
Eq. (\ref{saddle2}) involves accounting for the additional concave potential
term $-F^*(\lambda)$. We run through the complete argument with backward Euler iteration
once more using the step matrices from the beginning.

Repeating the saddle problem of interest
\begin{equation*}
\min_x \max_\lambda \left\{ \lambda^\top A x  -F^*(\lambda) \right\},
\end{equation*}
the $x$ and $\lambda$ updates for backward Euler iteration are
\begin{align*}
\lambda_{k+1} =& \lambda_k + \Sigma \left( A x_{k+1} - \partial F^*(\lambda_{k+1}) \right),\\
x_{k+1} =& x_k - T A^\top \lambda_{k+1}
\end{align*}
where $\Sigma$ and $T$ are symmetric positive definite matrices, i.e. positive eigenvalues.
Working toward the CPPD algorithm, $x_{k+1}$ in the $\lambda$ update is replaced by
an extrapolation and the differential term is moved to the left-hand side.
\begin{align}
\bar{x}_{k+1} =& 2 x_k  - x_{k-1}, \notag \\
\Sigma \partial F^*(\lambda_{k+1}) + \lambda_{k+1} =& \lambda_k + \Sigma A \bar{x}_{k+1} 
\label{cppdint1} \\
x_{k+1} =& x_k - T A^\top \lambda_{k+1}. \notag
\end{align}
The $\lambda$-update equation is implicit in the desired update variable $\lambda_{k+1}$,
and to obtain it explicitly involves the proximal mapping, also known as the resolvent
of $\Sigma \partial F^*$, see \ref{app:fixedpoint}.

The desired proximal mapping is expressed as the argument of a minimization problem,
which can be derived in a few steps.
Regarding the right-hand side of the $\lambda$-update as an argument to the desired
mapping
\begin{equation*}
\lambda_\text{arg} = \lambda_k + \Sigma A \bar{x}_{k+1},
\end{equation*}
this update equation becomes
\begin{equation}
\label{lupdate}
\Sigma \partial F^*(\lambda_{k+1}) + \lambda_{k+1} = \lambda_\text{arg}.
\end{equation}
Dropping the $k+1$ subscript from $\lambda_{k+1}$ for clarity,
rearranging terms, and multiplying through by $\Sigma^{-1}$ yields
\begin{equation*}
\partial F^*(\lambda) + \Sigma^{-1}( \lambda - \lambda_\text{arg}) = 0,
\end{equation*}
where inversion of $\Sigma$ is possible because it is assumed to be a symmetric positive
definite matrix.
The left-hand side can be written as a total differential
\begin{equation}
\label{proxgrad}
\frac{\partial}{\partial \lambda} \left( F^*(\lambda) + \frac{1}{2}
(\lambda - \lambda_\text{arg})^\top
\Sigma^{-1}( \lambda - \lambda_\text{arg}) \right) = 0.
\end{equation}
The second term inside the differential operation is a quadratic distance
function with metric $\Sigma^{-1}$
\begin{equation*}
\| \lambda - \lambda_\text{arg}\|^2_{\Sigma^{-1}} \equiv
(\lambda - \lambda_\text{arg})^\top
\Sigma^{-1}( \lambda - \lambda_\text{arg}),
\end{equation*}
because $\Sigma^{-1}$ is also a symmetric positive definite matrix.
Both the distance function and $F^*$ are convex functions, so their sum
is also a convex function. Accordingly, Eq. (\ref{proxgrad}) can
be viewed as specifying the minimizer of the function inside the differentiation;
setting the gradient of a convex function to zero identifies its minimizer.
Thus, the explicit expression for $\lambda_{k+1}$ from Eq. (\ref{lupdate})
is
\begin{equation}
\label{proxargmin}
\lambda_{k+1} = \argmin_\lambda \left\{ F^*(\lambda) + \frac{1}{2}
\| \lambda - \lambda_\text{arg}\|^2_{\Sigma^{-1}} \right\}.
\end{equation}
Technically, Eq. (\ref{proxgrad}) only applies when $F^*$ is smooth, but
Eq. (\ref{proxargmin}) does apply for the desired case where $F^*$ is convex
and possibly non-smooth.
For general $\Sigma$, the problem in Eq. (\ref{proxargmin}) can be challenging to solve
analytically. Thus, for the CPPD framework considered here, $\Sigma$ is restricted
to a diagonal matrix with all diagonal values set to the scalar $\sigma$, i.e.
\begin{equation*}
\Sigma = \sigma I,
\end{equation*}
where $I$ is the identity matrix.
Equation~(\ref{proxargmin}) becomes
\begin{equation*}
\lambda_{k+1} = \argmin_\lambda \left\{ F^*(\lambda) + \frac{1}{2 \sigma}
\| \lambda - \lambda_\text{arg}\|^2_{2} \right\} \equiv \prox_{\sigma F^*} (\lambda_\text{arg}),
\end{equation*}
where the argmin problem is defined as the prox mapping.

Using the explicit $\lambda$-update equation and reordering the steps
so that the $x$-update comes first, the CPPD update equations from Eq. (\ref{cppdint1})
become
\begin{align}
x^{(k+1)} =& x^{(k)} - T A^\top \lambda^{(k)}, \notag \\
\bar{x}^{(k+1)} =& 2 x^{(k+1)}  - x^{(k)}, \label{finalCPPD} \\
\lambda^{(k+1)} =& \prox_{\sigma F^*} (\lambda^{(k)} + \sigma \bar{x}^{(k+1)}).
\notag
\end{align}
In this algorithm framework, $T$ is still allowed to be a general positive definite
matrix. For the basic CPPD implementation as discussed in the original paper on the
primal-dual algorithm \cite{chambolle2011first}, $T$ is a diagonal matrix
\begin{equation*}
T= \tau I,
\end{equation*}
and the condition for convergence is
\begin{equation*}
\tau \sigma \le 1/\|A\|^2_2.
\end{equation*}
More generally, when $T$ is a positive definite matrix, the condition for convergence
is
\begin{equation*}
1/ \sigma \ge \|T A^\top A\|_2.
\end{equation*}
We demonstrate examples of both of these choices in Sec. \ref{sec:results}. The CPPD update
is proved to converge only in the case of strict inequality, but empirically, equality leads
to convergence.
A useful choice of step-sizes discussed by Pock and Chambolle \cite{Pock2011}
avoids the power method and allows for fast step-size computation, see \ref{app:steps}.
In the same Appendix, we also present a useful implementation of a non-diagonal step-matrix $T$
designed for image reconstruction problems in X-ray tomography with possibly non-standard scan
geometries.

\section{Results}
\label{sec:results}

In this section, instances of the CPPD algorithm are demonstrated on optimization problems
relevant to CT image reconstruction. The goals of the studies are to demonstrate 
use of the CPPD algorithm, to characterize the impact of algorithm parameters, to illustrate
convergence, and, in general, to motivate the necessity of comprehensive empirical studies.
All presented studies focus on sampling sufficiency, and the simulation studies model breast CT.
Various scanning geometries are considered and all of the simulations take the form of
an ``inverse crime'' \cite{wirgin04} study, where the projection data are perfectly
consistent with object and data models. 

The inverse crime studies have multiple purposes. First, they present a stringent test
for correct algorithm implementation, because we know that the data discrepancy can be
driven to zero. Second, this is the only practical way to demonstrate sampling sufficiency
for the discrete-to-discrete data model used in iterative image reconstruction; if the
sampling is sufficient the converged solution will be the test phantom,
and if not, the converged solution will differ from the test phantom.  Third, the inverse
crime set up is excellent for empirical demonstration of algorithm convergence rate, which
is useful for parameter tuning and optimization solver comparison.

\subsection{Optimization problems and solvers}

The optimization problems considered are designed to
invert the discrete-to-discrete data-model.
\begin{equation*}
Xf=g,
\end{equation*}
where $X$ represents discrete-to-discrete projection; $f$ is the vector of image pixels;
and $g$ is a vector of projection data. In implementing $X$ this matrix is the product
\begin{equation}
\label{activepixels}
X = X_\text{grid} M_\text{FOV}, 
\end{equation}
where $M_\text{FOV}$ is a diagonal matrix with entries 1 and 0 on the diagonal that
correspond respectively to pixels inside and outside of the FOV; $X_\text{grid}$
represents projection of the whole rectangular image grid. In the following
studies, $g$ is taken to be noiseless projection of the discretized digital phantom using the
system matrix $X$.

\subsection*{Least-squares}
Two optimization problems are considered for inversion of the data model.
The first is formulated as a least-squares (LSQ) optimization
\begin{equation*}
f^\star = \argmax_f \frac{1}{2} \|X f -g \|^2_2,
\end{equation*}
where the solution $f^\star$ is arrived at only through the available projection
data. In inverse crime studies, if $f^\star$ matches the test phantom, this
is evidence that the matrix $X$ encodes a CT system that is sufficiently sampled.
Inversion by LSQ provides the opportunity to demonstrate the CPPD algorithm on a ubiquitous
optmization problem that is amenible to solution by many standard algorithms.

\begin{algorithm}
\hrulefill
\begin{algorithmic}[1]
\State $f^{(k+1)} \gets f^{(k)} - \tau  X^\top \lambda^{(k)}$
\State $\bar{f} \gets 2f^{(k+1)} -f^{(k)}$
\State $\lambda^{(k+1)} \gets \left( \lambda^{(k)} + \sigma(X\bar{f} - g)\right)/ (1+\sigma)$
\end{algorithmic}
\hrulefill
\caption{Pseudocode for CPPD-LSQ inner loop for iteration
$k$. See \ref{app:instances} for derivation of the CPPD-LSQ updates.
The step length parameters are set to $\sigma = \rho/L$ and $\tau = 1/(\rho L)$, where
$L=\|X\|_2$. The step-ratio parameter $\rho$ and iteration number $k_{\rm max}$ are varied in the
studies. Although warm-starting is possible, we do not consider this here and initialize
$f^{(0)} =0$ and $\lambda^{(0)} = 0$.}
\label{alg:CPPD-LSQ}
\end{algorithm}

\begin{algorithm}
\hrulefill
\begin{algorithmic}[1]
\State $f^{(k+1)} \gets f^{(k)} - (\alpha/L^2)  X^\top (Xf^{(k)} -g)$
\end{algorithmic}
\hrulefill
\caption{Pseudocode for GD-LSQ inner loop for iteration
$k$. The step length is determined by $\alpha \in (0,2)$ and $L$, where
$L=\|X\|_2$. The parameters $\alpha$ and iteration number $k_{\rm max}$ are varied in the
studies. We initialize $f^{(0)} =0$.}
\label{alg:GD-LSQ}
\end{algorithm}

The derivation of the CPPD algoritm for least-squares optimization (CPPD-LSQ) is presented
in \ref{app:instances}, and the pseudo-code is listed in Algorithm \ref{alg:CPPD-LSQ}.
We compare performance of CPPD-LSQ
with gradient descent applied to the LSQ problem (GD-LSQ), as listed in \ref{alg:GD-LSQ},
and to conjugate gradients least-squares, as listed on page 57 of Ref. \cite{paige1982lsqr}.
The parameters of GD-LSQ and CPPD-LSQ are explained in their respective pseudo-codes. The only
parameter varied for CGLS is the total number of iterations; the starting image is initialized
to zero and no attempt at pre-conditioning is made. CGLS is chosen because it is the gold-standard
solver for large-scale LSQ problems. Basic GD-LSQ, with a fixed step-size,
is almost never used in practice
because there are many options that are more efficient,
see for example Ref. \cite{jensen2012implementation}, but it provides a common reference point.

\begin{figure}[t]
\centering
\includegraphics[width=0.7\textwidth]{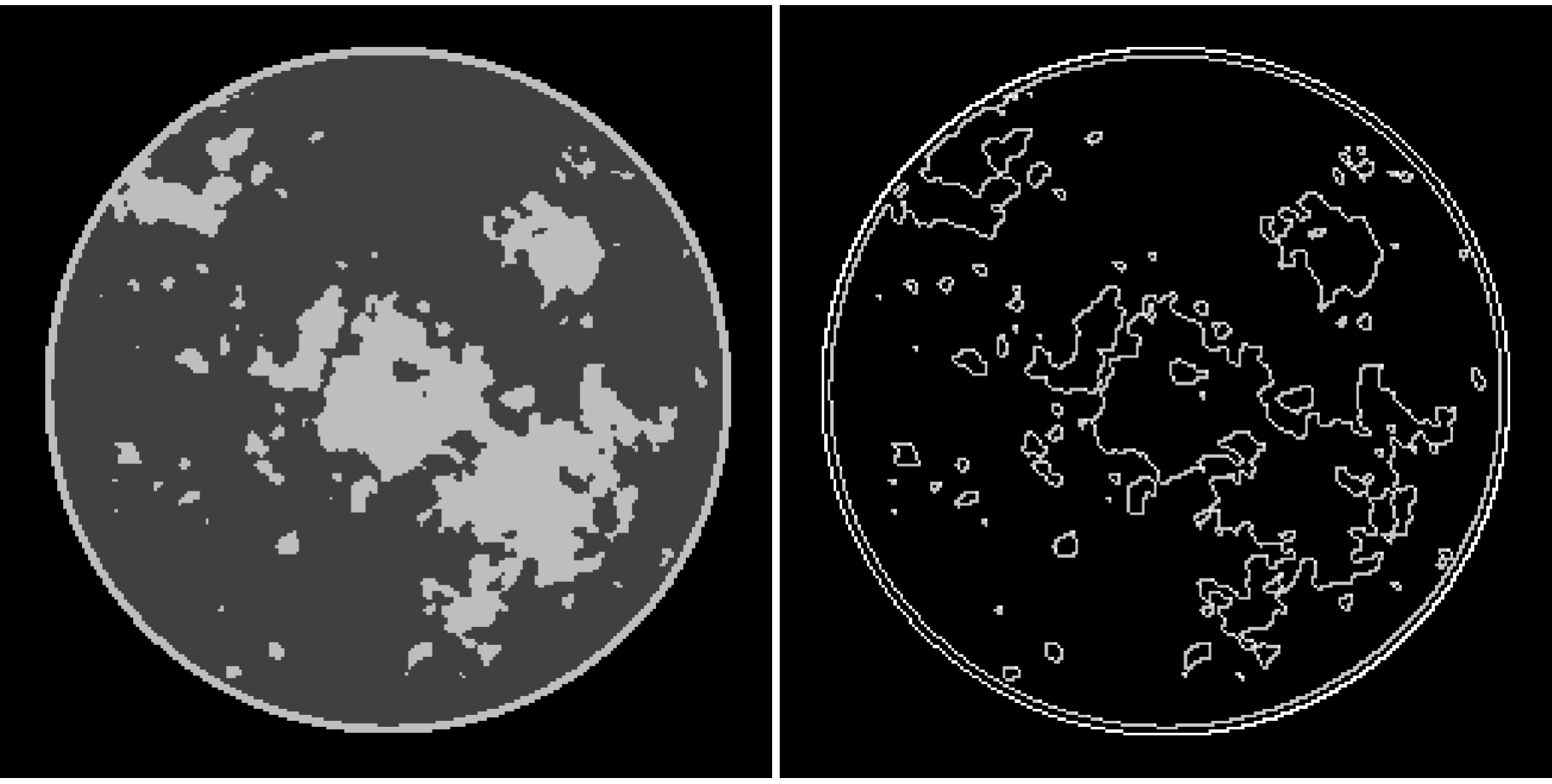}
\caption{\label{fig:bphantom} Computerized breast phantom and its gradient magnitude image (GMI).
In the GMI it is apparent that the phantom has a high degree of gradient sparsity even though
the tissue borders are highly irregular. The gradient sparsity is useful for testing ideal recovery
from under-sampled data by use of sparsity regularization with TV. The attenuation values in this
phantom or taken to be 0.194 cm$^{-1}$ and 0.233 cm$^{-1}$ for fat and fibro-glandular tissue, respectively.
The phantom and its GMI are shown in a gray scale window of [0.174, 0.253] cm$^{-1}$ and [0.0, 0.05] cm$^{-1}$,
respectively.}
\end{figure}

\subsection*{Total-variation constrained least-squares}
The second optimization problem considered is the
total-variation (TV) constrained least-squares (TVCLSQ) problem
\begin{equation*}
f^\star = \argmax_f \frac{1}{2} \|X f -g \|^2_2 \text{ such that } \|Df\|_1 \le \gamma,
\end{equation*}
where $D$ is the discrete gradient matrix operator and $\gamma$ is the constraint parameter.
For inverse crime studies, the phantom is known and $\gamma$ can be set to the phantom TV value
$\gamma_\text{ph}$. This optimization problem is useful for
Compressive Sensing \cite{candes2006robust,donoho2006compressed} investigations, where gradient
sparsity is exploited to reduce the scanning effort. In compressive sensing, image recovery
is governed by the number of samples, the size of the data, and the sparsity of the object.
If the former is sufficiently greater than the latter, the image can be recovered even for
configurations with reduced sampling.
For CT, reduced sampling usually means reduced numbers of projections or scanning angular range.
A theoretical, empirical compressive sensing study for CT
was carried out using this optimization problem in Ref. \cite{jorgensen2015little}.
The CPPD-TVCLSQ pseudo-code is derived in \ref{app:instances} where the
listing appears in Algorithm \ref{alg:CPPD-TVCLSQ}.
Compressive sensing has motivated a number of sparsity-exploiting image reconstruction
studies in medical imaging \cite{lustig2007sparse,sidky2008image,graff2015}.

\begin{figure}[t]
\centering
\includegraphics[width=0.6\textwidth]{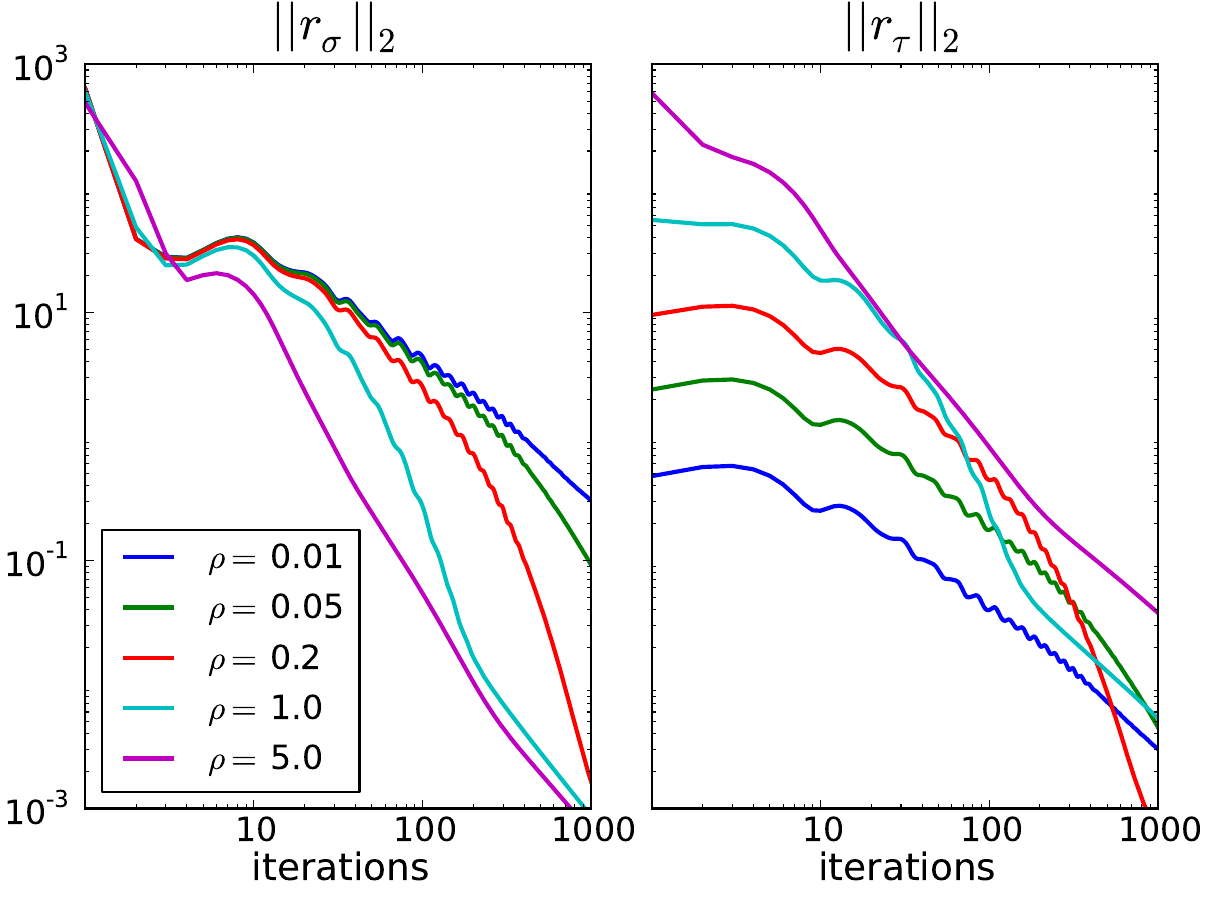}
\caption{\label{fig:lsq_sig_tau} Plots of $\|r_\sigma\|_2$ and $\|r_\tau\|_2$,
see Eqs.~(\ref{cc1}) and (\ref{cc2}) in \ref{app:convergence}, as a function of iteration number
and for different values of the step-size ratio parameter $\rho$.
}
\end{figure}

\subsection{simulation parameters}
The test set-up is based on breast CT in a 2D setting.
The phantom is shown in Fig. \ref{fig:bphantom} and
it is generated by a phantom model described in Ref. \cite{Reiser10}. The phantom and image reconstruction
grids are both 256$\times$256 pixels covering an area 18 cm $\times$ 18 cm.
The simulated CT sampling configuration is circular, fan-beam with varying numbers of projections
and scanning arc length.
The simulated projections employ a linear detector of 512 bins; the source-to-center
distance is 36 cm, and source-to-detector is 72 cm. The detector length is set so that the field-of-view
(FOV) has a diameter of 18 cm and matches the inscribed circle of the pixel grid.

\begin{figure}[t]
\centering
\includegraphics[width=0.6\textwidth]{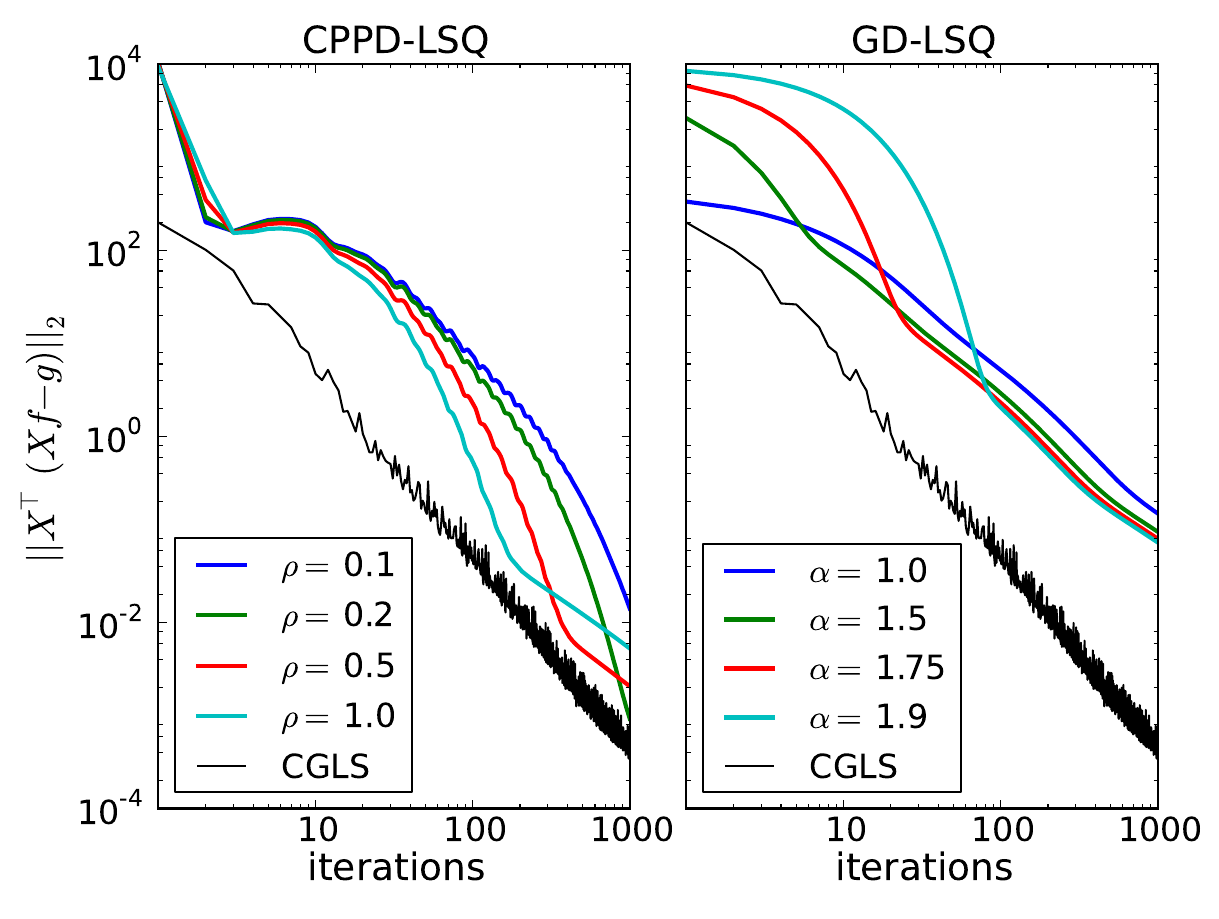}
\caption{\label{fig:lsq_grads} Plots of the LSQ objective function gradient magnitude
for CPPD-LSQ, GD-LSQ, and CGLS as a function of iteration number.
}
\end{figure}

The active pixels are those inside the FOV, and their 
number is 51,468 or approximately 79\% of the 256$\times$256 pixels in the square array.
The reason that only FOV pixels are selected for reconstruction is to make sure that all
active pixels are visible in all projections. Implementation-wise, a masking operation,
which multiplies all non-FOV pixels by zero, is applied prior to projection and after
back-projection, see Eq. (\ref{activepixels}).

In the presented empirical studies, three scan configurations are considered: full-sampling,
128 projections over a $2\pi$ scanning arc; sparse-view, 32 projections over a $2\pi$ scanning
arc; and limited angular-range, 128 projections over a $3 \pi/4$ scanning arc.
For both the full-sampling and limited angular-range scanning configurations the number
of samples is 128$\times$512, which is equal to the pixel grid size of 256$\times$256.
From the perspective of number of samples versus number of unknowns,
these configurations represent slight over-sampling, because only 79\% of the pixels are
active due to the FOV masking. This amount of over-sampling, by itself, presents a challenge.
We have found that the X-ray transform matrix
conditioning is relatively poor \cite{Jakob13}, when the number of unknown pixel values
is similar to the number of samples, irrespective of scanning arc length. 
The limited angular-range configuration provides an interesting test, because it is
sufficiently sampled from the perspective of number of unknowns and samples but the
corresponding continuous X-ray transform model is known to be insufficiently sampled
since the scanning arc length is less than $\pi$ plus fan-angle. The sparse-view
scan configuration is clearly under-sampled if no prior information on the object is
exploited, because the number of samples is less than the number of unknown pixel
values. For this configuration and phantom, exploiting gradient sparsity is known
to be an effective strategy. These sampling conditions are chosen to illustrate
various aspects of CPPD convergence.

\begin{figure}[h!]
\centering
\includegraphics[width=0.6\textwidth]{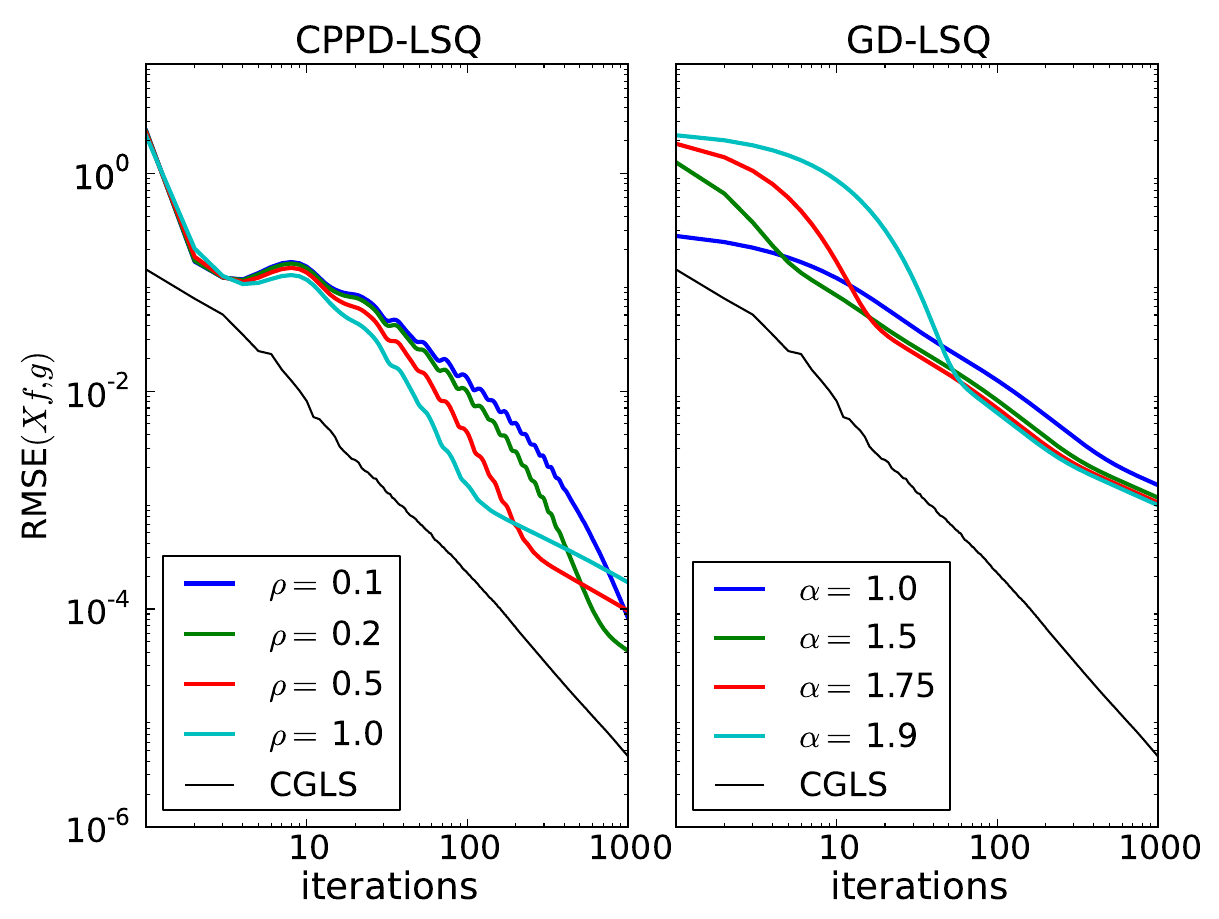}
\caption{\label{fig:lsq_derrs} Plots of the LSQ objective function
for CPPD-LSQ, GD-LSQ, and CGLS as a function of iteration number.
}
\end{figure}

\subsection{Least-squares convergence studies for full-sampling}
\label{sec:LSQconvergence}

The convergence plots in Fig. \ref{fig:lsq_sig_tau} show the splitting gap and transversality
evolution for CPPD-LSQ as a function of iteration number for different values
of the step-size ratio $\rho$, where
\begin{equation*}
\rho = \sqrt{\frac{\sigma}{\tau}}.
\end{equation*}
The transversality, $r_\tau$, and splitting gap, $r_\sigma$, convergence metrics for CPPD
are presented in detail in \ref{app:convergence}.  When $\rho$ is larger
than 1, $\sigma >\tau$, and as a general overall trend larger $\rho$ correlates with a faster
decrease in $\|r_\sigma\|_2$ while smaller $\rho$ tends to improve convergence in $\|r_\tau\|_2$.
Overall the curves from the shown intermediate values of $\rho=0.2$ and $\rho=1.0$
seem most promising when considering both convergence plots together.

\begin{figure}[t]
\centering
\includegraphics[width=0.6\textwidth]{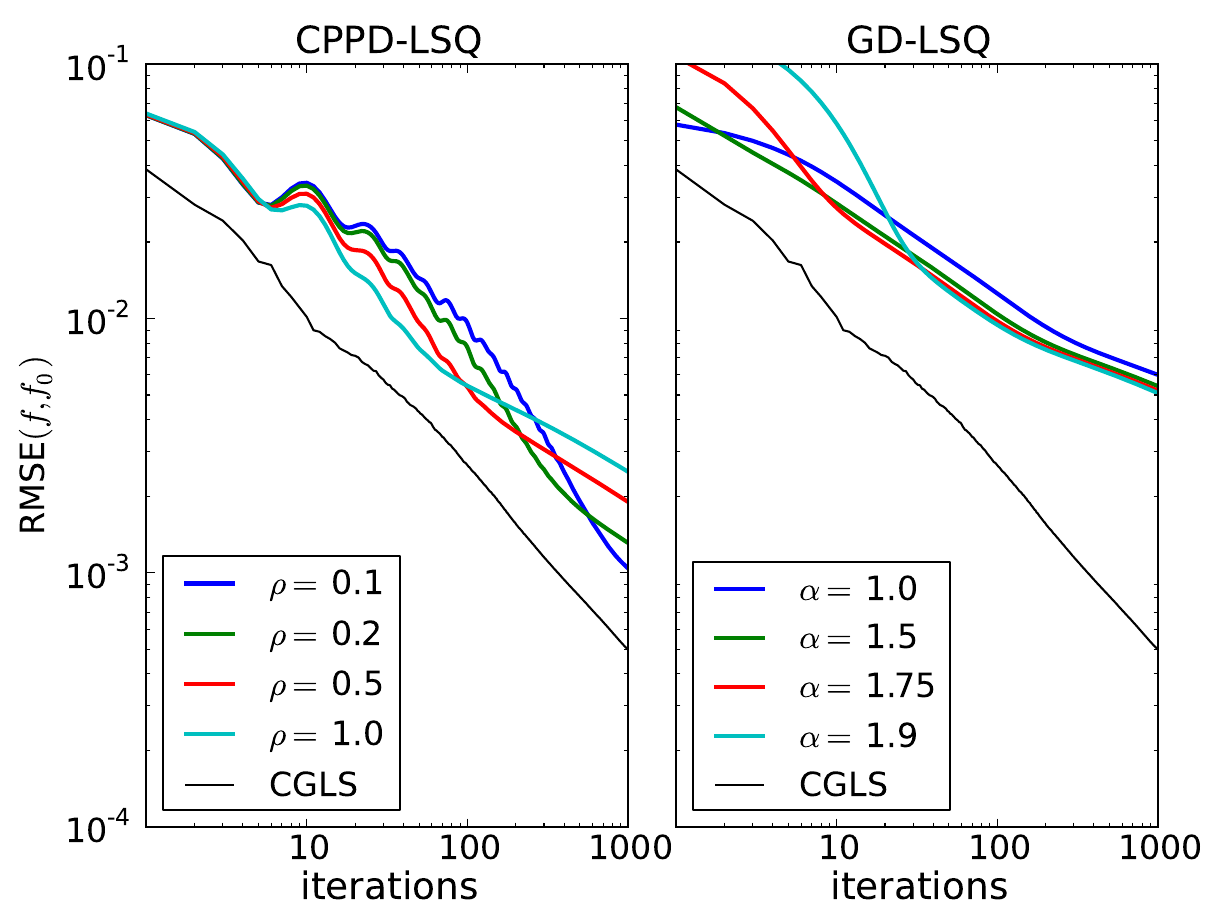}
\caption{\label{fig:lsq_ierrs} Plots of the image root-mean-square-error (RMSE)
discrepancy
for CPPD-LSQ, GD-LSQ, and CGLS as a function of iteration number.
}
\end{figure}

For the LSQ problem, which has a differentiable objective function, we can
also plot the magnitude of the LSQ objective gradient as a function of iteration number.
This metric is a complete first-order convergence condition and it allows comparison with GD-LSQ and CGLS.
The gradient of the LSQ objective function
\begin{equation*}
\phi(f) = \frac{1}{2}\|Xf-g\|_2
\end{equation*}
is 
\begin{equation*}
\partial \phi(f) = X^\top (Xf -g).
\end{equation*}
The comparison of CPPD-LSQ with GD-LSQ and CGLS is shown in Fig. \ref{fig:lsq_grads}.
By this metric and for this particular CT configuration, the convergence rate
of CPPD-LSQ is better than GD-LSQ and not as good as CGLS.

\begin{figure}[t]
\centering
\includegraphics[width=0.8\textwidth]{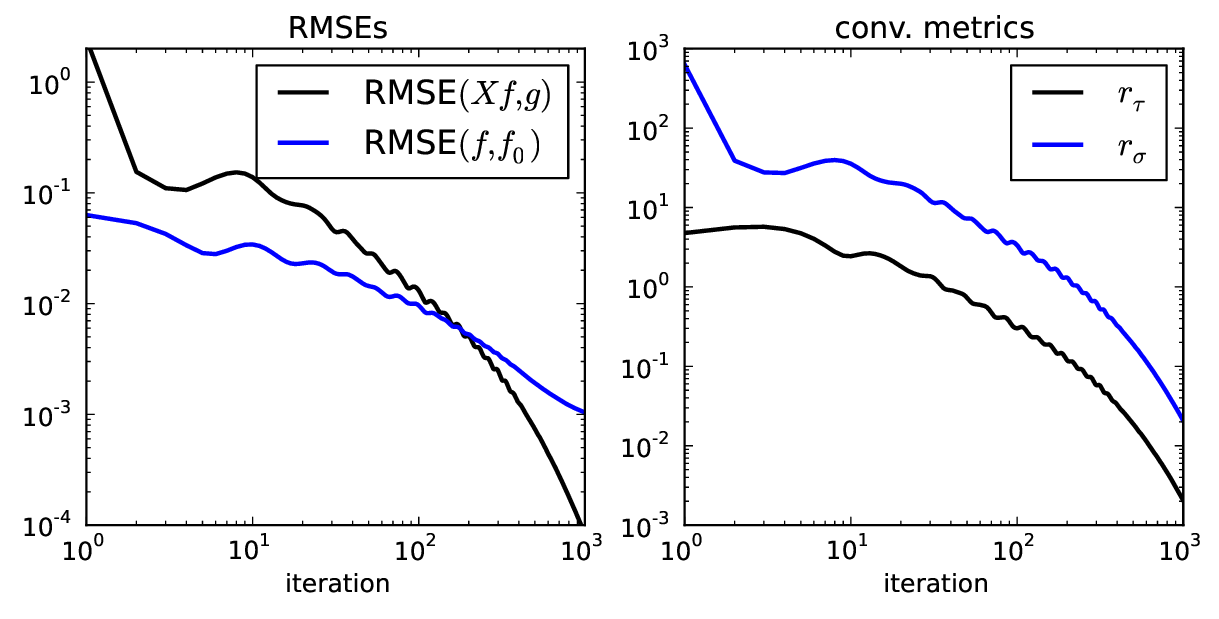}
\caption{\label{fig:lsq_conv} Convergence metrics for CPPD-LSQ for a step-size
ratio of $\rho=0.1$. The left panel shows the image and data RMSE, and the right
panel displays the transversality condition and splitting gap.
}
\end{figure}

Because we are performing an inverse crime study there are also two more convergence metrics
available. The LSQ objective value itself can be driven to zero, and the corresponding
results are shown in Fig. \ref{fig:lsq_derrs}.
The LSQ objective function should tend to zero no matter what is the data
sampling configuration for consistent data. To test sampling sufficiency,
it is necessary to demonstrate convergence of the image estimate to the test phantom.
The plots of the image estimate discrepancy 
in Fig. \ref{fig:lsq_ierrs} show
curves that all tend to zero, thus the present configuration
does provide sampling sufficiency.

\begin{figure}[h!]
\centering
\includegraphics[width=0.5\textwidth]{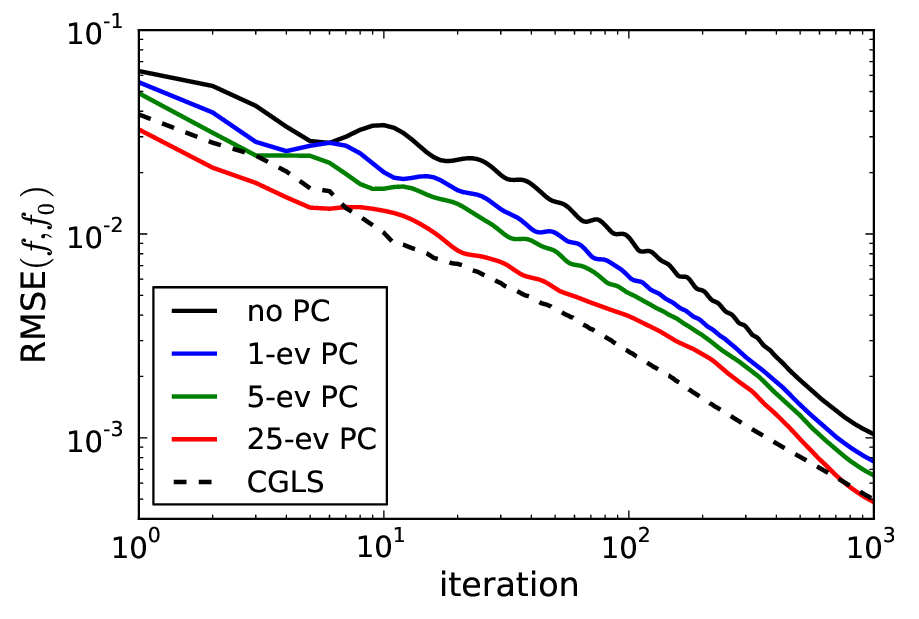}
\caption{\label{fig:lsq_pc} Impact of non-diagonal preconditioning (PC) on image RMSE
for CPPD-LSQ applied to image reconstruction for the 128-view data configuration.
The legend indicates the number of eigenvectors (ev) of $X$ used in forming the preconditioner.
Also shown are the traces for scalar $\sigma$ and $\tau$, i.e. no preconditioning and
generic CGLS.
}
\end{figure}

The convergence curves are all plotted on a log-log scale because all of these quantites
are expected to tend to zero. 
By using a log-log plot, it is possible to visually estimate at which iteration number
a desired accuracy can be obtained. Asymptotically, first-order algorithms converge
to the solution by a power of the iteration number; a quantity that tends to zero will
thus appear linear on a log-log scale when this asymptotic behavior is reached.
For the LSQ problem, the metrics $r_\sigma$, $r_\tau$, and object gradient magnitude
should all tend to zero irrespective of data consistency or sampling. The LSQ objective
function tends to zero with perfectly consistent data, and the image RMSE
tends to zero with consistent data and sufficient sampling.

Taking in all these curves, there is the general impression that CPPD-LSQ performance in terms
of convergence lies somewhere between GD-LSQ and CGLS. The range of interesting $\rho$ settings
for CPPD-LSQ appears to be between 0.1 and 1.0. Focusing on the convergence metrics for CPPD-LSQ,
they clearly all tend to zero. The $\rho$-rank-order of the convergence is seen to depend on iteration
number and which metric is used.
Because the presented results center on inverse crime studies,
the image RMSE metric is most relevant and $\rho$ is selected based on rank-order at 1000 iterations,
which is the number of iterations used in the simulations. For the $\rho$ values tested,
$\rho=0.1$ leads to the lowest image RMSE after 1,000 iterations and the corresponding
convergence metric traces are displayed in Fig. \ref{fig:lsq_conv}.
Whichever convergence metric is used
for selecting $\rho$, it is clear that the convergence rate does depend strongly on $\rho$ and
this parameter must be tuned in order to maximize CPPD performance.

\begin{figure}[t]
\centering
\includegraphics[width=0.8\textwidth]{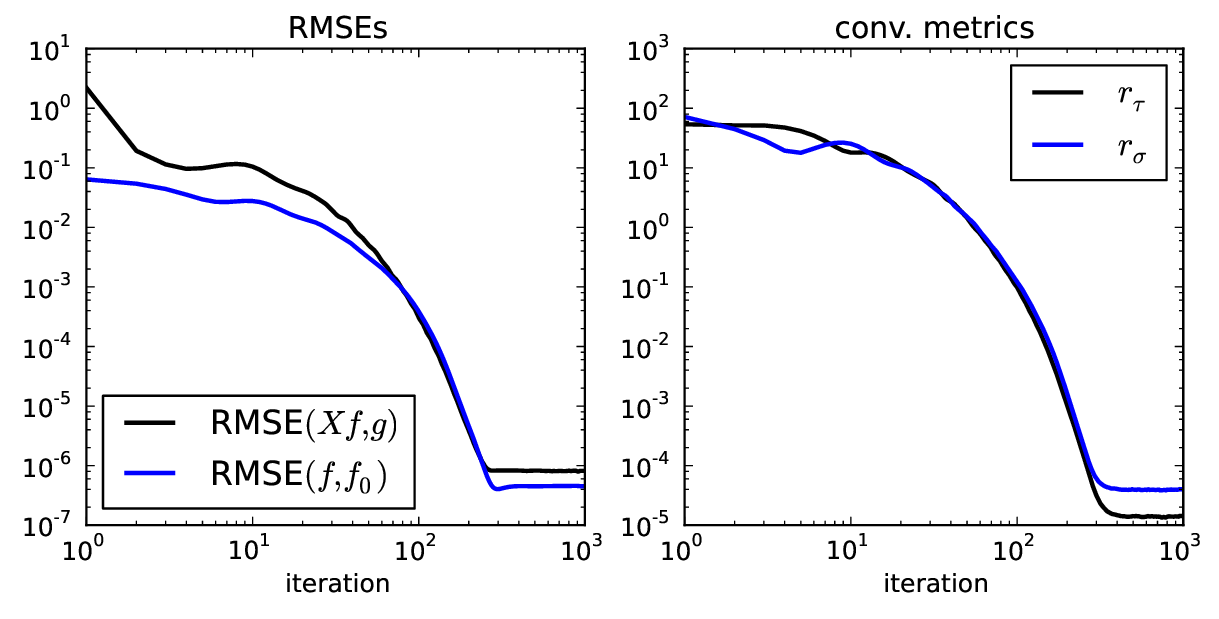}
\caption{\label{fig:tvlsq_conv} Convergence metrics for CPPD-TVCLSQ for a step-size
ratio of $\rho=1.0$. The left panel shows the image and data RMSE, and the right
panel displays the transversality condition and splitting gap.
}
\end{figure}
\begin{figure}[h!]
\centering
\includegraphics[width=0.5\textwidth]{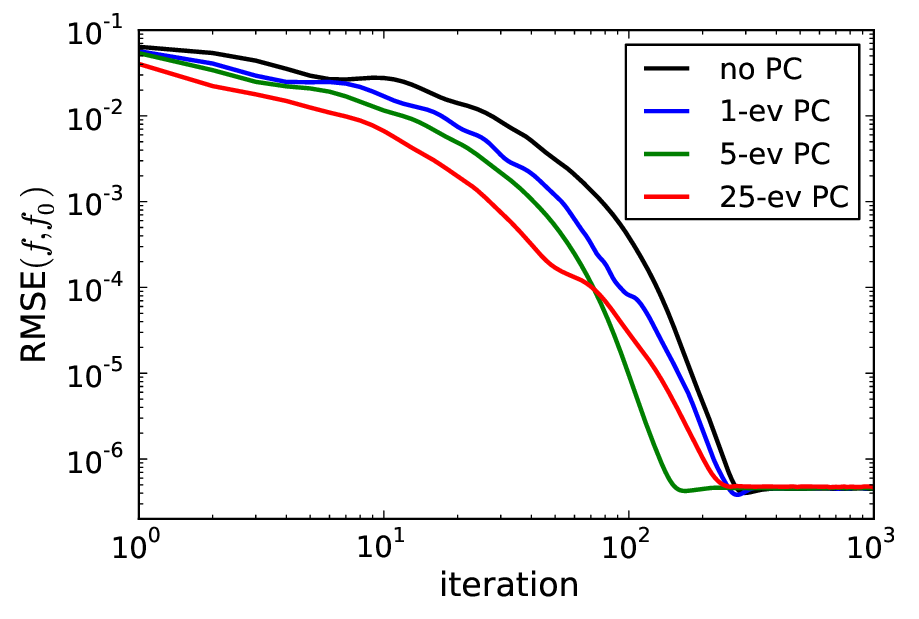}
\caption{\label{fig:tvlsq_pc} Impact of non-diagonal preconditioning (PC) on image RMSE
for CPPD-TVCLSQ applied to image reconstruction for the 128-view data configuration.
The legend indicates the number of eigenvectors (ev) of $X$ used in forming the preconditioner.
Also shown are the traces for scalar $\sigma$ and $\tau$, i.e. no preconditioning and
generic CGLS.
}
\end{figure}
\begin{figure}[h!]
\centering
no preconditioning
\includegraphics[width=1.0\textwidth]{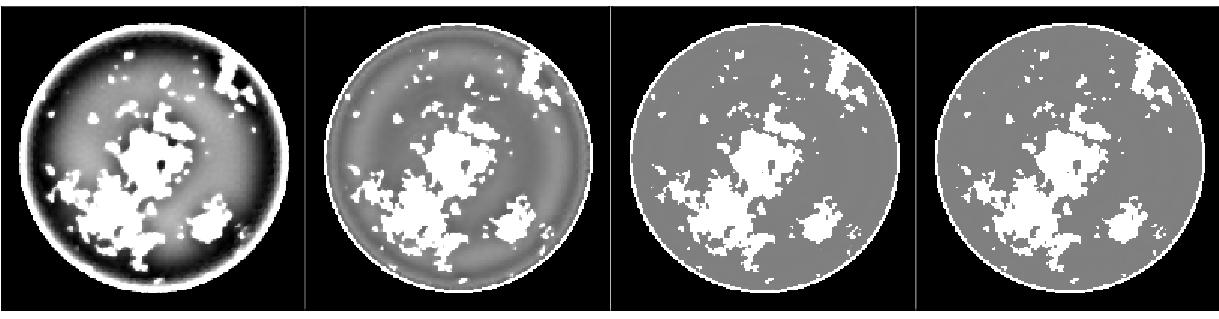}
non-diagonal preconditioning with five eigenvectors
\includegraphics[width=1.0\textwidth]{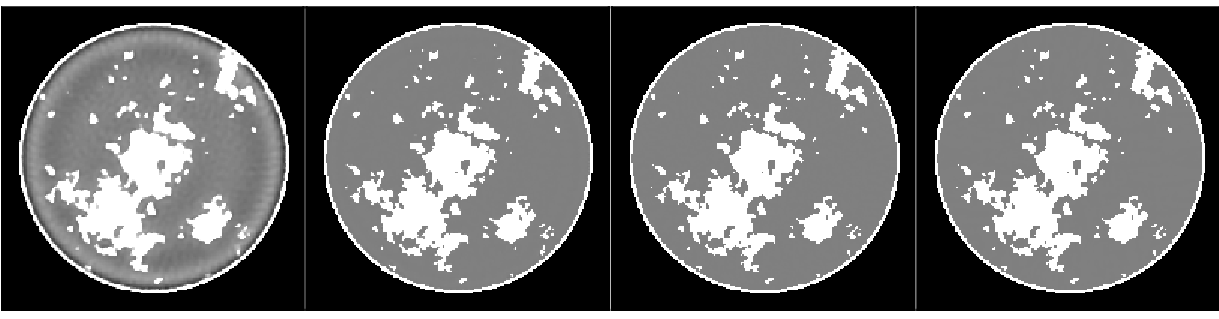}
iterations:~20~~~~~~~~~~~~~~~~~~50~~~~~~~~~~~~~~~~~~~~~~~100~~~~~~~~~~~~~~~~~~~~~~~~200~~~~~~~~
\caption{\label{fig:tvlsq_seq} Sequence of image estimates for different iteration
numbers of CPPD-TVCLSQ. Top row shows results for no preconditioning, and the bottom row displays
images for non-diagonal preconditioning. The gray scale is [0.174, 0.214]
cm$^{-1}$, which is centered on the background adipose
attenuation of 0.194 cm$^{-1}$ so that non-unformity in the background is easily seen.
}
\end{figure}

\subsection*{Improving CPPD-LSQ performance with non-diagonal preconditioning}

In case that greater efficiency is needed for the CPPD algorithm, non-diagonal step matrices
can be exploited as alluded to in Sec. \ref{sec:matrix-mapping}. To demonstrate such preconditioning,
the scalar step-size $\tau$ in the first line of Algorithm \ref{alg:CPPD-LSQ} is replace by a matrix $T$,
which, in the case of the LSQ problem, is set to an approximate inverse of $X^\top X$. For a
low-rank approximate inverse, see Eq. (\ref{nondiagtau}) in \ref{app:steps}, where $T$ is constructed
from the leading eigenvectors of $X^\top X$. The gain in convergence of the image RMSE metric
by use of the low-rank preconditioner is shown in Fig. \ref{fig:lsq_pc}. Interestingly, using only
the leading eigenvector by itself to form $T$ results in a significant drop in image RMSE.
Including additional eigenvectors does improve convergence, but the gain with each additional
eigenvector becomes smaller and smaller. For the shown example, it is possible to reach
the convergence rate of un-preconditioned CG by using 25 leading eigenvectors to form $T$.
The low-rank preconditioner is also shown to be useful for the results from
 other configurations considered here.

\subsection{Total-variation constrained least-squares convergence studies for full-sampling}
\label{sec:TVCLSQconvergence128-2pi}

In the following sub-sections, application of CPPD-TVCLSQ to reduced sampling configurations is
presented. Here, we apply CPPD-TVCLSQ to the same full sampling system that is studied previously
with CPPD-LSQ. For CPPD-TVCLSQ there are three algorithm parameters; as with CPPD-LSQ the iteration
number and step-size ratio $\rho$ need to be specified. Additionally, the TV constraint parameter,
$\gamma$, must also be set for TVCLSQ problem. The present results focus on image recovery
from ideal data, thus the constraint parameter is set to the phantom TV, i.e. $\gamma=\gamma_\text{ph}$.

As with CPPD-LSQ, 1000 iterations of CPPD-TVCLSQ are executed for several values of $\rho$.
Shown in Fig. \ref{fig:tvlsq_conv} are the curves corresponding to the $\rho$ value that
showed the most rapid convergence in the image RMSE.  Note that
the RMSE values reached are much lower than that of LSQ in Fig.~\ref{fig:lsq_conv}
in particular the image RMSE convergence is much more rapid for CPPD-TVCLSQ.
Both the RMSE curves and the CPPD convergence metrics $r_\sigma$ and $r_\tau$
show a rapid decreasing tend to zero, but then all curves hit a
hard plateau. The source of the plateauing behavior has been traced
to the finite precision of the computation, which is performed in single-precision
floating point arithmetic. This has been verified by changing the computer variable precision from 4-byte
to 8-byte floating point representation and by shrinking the overall size of the modelled
tomographic system. There are ways to increase the numerical
accuracy from fixed precision computation
\cite{kahan1965pracniques,Higham},
but we have not attempted this with the present CPPD-TVCLSQ implementation.

Applying the low-rank preconditioner in the form of a non-diagonal step matrix $T$ is
more complicated than the CPPD-LSQ case, because CPPD-TVCLSQ involves
the stacking of two matrices $X$ and $D$. In \ref{app:steps} one possible implementation
of this preconditioner is presented for CPPD applied the TV-penalized least-squares.
The same strategy also applies to CPPD-TVCLSQ, and the results of applying this preconditioner
is shown in Fig. \ref{fig:tvlsq_pc}. As with CPPD-LSQ, the low-rank preconditioner 
improves the convergence rates seen in the image RMSE curves.

For inverse crime studies the CPPD convergence metrics $r_\sigma$ and $r_\tau$ are not
as directly meaningful as the image and data RMSE, but it is important to note
the relative convergence between the RMSE curves and the CPPD convergence metrics

For inverse crime studies the CPPD convergence metrics $r_\sigma$ and $r_\tau$ are not
as directly meaningful as the image and data RMSE. It is, however, important to note
the relative convergence between the RMSE curves and the CPPD convergence metrics,
because the CPPD metrics will always tend to zero including the realistic conditions
where the data contains inconsistencies, while the RMSE curves will in general not do so.

\section{End of the main thread for now ...}
The next two unlabelled subsections have results for sparse-view and limited angular range scanning.
These results can be improved upon by optimizing the matrix stacking methodology. This topic will
be added in a future version of these notes.

\clearpage 

\subsection*{CPPD-TVCLSQ convergence studies, 32 projections over $2\pi$ scanning}
\label{sec:TVCLSQconvergence32-2pi}

\begin{figure}[h!]
\centering
\includegraphics[width=0.8\textwidth]{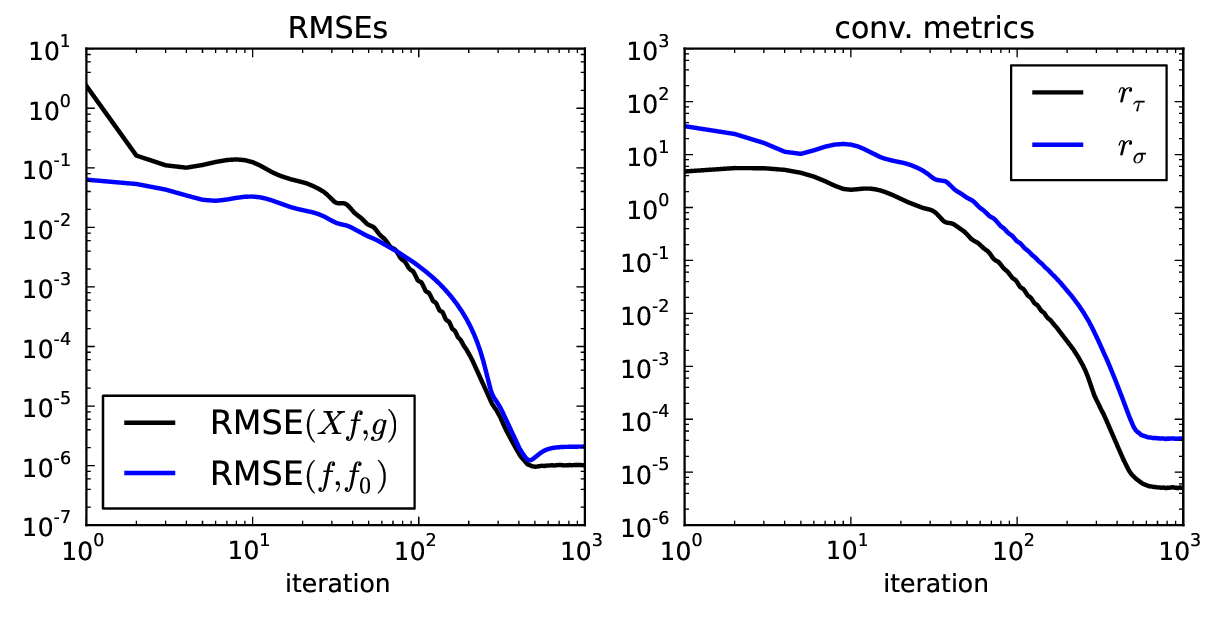}
\caption{\label{fig:tvlsq32_conv} Convergence metrics for CPPD-TVCLSQ for a step-size
ratio of $\rho=0.2$. The left panel shows the image and data RMSE, and the right
panel displays the transversality condition and splitting gap.
}
\end{figure}
\begin{figure}[h!]
\centering
\includegraphics[width=0.5\textwidth]{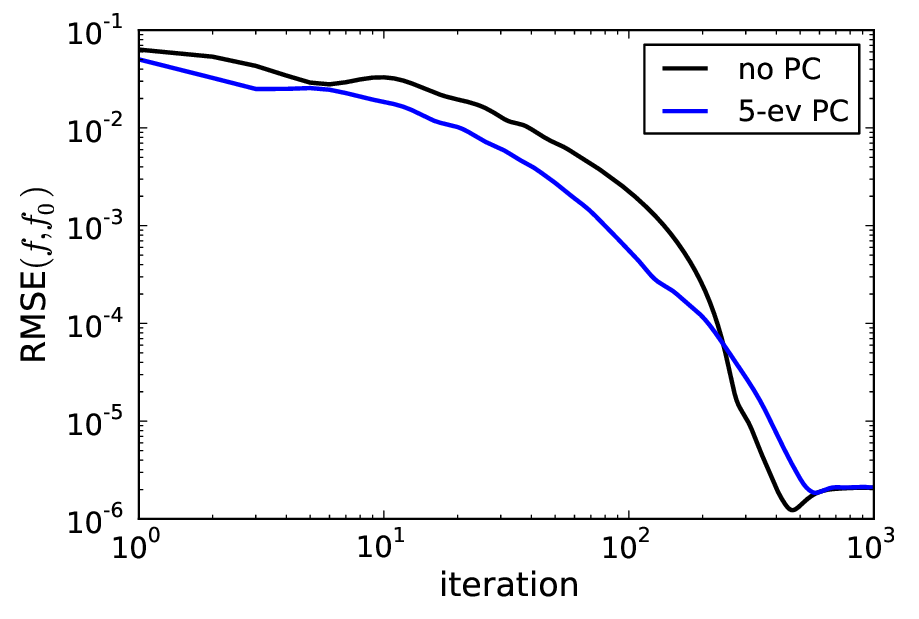}
\caption{\label{fig:tvlsq32_pc} Impact of non-diagonal preconditioning (PC) on image RMSE
for CPPD-TVCLSQ applied to image reconstruction for the 32-view data configuration.
The legend indicates the number of eigenvectors (ev) of $X$ used in forming the preconditioner.
Also shown are the traces for scalar $\sigma$ and $\tau$, i.e. no preconditioning and
generic CGLS.
}
\end{figure}
\begin{figure}[h!]
\centering
no preconditioning
\includegraphics[width=1.0\textwidth]{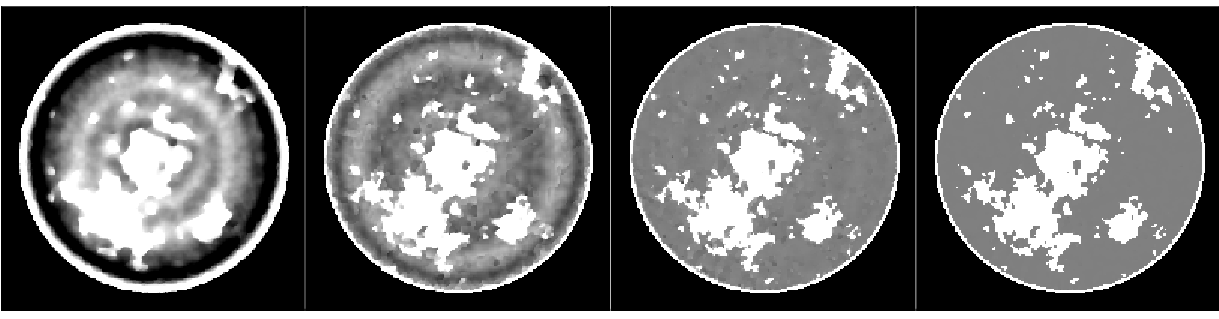}
non-diagonal preconditioning with five eigenvectors
\includegraphics[width=1.0\textwidth]{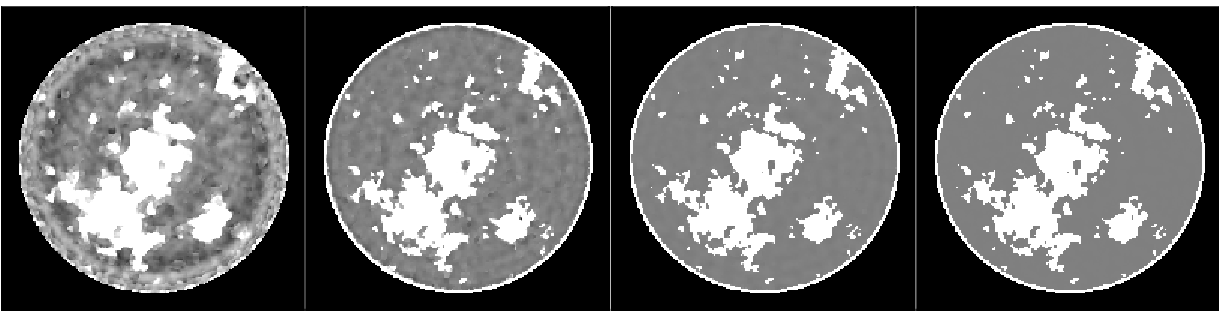}
iterations:~20~~~~~~~~~~~~~~~~~~50~~~~~~~~~~~~~~~~~~~~~~~100~~~~~~~~~~~~~~~~~~~~~~~~200~~~~~~~~
\caption{\label{fig:tvlsq32_seq} Sequence of image estimates for different iteration
numbers of CPPD-TVCLSQ. Top row shows results for no preconditioning, and the bottom row displays
images for non-diagonal preconditioning. The gray scale is [0.174, 0.214]
cm$^{-1}$, which is centered on the background adipose
attenuation of 0.194 cm$^{-1}$ so that non-unformity in the background is easily seen.
}
\end{figure}

\clearpage

\subsection*{CPPD-TVCLSQ convergence studies, 128 projections over $3\pi/4$ scanning}
\label{sec:TVCLSQconvergence128-3pi4}

\begin{figure}[h!]
\centering
\includegraphics[width=0.8\textwidth]{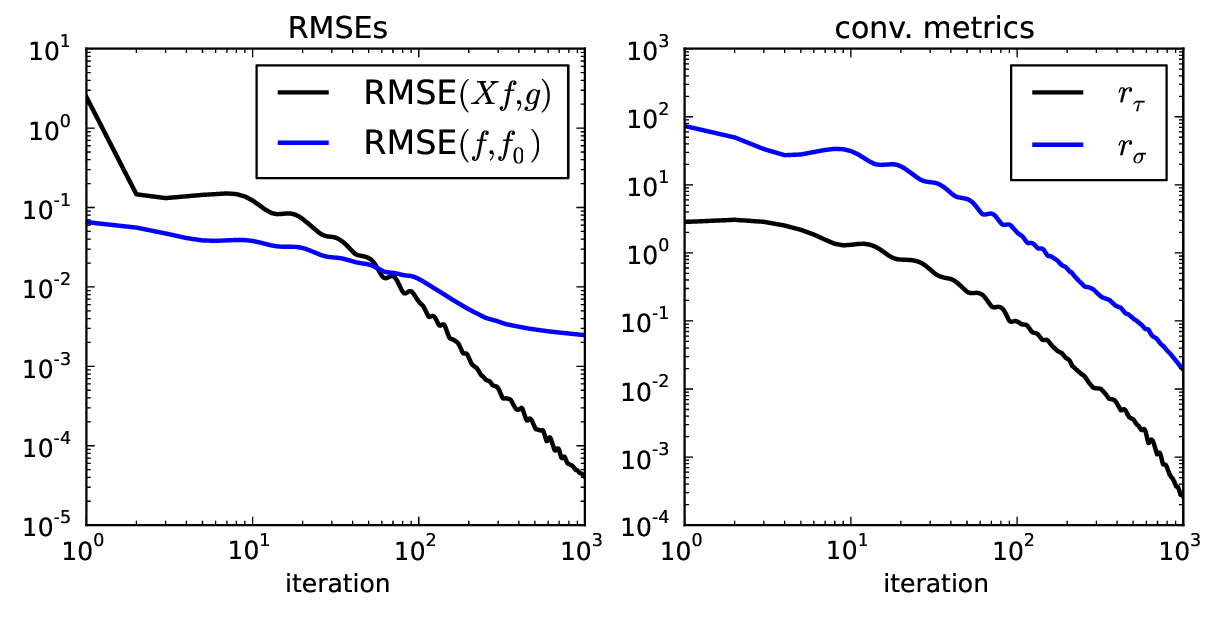}
\caption{\label{fig:tvlsq135deg_conv} Convergence metrics for CPPD-TVCLSQ for a step-size
ratio of $\rho=0.5$. The left panel shows the image and data RMSE, and the right
panel displays the transversality condition and splitting gap.
}
\end{figure}
\begin{figure}[h!]
\centering
\includegraphics[width=0.5\textwidth]{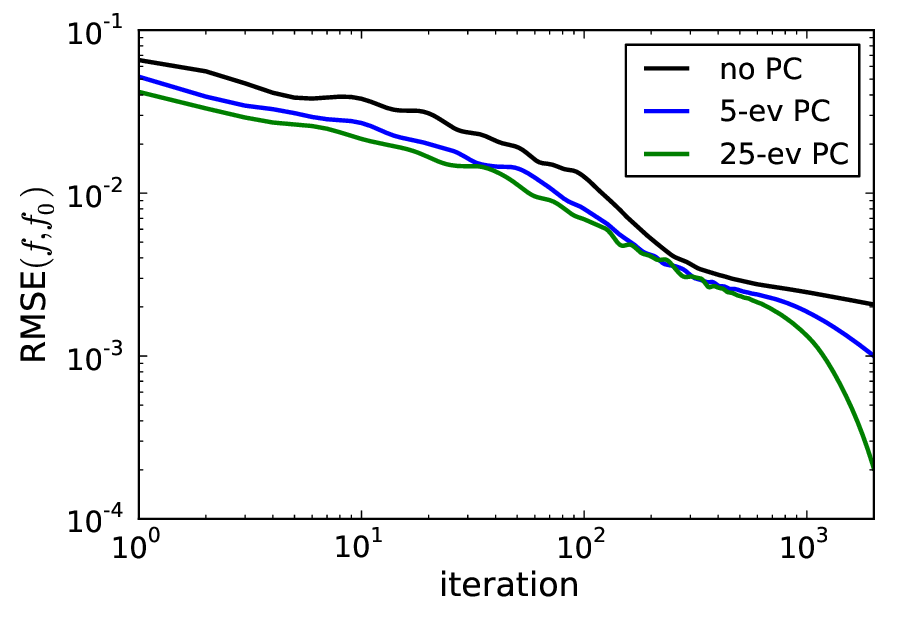}
\caption{\label{fig:tvlsq135deg_pc} Impact of non-diagonal preconditioning (PC) on image RMSE
for CPPD-TVCLSQ applied to image reconstruction for the 128-view, $3\pi/4$
angular range data configuration.
The legend indicates the number of eigenvectors (ev) of $X$ used in forming the preconditioner.
Also shown are the traces for scalar $\sigma$ and $\tau$, i.e. no preconditioning and
generic CGLS.
}
\end{figure}
\begin{figure}[h!]
\centering
no preconditioning
\includegraphics[width=1.0\textwidth]{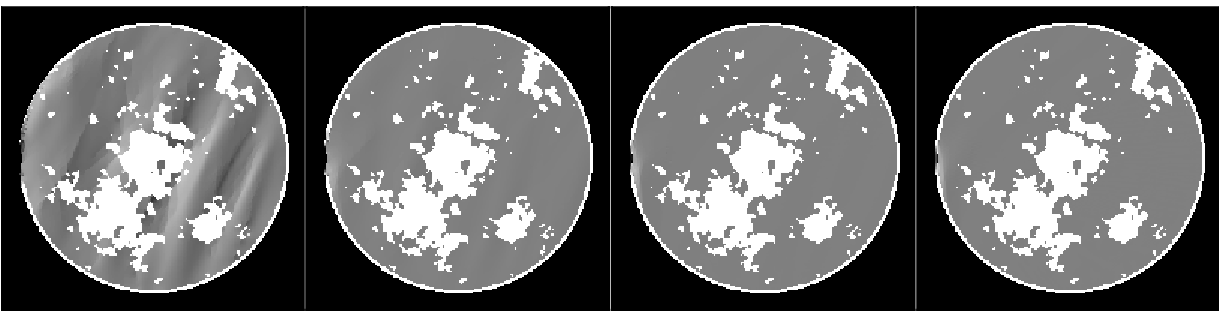}
non-diagonal preconditioning with five eigenvectors
\includegraphics[width=1.0\textwidth]{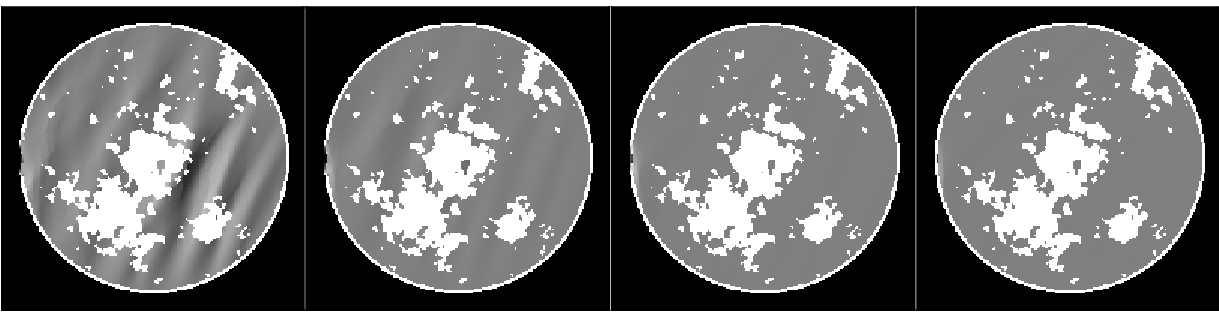}
non-diagonal preconditioning with 25 eigenvectors
\includegraphics[width=1.0\textwidth]{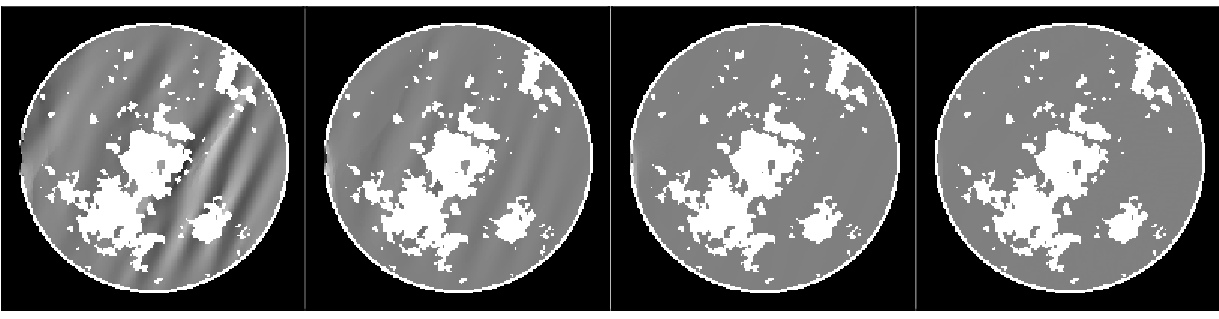}
iterations:~200~~~~~~~~~~~~~~~~~~500~~~~~~~~~~~~~~~~~~~~~~~1000~~~~~~~~~~~~~~~~~~~~~~~~2000~~~~~~~~
\caption{\label{fig:tvlsq135deg_seq} Sequence of image estimates for different iteration
numbers of CPPD-TVCLSQ. Top row shows results for no preconditioning; the middle row displays
images for 5-ev non-diagonal preconditioning; and the bottom row has the 25-ev results.
 The gray scale is [0.174, 0.214]
cm$^{-1}$, which is centered on the background adipose
attenuation of 0.194 cm$^{-1}$ so that non-unformity in the background is easily seen.
}
\end{figure}

\clearpage

\appendix

\section{Critical points and optimization}
\label{app:critical}

This appendix briefly summarizes the connection between critical points and optimization for smooth
functions. For a function $f(x)$ the critical points occur at
\begin{equation*}
\partial f(x_c) = 0.
\end{equation*}
Such points are related to extrema of $f(x)$, which can be specified by optimization: either minimization,
maximization, or some combination thereof.
All extrema of $f(x)$ are critical points, but not all critical points are extrema.
We ignore, here, the complicating issue of local and global extrema.
We need the connection between critical points and extrema for two reasons:
(1) to write down the solution of an optimization
problem as an equation; e.g. to be able to use
\begin{equation*}
\partial f(x^\star) = 0
\end{equation*}
as the solution of
\begin{equation*}
x^\star = \min_x f(x);
\end{equation*}
and (2) to go in the other direction, if we have a large-scale equation, where
the solution can be viewed as a critical point of a potential, it can be helpful
to write the problem as an optimization. The latter purpose is particularly
useful when trying to develop iterative algorithms to solve an equation, as opposed
to its direct solution.

To determine what type of extremum a critical point is or if it is an extremum at all,
it is necessary to examine higher order derivatives.
For example, for
\begin{equation*}
f(x) = x^2,
\end{equation*}
$\partial f(x=0) = 0$, and $\partial^2 f(x=0) = 2$. Because the second derivative is positive
we know that the critical point at $x_c=0$ is a minimum, and this critical point can be
specified by minimization
\begin{equation*}
x_c = \min_x x^2.
\end{equation*}
Clearly, for
\begin{equation*}
f(x) = -x^2,
\end{equation*}
$x_c=0$ is a maximum because the second derivative is negative, and its critical point can
be specified by maximization.
But then there are more ambiguous situations such as
\begin{equation*}
f(x) = x^3,
\end{equation*}
where the first and second derivatives are zero at $x_c=0$. In this case,
the lowest-order derivative, which is non-zero, is odd (namely the third-order derivative)
so $x_c=0$ is not an extremum and it cannot be specified by an optimization problem.
The above examples apply no matter what is the dimension of $x$.

There is a type of extremum that is possible only if $x$ is at least a 2-dimensional vector;
namely a critical point in $n$-dimensions, $n\ge2$, can be a saddle point - concave in some
directions and convex in others.
A clear example for 2-dimensional $x$ is
\begin{equation*}
f(x_1,x_2)=x_1^2 - x_2^2,
\end{equation*}
which is clearly convex as a function of $x_1$ and concave as a function of $x_2$.
Thus, this critical point can be specified by a combination of minimization and maximization
\begin{equation*}
\min_{x_1} \max_{x_2} \left\{x_1^2 - x_2^2\right\}.
\end{equation*}

A less clear example, however, is
\begin{equation*}
f(x_1,x_2)=x_1 x_2,
\end{equation*}
which is linear in $x_1$ and $x_2$. Linear functions can be taken as convex or concave.
To classify this critical point, we analyze it with second-order derivatives.
The Hessian provides second-order characterization of multi-dimensional functions
\begin{eqnarray*}
H= \partial^2 f(x)&=
\left( \begin{array}{cc} \partial^2 f/\partial x_1^2 &
\partial^2 f/\partial x_2 \partial x_1 \\
\partial^2 f/\partial x_1 \partial x_2  &
\partial^2 f/\partial x_2^2 \end{array} \right) \\
&=
\left( \begin{array}{cc} 0 &
1 \\
1  &
0 \end{array} \right).
\end{eqnarray*}
The critical point classification can then be determined by diagonalization of $H$.
The eigenvalues of $H$ yield the curvature of $f$ in the directions of the eigenvectors.
If there are negative and positive eigenvalues of $H$, the critical point is a saddle
point.
In this particular example, the eigenvalues are $1$ and $-1$ corresponding to the
eigenvectors $(1,1)^\top$ and $(1,-1)^\top$, respectively. Thus, we know that $f(x) = x_1 x_2$
is a saddle point. Changing coordinates as suggested by the diagonalization
\begin{eqnarray*}
s= x_1+x_2,\\
t= x_1-x_2,
\end{eqnarray*}
makes this abundantly clear
\begin{equation*}
x_1 x_2 = (s^2 - t^2)/4.
\end{equation*}
We know that optimization of $f(x)=x_1x_2$ to find the critical point involves min-max
because the critical point is a saddle point.
Interestingly, this critical point can be specified by either
of the following min-max problems
\begin{eqnarray*}
\min_{x_1}\max_{x_2} x_1x_2,\\
\min_{x_2}\max_{x_1} x_1x_2,
\end{eqnarray*}
where actually being able to carry out the calculation of these min-max problems requires
non-smooth analysis and this is discussed in \ref{app:nonsmooth}.

Analyzing the only the bilinear term in the saddle-point optimizations in Eqs. (\ref{saddle1}) and (\ref{saddle2})
the stationary saddle-point can be identified in two ways. For the bilinear  term in Eq. (\ref{saddle1}), if maximization 
is selected for the variables $x$ and $y$, minimization must be performed over $\lambda$. If minimization
is chosen for $x$ and $y$, maximization must be selected for $\lambda$. The reason why minimization over $x$ and $y$
is used when considering the complete Lagrangian expression in Eq. (\ref{saddle1}) is that $\phi(y)$ is convex.
Likewise, for Eq. (\ref{saddle2}), maximization over $\lambda$ is selected because $-\phi^*(\lambda)$ is concave.

\section{Non-smooth convex functions}
\label{app:nonsmooth}

The workings of the Chambolle-Pock primal-dual (CPPD) algorithm can be mostly understood within
the context of smooth optimization, but one of the main motivations for using the CPPD algorithm
is to perform optimization with non-smooth convex functions. Accordingly, we do need to cover
this topic, but we attempt to do so with a bare minimum of material. For more in depth presentation,
the reader is referred to the classic text by Rockafellar \cite{rockafellar1970convex}. There are
also a number of other textbooks, e.g. \cite{hiriart1993convex}, tutorial papers, and online reference material
on this topic.
The approach to the presentation here is greatly simplified by a comment made by Marc Teboulle,
co-author of FISTA \cite{Beck09}, at the 2014 SIAM Imaging Science conference in Hong Kong.
In the discussion after a presented paper on non-smooth optimization, Marc pointed out that in dealing
with non-smooth convex optimization the vast majority of cases where it is used
center on the absolute value and the indicator functions. The absolute value function
is well-known, but the indicator function might be less familiar to a medical physics audience.

Indicator functions are a convenient construct in convex analysis for converting a convex
set into a convex function
\begin{equation*}
\delta_C(x) = \left\{
\begin{array}{cc}
0 & x \in C \\
\infty & x \notin C
\end{array} \right. ,
\end{equation*}
where $C$ is a convex set. Clearly, the absolute value and indicator functions are not differentiable
everywhere. The concept of differentiation can be generalized to accommodate non-smoothness, and this
topic is taken up in \ref{app:subdiff}.

For optimization, the indicator function allows the restriction of possible solutions
to various convex constraints by adding infinite walls to the objective function.
In performing addition of convex functions,
we need an additional rule for handling infinity
\begin{equation}
\label{infadd}
a + \infty = \infty,
\end{equation}
where $a$ is any scalar, and for $a > 0$
\begin{equation}
\label{infmult}
a \cdot \infty = \infty,
\end{equation}
and
\begin{equation}
\label{infmultzero}
0 \cdot \infty = 0.
\end{equation}
Using the algebra of convex functions, constrained optimization can
be made to look like unconstrained optimization.

For example, the convex constrained minimization problem
\begin{equation*}
\min_x \frac{1}{2} x^2 \; \; {\rm such\;that} \; \; 1\le x \le 2
\end{equation*}
can be written as the convex minimization
\begin{equation*}
\min_x \left\{ \frac{1}{2} x^2 +\delta_S (x) \right\} \text{ where } S=\{x| x \in [1,2]\} ,
\end{equation*}
using Eq. (\ref{infadd}).

We are now in a position to analyze the saddle-point optimization
\begin{equation*}
\min_{x_1}\max_{x_2} x_1x_2.
\end{equation*}
Performing the maximization over $x_2$ first
\begin{equation*}
\max_{x_2} x_1 x_2,
\end{equation*}
three cases need to be considered: $x_1<0$, $x_1=0$, and $x_1>0$.
For $x_1<0$, the maximizer is $x_2=-\infty$ and the maximum for this case is $\infty$.
For $x_1=0$, the maximum for this case is $0$.
For $x_1>0$, the maximizer is $x_2=\infty$ and the maximum is $\infty$.
Putting these cases together, we have
\begin{equation*}
\max_{x_2} x_1 x_2 = \delta_Z(x_1) \text{ where } Z=\{0\}.
\end{equation*}
Minization over $\delta_Z(x_1)$ is trivial; the minimizer is $x_1=0$.
Thus the saddle-point is identified to be $x_1=x_2=0$.
Note that in performing this saddle-point optimization it is necessary to use the rules for
multiplication with $\infty$.

\section{The Legendre-Fenchel transform}
\label{app:legendre}

The Legendre-Fenchel (LF) transform, also known as convex conjugation, is one of the main
operations in convex analysis. It essentially provides a way to represent a convex function
in terms of linear functions that support its epigraph. For scalar functions of a vector,
this concept generalizes straight-forwardly to planar support functions.
It is useful for manipulating and simplifying optimization
problems when the order of optimization
operations can be interchanged. In such cases minimization/maximization can be performed on the component
lines instead of the function itself. This capability with convex analysis is analogous to Fourier
analysis, where functions are decomposed in plane waves and integration/differentiation operations can
be performed on the individual plane wave components. See Table 1 of Ref. \cite{Komodakis15} for a more
complete comparison of Fourier and convex analysis. In the following, we present the discussion in terms
of a convex function $f(x)$ in 2D, i.e. $x$ is a scalar, and occasionally insert remarks on generalizing
$x$ to a vector. All formulas, however, are written so that they apply to the multi-dimensional vector case.

\begin{figure}[h!]
\centering
\includegraphics[width=0.5\textwidth]{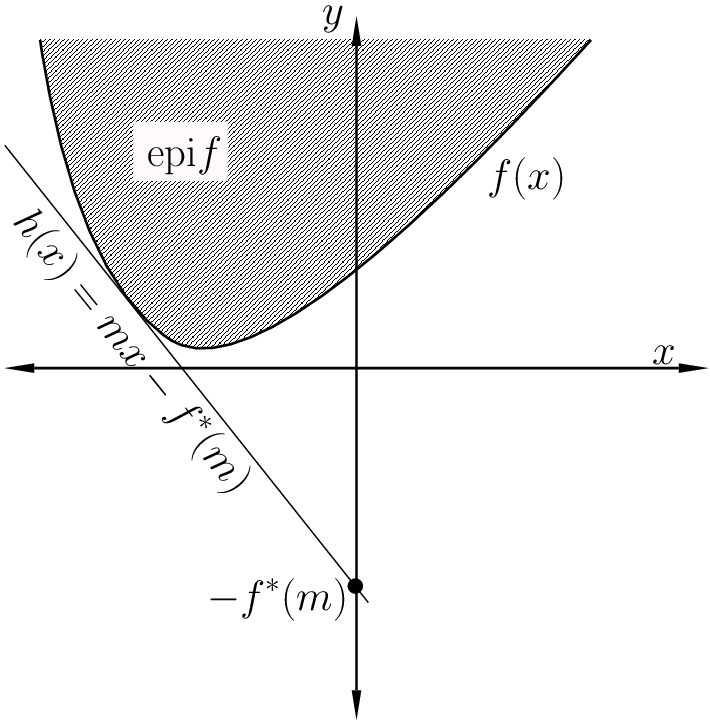}
\caption{\label{fig:legendre} Schematic of a convex function $f(x)$ and its epigraph.
The indicated line satisfies $h(x^\star) = f(x^\star)$ and $h(x) \le f(x)$. 
The line $h(x)$ intercepts the $y$-axis at $-f^*(m)$.}
\end{figure}

The LF transform essentially yields a parameterization of the epigraph supporting
lines in terms of their slope.
The convex function $f(x)$ can be represented as the max over all lines that lie
below $f(x)$. Parameterizing such a line by a slope $m$ and ``$y$-intercept'' $-b$
(in the multi-dimensional case, $m$ is a vector of the same dimension as $x$, and $b$ is still a scalar)
\begin{equation}
\label{hineq}
h(x) = m^\top x - b,
\end{equation}
the constraint that $h(x)$ lies below $f(x)$ is enforced by constraining $b$
to satisfy
\begin{equation}
\label{fstarcond}
b \ge  m^\top x - f(x).
\end{equation}
The set of $m$ and $b$ that specify lines that satisfy this inequality are denoted by the set $F^*$
\begin{equation*}
(m,b) \in F^*,
\end{equation*}
if Eq.~(\ref{hineq}) holds for all $x$.
We then can write the function $f(x)$ as a maximization 
\begin{equation}
\label{linerep}
f(x) = \max_{(b,m) \in F^*} \left\{ m^\top x - b \right\}.
\end{equation}
This representation of $f(x)$ can be reduced to just the supporting lines, i.e. the lines
that actually intersect $f(x)$ at at least one point, by restricting $b$ to the smallest
possible value for a given slope $m$.
\begin{eqnarray*}
b \ge  m^\top x - f(x), \\
b \ge \max_x \left\{m^\top x - f(x)\right\},\\
b_{\rm min} (m) = \max_x \left\{m^\top x - f(x)\right\}.
\end{eqnarray*}
The LF transform of $f(x)$, $f^*(x)$, is defined to be $b_{\rm min} (m)$, i.e.
\begin{equation}
\label{legfendef}
f^*(m) = \max_x \left\{m^\top x - f(x)\right\},
\end{equation}
and by construction $F^*$ is the epigraph of $f^*(m)$.
The LF transform has a clear geometric meaning, shown in Fig. \ref{fig:legendre}, which can be exploited to compute 
transforms of specific functions in addition to analysis techniques.

Because $F^*$ is the epigraph of $f^*(m)$, we can perform the maximization over $b$ in
Eq.~(\ref{linerep}) yielding
\begin{equation}
\label{legfeninv}
f(x) = \max_{m} \left\{ m^\top x - f^*(m) \right\}.
\end{equation}
Comparing Eqs.~(\ref{legfendef}) and (\ref{legfeninv}), we see that the LF transform is its
own inverse, and
\begin{equation*}
f^{**} = f,
\end{equation*}
provided that $f$ and $f^*$ are convex.
Because of this relation the LF transform is oftened referred to as convex conjugation.
If $f$ is not convex, we can not write Eq.~(\ref{linerep}) and
Eq.~(\ref{legfeninv}) does not hold.
Whether or not $f$ is convex, we can still compute its LF transform $f^*$ with
Eq.~(\ref{legfendef}) and in doing so, $f^*$ will be convex.
Also, if $f$ is non-convex, $f^{**}$ is the tightest convexification
of $f$. 

We present the argument that $f^*(m)$ is convex, if it is computed from Eq.~(\ref{legfendef}).
Considering the points $(x,z) \in F$, where $F$ is the epigraph of $f(x)$,
the maximization in Eq.~(\ref{legfendef}) can be written
\begin{equation}
\label{linerepdual}
f^*(m) = \max_{(x,z) \in F} \left\{m^\top x - z \right\},
\end{equation}
where we have used the fact that point $(x,z)$ lies above $(x,f(x))$.
Thus, we see that $f^*(m)$ is a maximization over lines in $m$-space
\begin{equation*}
g(m) = x^\top m - z,
\end{equation*}
where $x$ plays the role of a slope and $-z$ is the corresponding $y$-intercept.
Because $f^*(m)$ is a maximization over convex functions it is itself convex.

 Nevertheless, even if $f$ is non-convex, $f^*$ is convex
because Eq.~(\ref{linerepdual}) still holds and in this case $f^{**}$ is the tightest convexification
of $f$. Because $f^{**}=f$ for convex functions, the LF transform is often referred to as convex conjugation.

\subsection{LF transform examples}
In order to illustrate the main approaches to computing the LF transform, we find $f^*$
for quadratic,  absolute value, linear, and indicator functions.

\subsubsection*{(I) $f(x) = a x^2/2$: }
For a differentiable convex function
one can use the standard optimization technique of setting the gradient of the objective function
to zero, solving for the maximizer, and then plugging the maximizer back into the objective function
to obtain the function max.
Starting from Eq.~(\ref{legfendef}), we differentiate the objective function with respect to $x$ and set to 0
\begin{equation*}
\partial f(x^\star) = m ,
\end{equation*}
where $x^\star$ denotes the maximizer of Eq.~(\ref{legfendef}), i.e. the value of $x$ where the slope of $f(x)$ is $m$.
For the quadratic example, we have
\begin{eqnarray*}
a x^\star = m ,\\
x^\star = m/a.
\end{eqnarray*}
We plug this back into the LF objective function
\begin{equation*}
f^*(m) = m^\top x^\star - f(x^\star) = m^2/(2 a),
\end{equation*} 
obtaining another quadratic, with inverse width, as $a$ now appears in the denominator.
This example also hints at the duality nature of the LF transform; as an exercise it is worthwhile
to show that applying the LF transform again will yield a quadratic with $a$ back in the
numerator -- the same function we started with.

The multi-dimensional case is a trivial extension, because the LF objective function separates.
Generalizing the one-dimensional quadratic, we have
\begin{equation*}
f(x) = x^\top A x/2,
\end{equation*}
where $A=\diag{a}$ is a diagonal matrix with all positive diagonal elements $a>0$.
The LF transform is
\begin{eqnarray*}
f^*(m) &= \max_x \left\{ m^\top x - x^\top A x/2 \right\}, \\
       &= \sum_i \max_{x_i} \left\{ m_ix_i - a_i x_i^2 \right\}, \\
       &= \sum_i  m_i^2/(2 a_i),\\
       &= m^\top A^\prime m/2,
\end{eqnarray*}
where $A^\prime =\diag{1/a}$.

\begin{figure}[h!]
\centering
\includegraphics[width=0.5\textwidth]{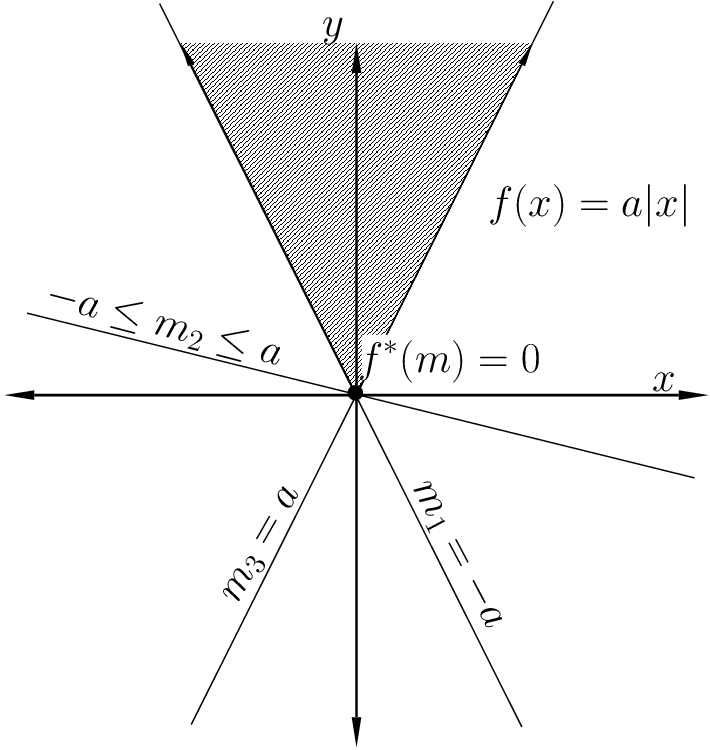}
\caption{\label{fig:ccabs} Schematic of a convex function $f(x)=a|x|$ and its LF transform.
Three support lines for this function are indicated. The lines labelled
$m_1$, $m_2$, and $m_3$ are respectively the lines with lowest possible slope, a generic
support line, and the line with highest possible slope. All support lines
intersect the $y$-axis at 0; hence the LF transform is zero for the allowed
slopes $-a \le m \le a$.}
\end{figure}

\subsubsection*{(II) $f(x) = a |x|$: }
The LF transform of the absolute value can be handled analytically, considering the cases where the minimizer is 
positive or negative.
There is, however, a much simpler geometric approach using Fig. \ref{fig:legendre}, which
we adapt to the function of interest in Fig. \ref{fig:ccabs}.
From this figure, the support lines to $a |x|$ all have slope between $-a$ and $a$, and
the $y$-intercept of all the supporting lines is 0, hence
\begin{equation*}
f^*(m) = \delta(-a \le m \le a).
\end{equation*}
The use of the indicator allows us to restrict the domain of slopes to those of the support lines,
and when the slope is between $-a$ and $a$ the indicator's value is zero.

For the multi-dimensional generalization, we consider
\begin{equation*}
f(x) = a \|x\|_1,
\end{equation*}
where $a$ in this case is still a scalar.
It is possible to employ the purely geometric approach, but
it is simpler to take advantage of the separable objective function
\begin{eqnarray*}
f^*(m) &= \max_x \left\{ m^\top x - \|x\|_1 \right\}, \\
       &= \sum_i \max_{x_i} \left\{ m_ix_i - a| x_i| \right\}, \\
       &= \sum_i  \delta (-a \le m_i \le a),\\
       &=   \delta (\|m\|_\infty \le a),
\end{eqnarray*}
where $\| \cdot \|_\infty$ yields the magnitude of the largest component of its argument
\begin{equation*}
\| m \|_\infty = \max \left( |m|_1, |m|_2, \dots, |m|_i,\dots \right).
\end{equation*}
In summary, the $\ell_1$ norm LF transform is most easily dealt with by a combination of geometric
and analytic methods; the LF optimization problem is seperated analytically and the individual
one-dimensional optimation problems are handled geometrically.

\subsubsection*{(III) $f(x) = a x + c$: }
Using the geometric approach for a linear function is also quite straight-forward. There is only one
support line with slope $m=a$. From Fig. \ref{fig:legendre} the value of $f^*$ is the negative $y$-intercept,
and the $y$-intercept for the line of interest is $c$. As a result, $f^*$ for $m=a$ has the value of $-c$
\begin{equation*}
f^*(m) = \delta(m=a) - c.
\end{equation*}
The indicator function reduces the domain of the function to a single point in this case, and subtracting
$c$ does not alter this because of Eq.~(\ref{infadd}). This same argument and formula applies for
the multi-dimensional case of $f(x) = a^\top x + c$.

\begin{figure}[h!]
\centering
\includegraphics[width=0.6\textwidth]{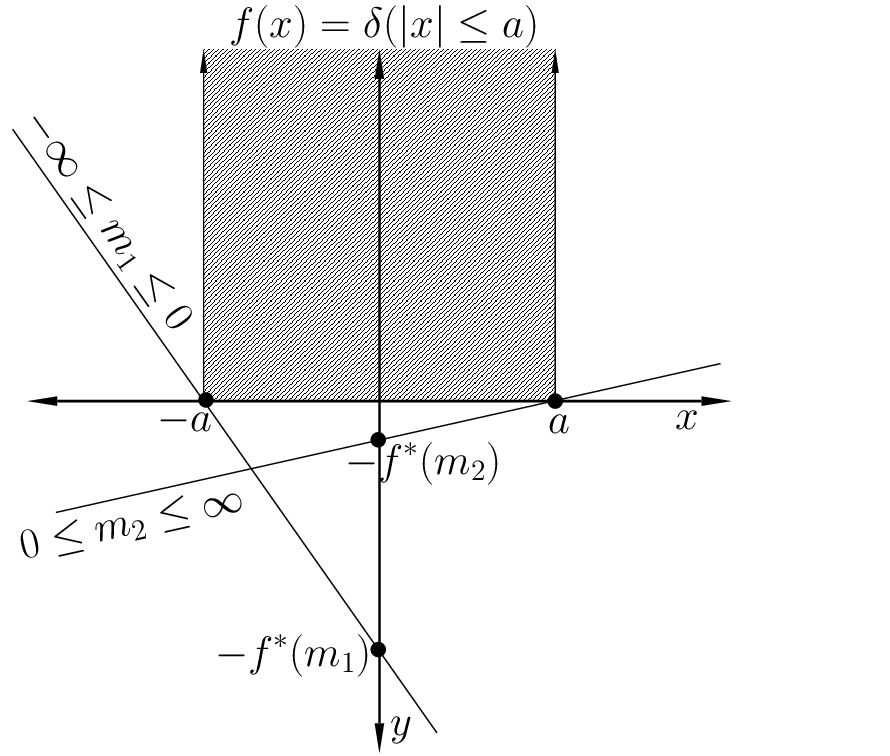}
\caption{\label{fig:indseg} Schematic of a convex function $f(x)=delta(|x| \le a)$ and its LF transform.
Two support lines for this function are indicated. The lines labelled
$m_1$ and $m_2$ are the generic support lines that intersect $f(x)$ at $x=-a$ and $x=a$, respectively.
Note that one support line with $m=0$ intersects $f(x)$ over the whole segment between $x=-a$ and $x=a$.
Also shown are the $y$-intercepts for the lines $m_1$ and $m_2$. By geometric reasoning it is clear
that $f^*(m_1) = -a m_1$ and $f^*(m_2) = a m_2$. Putting these cases together yields the result $f^*(m) = a|m|$.}
\end{figure}

\subsubsection*{(IV) indicator functions $f(x) =\delta(-a \le x \le a)$ and $f(x)= \delta(x=a) - c$: }
For the indicator function examples we choose exactly the same functions we arrived at from the previous
two examples of LF transforms. In both cases it is simplest to use the geometric approach. For
\begin{equation*}
f(x) =\delta(-a \le x \le a),
\end{equation*}
we see in Fig. \ref{fig:indseg} that lines of all slopes $m$ contribute to the function support
and the desired LF transform can be extracted from the negative $y$-intercept of the drawn lines.
We thus obtain
\begin{equation*}
f^*(m) = a |m|.
\end{equation*}
For
\begin{equation*}
f(x) = \delta(x=a) -c,
\end{equation*}
the same approach yields its LF transform
\begin{equation*}
f^*(m) = a m + c.
\end{equation*}
In both cases we observe that the LF transform has inverted the previous two examples.

\section{The subdifferential and subgradient}
\label{app:subdiff}

For convex functions the subdifferential is a useful generalization
of standard differentiation of smooth functions.
The subdifferential $\partial f(x)$ is a set-valued mapping defined by
the following inequality
\begin{equation}
\label{subdiff}
\partial f(x) = \left\{ m \, |\, \forall x^\prime : \, f(x^\prime) \geq f(x) + m^\top (x^\prime-x) \right\}.
\end{equation}
For functions of a scalar, the subdifferential at any point $x$ is the set of slopes of lines
that pass through $f(x)$ but lie completely underneath $f(x)$.
Equation (\ref{subdiff}) expresses the $n$-dimensional generalization of this idea.
For differentiable $f(x)$, $\partial f(x)$ yields the usual gradient.
A sub-gradient $g$ is one of the elements
of $\partial f(x)$
\begin{equation*}
g \in \partial f(x).
\end{equation*}
The subdifferential is useful for convex functions because there will always be at least one linear
function that goes through $f(x)$ and lies completely beneath $f(x)$. Also, in this paper we
are primarily concerned with first-order algorithms and optimality conditions, and first-order
subdifferentiation is uncomplicated.

\begin{figure}[h!]
\centering
\includegraphics[width=0.7\textwidth]{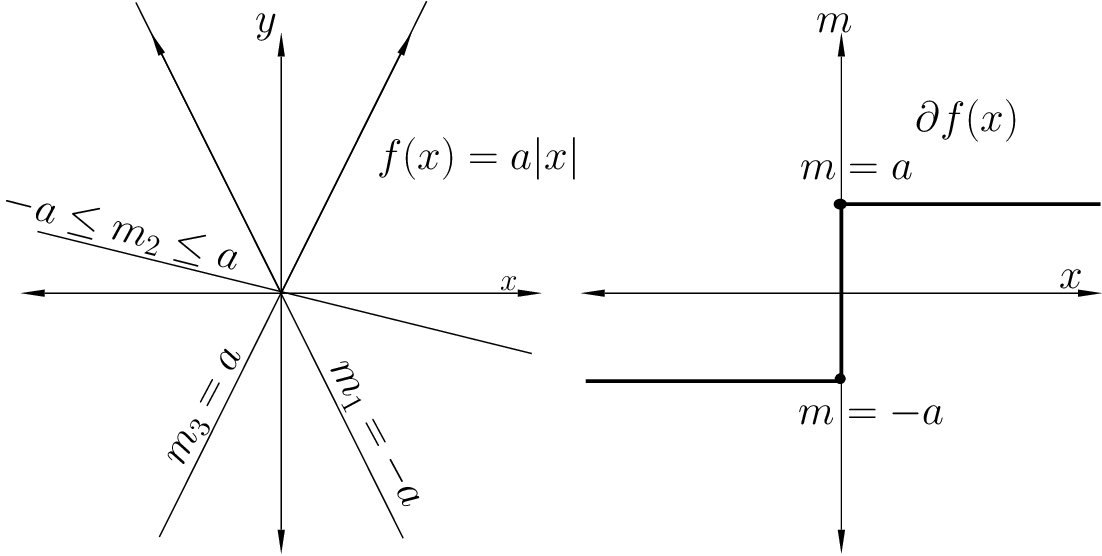}
\caption{\label{fig:sdabs} Schematic of a convex function $f(x)=a|x|$ (left) and its subdifferential (right).
The lines indicated by the slope $m_1$, $m_2$, and $m_3$ yield the subdifferential for $x<0$, $x=0$,
and $x>0$, respectively. For $x \neq 0$ this function is smooth, and accordingly the subdifferential
is the standard gradient. For $x =0$, lines with slope $|m_2| \leq a$ all fit underneath the function
while intersecting $f(0)=0$. Thus for this point the subdifferential is the set $[-a,a]$.}
\end{figure}
The classic example illustrating the subdifferential for a non-smooth function
is $\partial f(x)$ when
\begin{equation*}
f(x) = a|x|.
\end{equation*}
This function and its subdifferential are illustrated in Fig. \ref{fig:sdabs}.
The subdifferential is multi-valued only at $x=0$ where the function is not-differentiable.
In particular, we note that the subgradient $g=0$ is a member of the subdifferential $\partial f(x=0)$
and at the same time $x=0$ is the minimum of this function.
In fact, the first-order condition of optimality for a convex, non-smooth function
is
\begin{equation*}
0 \in \partial f(x^\star),
\end{equation*}
which is the generalization of the condition
\begin{equation*}
0 = \partial f(x^\star),
\end{equation*}
for differentiable $f$.
We note also that a minimum at a non-smooth point lends robustness to the minimizer; i.e. perturbing
the function by a smooth function will not change the minimizer unless the slope of the perturbation
exceeds the extreme values of the subdifferential at the minimizer.

In the main text, the subdifferential concept is not absolutely necessary for the development.
For non-smooth
$F$ or $F^*$, $\partial F$ and $\partial F^*$ can be multivalued in these equations.
Also, the maximizer in the LF transform objective function is specified by the subdifferential.
This is seen in the similarity between Figs. \ref{fig:ccabs} and the left-hand graph in
\ref{fig:sdabs}. While we do not absolutely need the concept of the subdifferential to compute the
LF transform as we saw in \ref{app:legendre}, we can now express the maximizer of the LF objective function
in Eq.~(\ref{legfendef})
formally as $x^*$ satisfying the equation
\begin{equation*}
m \in \partial f(x^*).
\end{equation*}

\section{Saddle point solver intuition}
\label{app:intuition}

Saddle point optimization is fundamentally different than convex function minimization
or concave function maximization. A saddle potential has at least one direction of
negative curvature and one direction of positive curvature, and accordingly the minimum
dimension for saddle point optimization is two, while both minimization and maximization
can be performed for functions of scalars. As a result,
the intuition for saddle point optimization is more complex than that of minimization/maximization.
Specifically, for the latter one can imagine computing the function gradient and taking a step
in that direction for maximization or in the opposite direction for minimization.
For saddle point optimization, using only first-order or gradient information, an algorithm
needs to go uphill, with the gradient, for variables that are being maximized over and downhill,
against the gradient, for variables that are being minimized over. While the decomposition
of the gradient to form a step-direction for a saddle point optimization algorithm seems straightforward,
complications arise when the coordinates are not well aligned with the directions
of curvature, and this is the reason why the CPPD algorithm is significantly more complex
than basic gradient descent/ascent for minimzation/maximization.

In this appendix, we present examples of saddle point optimization that serve to motivate
the particular form of the CPPD update steps and to illustrate convergence behavior as a function
of algorithm parameters. This appendix contains examples of forward Euler iteration that gives
a more complete picture of its behavior with saddle point problems; a couple examples of the
approximate backward Euler iteration then motivate the particular choice of CPPD algorithm
parameter settings; and then finally the convergence behavior of the CPPD applied to
least-squares minimization is explained by use of an eigenvector decomposition.

\subsection{Saddle point optimization with forward Euler iteration}
As discussed in Sec. \ref{sec:fei} the forward Euler iteration for finding the saddle point
of
\begin{equation*}
s(x,\lambda) = \lambda^\top Ax
\end{equation*}
does not converge for any step size. Recall that the forward Euler iteration
for this problem is
\begin{align*}
x_{k+1}&=x_k- \alpha A^\top \lambda_k \\
\lambda_{k+1}&=\lambda_k+ \alpha A x_k,
\end{align*}
where $k$ is the index number for the iteration.
The simplest special case of this problem occurs when both $x$ and $\lambda$
are scalars and $A=1$
\begin{equation*}
s_0(x,\lambda) = \lambda x,
\end{equation*}
and the corresponding forward Euler iteration is
\begin{align*}
x_{k+1}&=x_k- \alpha  \lambda_k \\
\lambda_{k+1}&=\lambda_k+ \alpha  x_k.
\end{align*}
The saddle point of $s_0(x,\lambda)$ is at $x=\lambda=0$, but if the
forward Euler iteration is initialized away from this saddle point, the
subsequent iterations will spiral away from the origin of $x\lambda$-plane.
If the current iterate is $x_k,\lambda_k$, its distance from the origin
is
\begin{equation*}
r_k = \sqrt{x_k^2 + \lambda_k^2}.
\end{equation*}
Using the update equations, one can show that the distance of the next
iterate is related to the current distance by
\begin{equation*}
r_{k+1} = \sqrt{1+\alpha^2} \; r_k.
\end{equation*}
Clearly, the distance of the iterates from the origin will increase for any
value of the step size parameter $\alpha$.
This example, however, does not mean that forward Euler iteration always
fails to converge for saddle point problems.

Consider a different potential
\begin{equation*}
s_1(x,\lambda) = x^2 - \lambda^2,
\end{equation*}
where $x$ and $\lambda$ are scalars.
This potential also has a critical point at $x=\lambda=0$, which is a saddle point.
In this case the forward Euler updates, derived by a taking a step in the direction
of the derivative of $\lambda$ and opposite to the direction of the derivative in $x$,
are
\begin{align*}
x_{k+1}&=x_k- 2 \alpha  x_k = (1-2\alpha) x_k \\
\lambda_{k+1}&=\lambda_k- 2\alpha  \lambda_k = (1-2\alpha)\lambda_k.
\end{align*}
By inspection, it is clear that the relationship between successive distances
to the origin for these update formulas is
\begin{equation*}
r_{k+1} = (1-2\alpha) r_k,
\end{equation*}
and the choice of $0<\alpha<1$ leads to convergence to the saddle point
at the origin of the $x \lambda$-plane. Interestingly, the potentials $s_0$ are $s_1$
related by rotation of 45 degrees; i.e. the variable substitution $x^\prime = x+\lambda$
and $\lambda^\prime = x -\lambda$ turns one of these potentials into the other, up to a scalar
multiple.

The fact that $s_0$ and $s_1$ are related by simple rotation would seem to suggest 
a possible means to attack the general saddle point problem
\begin{equation*}
\min_x \max_\lambda \left\{ \lambda^\top A x - \phi^*(\lambda) \right\},
\end{equation*}
where maybe a coordinate change would allow forward Euler iteration to be successfully
applied. The barrier to this strategy is that the
equivalent to performing the rotation of $s_0$ to $s_1$ involves second-order
information, namely computing the eigendecomposition of the Hessian of this saddle potential.
This can be prohibitively expensive for large-scale optimization.

Fortunately, there is a
way to address the saddle point optimization problem of interest using only first order
information as given by the approximate backward Euler iteration in Sec. \ref{sec:abei}.
As there are a number of parameters in the original form of this algorithm, it is
illustrative to examine it in the special case saddle point potentials $s_0(x,\lambda)$
and another one corresponding to one-dimensional quadratic optimization.

\subsection{Saddle point optimization with approximate backward Euler iteration}
Recall from Sec. \ref{sec:abei} that the approximate 
backward Euler step can optimize the saddle potential $s(x,\lambda)$, and the
update steps given in Eqs. (\ref{cps2A})-(\ref{cps2C}) are repeated here
\begin{eqnarray*}
x_{k+1} =& x_k - \tau A^\top \lambda_{k},\\
\bar{x}_{k+1} =& x_{k+1} + \theta (x_{k+1}  - x_{k}), \\
\lambda_{k+1} =& \lambda_k + \sigma A \bar{x}_{k+1}.
\end{eqnarray*}
There are three algorithm parameters $\sigma$, $\tau$, and $\theta$.
According to Ref. \cite{chambolle2011first}, these parameters should be
chosen, respecting the following inequalities
\begin{equation*}
0 \le \theta \le 1, \; \; \; \sigma \tau < \|A\|_2,
\end{equation*}
where the latter strict inequality is can be taken as $\ge$ in most cases
of practical interest.

Applying this algorithm to the two-dimensional saddle point problem $s_0(x,\lambda)$,
$A=1$ and the update steps reduce to
\begin{eqnarray*}
x_{k+1} =& x_k - \tau \lambda_{k},\\
\lambda_{k+1} =& \lambda_k + \sigma ((1+\theta) x_{k+1} - \theta x_k),
\end{eqnarray*}
where the extrapolation step is absorbed into the $\lambda$ update.
Parameterizing the step-lengths in terms of their product $a$ and allowing
the equality case for this product, $\tau$ can be written in terms
of $\sigma$ and $a$
\begin{equation*}
\sigma \tau = a \le 1.
\end{equation*}
The update steps can thus be manipulated into the following matrix-vector product form
\begin{equation*}
\left( \begin{array}{c} x_{k+1} \\ \lambda_{k+1} \end{array} \right) =
\left( \begin{array}{cc} 1 & -a/\sigma \\ \sigma & 1 - a - \theta a \end{array} \right) 
\left( \begin{array}{c} x_{k} \\ \lambda_{k} \end{array} \right).
\end{equation*}
Considering a couple of special cases provides some orientation on the parameter
dependences of the iteration with this update.

\subsubsection*{$\theta=1$ and $a=1$:}
These parameter settings yield
\begin{equation*}
\left( \begin{array}{c} x_{k+1} \\ \lambda_{k+1} \end{array} \right) =
\left( \begin{array}{cc} 1 & -1/\sigma \\ \sigma & -1 \end{array} \right) 
\left( \begin{array}{c} x_{k} \\ \lambda_{k} \end{array} \right),
\end{equation*}
and it is straight-forward to verify that
\begin{equation*}
\left( \begin{array}{c} x_{k+2} \\ \lambda_{k+2} \end{array} \right) =
\left( \begin{array}{cc} 0 & 0 \\ 0 & 0 \end{array} \right) 
\left( \begin{array}{c} x_{k} \\ \lambda_{k} \end{array} \right).
\end{equation*}
This result means that no matter what $x_0$ and $\lambda_0$ are initialized
to, the updates will converge in two steps to the saddle point at $(x,\lambda) = (0,0)$,
independent of $\sigma$.

\subsubsection*{$\theta=0$ and $a=1$:}
These parameter settings yield
\begin{equation*}
\left( \begin{array}{c} x_{k+1} \\ \lambda_{k+1} \end{array} \right) =
\left( \begin{array}{cc} 1 & -1/\sigma \\ \sigma & 0 \end{array} \right) 
\left( \begin{array}{c} x_{k} \\ \lambda_{k} \end{array} \right).
\end{equation*}
From this update, one can show that
\begin{equation*}
\left( \begin{array}{c} x_{k+8} \\ \lambda_{k+8} \end{array} \right) =
\left( \begin{array}{cc} 1 & 0 \\ 0 & 1 \end{array} \right) 
\left( \begin{array}{c} x_{k} \\ \lambda_{k} \end{array} \right),
\end{equation*}
and it is clear that the updates will repeat values every eighth iteration and
not converge to the saddle point, again, independent of $\sigma$'s value.
Thus, we encounter a case where the step-size
product must respect the strict inequality $a<1$ in order to converge.

These two specific cases provide a heuristic for selecting $\theta=1$, but this particular
saddlepoint problem is too simple to provide intuition on why it is important
to tune the step-size ratio
\begin{equation*}
\rho = \sqrt{\sigma/\tau}.
\end{equation*}
No matter what value is chosen for $\sigma$, the two particular parameter settings
shown yield the same results.

\subsection{Tuning the step-size ratio of the CPPD algorithm}

In order to appreciate the impact of the step-size ratio in the CPPD algorithm
the minimum complexity problem to consider is one-dimensional quadratic optimization.
Specifically, consider the minimization
\begin{equation*}
\min \frac{1}{2} x^2.
\end{equation*}
Repeating the steps described in Sec. \ref{sec:CPPDback}, allows this minimization
to be generalized to the saddle point problem
\begin{equation*}
  \min_x \max_\lambda  \left\{x \lambda - \frac{1}{2} \lambda^2\right\},
\end{equation*}
which is a special case of Eq. (\ref{saddle2}). The solution for the 
minimization is clearly $x= 0$ and the saddle point for the second potential
is $(x,\lambda)=(0,0)$. The approximate backward Euler iteration for this
saddle potential is
\begin{equation*}
\left( \begin{array}{c} x_{k+1} \\ \lambda_{k+1} \end{array} \right) =
\left( \begin{array}{cc} 1 & \frac{-a}{\sigma} \\ \frac{\sigma}{1+\sigma} &
                                          \frac{1 - 2a}{1+\sigma}  \end{array} \right) 
\left( \begin{array}{c} x_{k} \\ \lambda_{k} \end{array} \right),
\end{equation*}
where $\theta=1$, $\sigma \tau=a$, and $a\le1$. Fixing $a$ and varying $\sigma$ effectively
varies the step-size ratio $\rho$. Even with $\theta$ set to one this update step is
complicated to analyze with analytic methods; a given update trajectory depends non-trivially
on the initial values $(x_0,\lambda_0)$, $a$, and $\sigma$.
For comparison, consider the gradient descent (GD)
update for the corresponding quadratic minimization
\begin{equation*}
x_{k+1} = x_k - a x_k = (1-a) x_k,
\end{equation*}
where convergence is attained for $0<a<2$ and
any iterate can be computed directly from the initial value by
\begin{equation*}
x_k = (1-a)^k x_0.
\end{equation*}
The trajectory of the iterates $x_k$ is clearly easier to characterize than the
trajectory of $x_k,\lambda_k$ for the CPPD algorithm.

\begin{figure}[h!]
\centering
\includegraphics[width=0.75\textwidth]{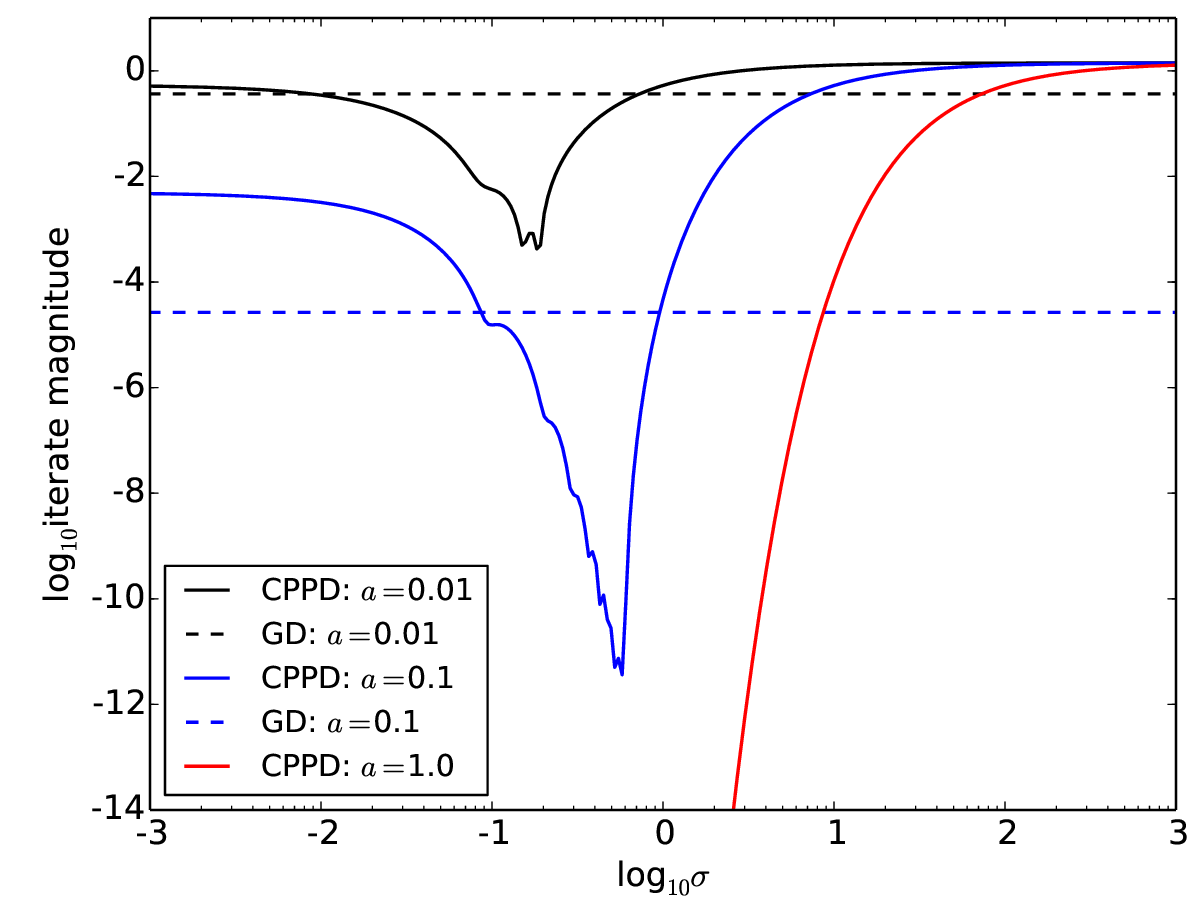}
\caption{\label{fig:cppdtraj}
Plotted is the magnitude of $(x_k,\lambda_k)$ after 100 iterations. The initial values
are $(x_0,\lambda_0) = (1,0)$; $\sigma$ is varied from $10^{-3}$ to $10^3$; and curves
are shown for $a=0.01$, $a=0.1$, and $a=1.0$. For comparison, the gradient descent (GD)
iterate magnitude $|x_k|$ after 100 iterations for the corresponding $a$ values. Note
that GD does not depend on $\sigma$. Also, GD is not shown for $a=1$, because convergence
is obtained in one iteration for this case.
}
\end{figure}
A specific case is shown in Fig. \ref{fig:cppdtraj} for one set of initial values,
where the magnitude of the iterate vector $(x_k,\lambda_k)$ is shown after 100
iterations and compared with the GD algorithm.
The iterate magnitude is an indication of convergence, because the solution to the
quadratic minimization and saddle point problem is $0$ and $(0,0)$, respectively.
For the case of $a=1$, CPPD converges
rapidly for small $\sigma$, but GD converges in one iteration. For $a=0.1$,
there is a clear minimum in the iterate magnitude at $\sigma \approx 0.5$ and that
magnitude is much smaller than what is attained by GD. Similarly for $a=0.01$,
there is a minimum at $\sigma \approx 0.2$, which is much smaller than the corresponding
GD result.  For $a=1$, GD outperforms CPPD in terms of convergence since it converges
in one iteration. But for the other $a$ values, CPPD outperforms GD provided that
$\sigma$ is tuned. If the $\sigma$ parameter is not tuned, the CPPD iterations can
converge very slowly, even more slowly than GD.

This example is actually relevant to large-scale least-squares optimization. Performing
an SVD of the system matrix leads to a set of uncoupled one-dimensional quadratic optimizations.
In such a decomposition and normalizing the system matrix to 1, the parameter $a$ takes on the
role of the eigenvalues of the system matrix. Performing the CPPD iteration for the least-squares
system is equivalent to selecting the same $\sigma$ value for all of the one-dimensional
quadratic sub-problems obtained by the SVD analysis. For example, the plot
in Fig. \ref{fig:cppdtraj} is useful for analyzing a system matrix with only three eigenvectors
with eigenvalues $0.01$, $0.1$, and $1.0$. In this case, we note that the $a=0.01$
curve is the largest for all shown $\sigma$ values; $\sigma$ should be chosen
by finding the minimum of this curve. 

This analysis, however, is only for illustrative purposes. The CPPD trajectory has a lot
of complexity not captured by Fig. \ref{fig:cppdtraj}. Furthermore, SVD of large-scale
system matrices may not be practical. In practice, it is simpler to directly tune
$\sigma$ or equivalently the step-size ratio $\rho$, for fixed $\sigma \tau$.

\section{Fixed-point iteration, the proximal point algorithm, resolvents, and monotone operators}
\label{app:fixedpoint}

To provide further insight into the CPPD algorithm, we summarize the fixed-point iteration
formalism in which convergence of this and other recent first-order algorithms can be readily
shown.  In this appendix, fixed-point iteration is explained with a couple of one-dimensional
examples. First, the familiar gradient descent algorithm is cast as a fixed-point iteration
in order to obtain intuition on how fixed-point iteration relates to optimization. Second,
a specific form of fixed-point iteration called the proximal point algorithm is presented
in order to appreciate its effectiveness with non-smooth optimization.
Finally,
the generalization of the proximal point algorithm to the CPPD and other related algorithms
requires the concept of the resolvent and monotone operators, which are briefly explained
here. For further reading on how the CPPD can be framed as a generalized proximal point
algorithm, please consult the article by He and Yuan \cite{he2012convergence}. An excellent
tutorial paper on fixed-point iteration by Ryu and Boyd \cite{ryu2016primer} explains a number
of recent first-order algorithms, including ADMM, ISTA, Douglas-Rachford, and CPPD.

A fixed point, $x^\star$, of an operator $M$ obeys the equation
\begin{equation*}
x^\star = M(x^\star).
\end{equation*}
An algorithm for finding such points is the fixed-point algorithm
\begin{equation*}
x^{(k+1)}  = M(x^{(k)}),
\end{equation*}
where an initial guess $x^{(0)}$ is made and subsequent iterations indexed by $k$
are obtained by feeding the output of $M$ back into $M$. Convergence of this algorithm
to $x^\star$ is discussed after showing gradient descent and proximal point examples.

\begin{figure}[h!]
\centering
\includegraphics[width=0.45\textwidth]{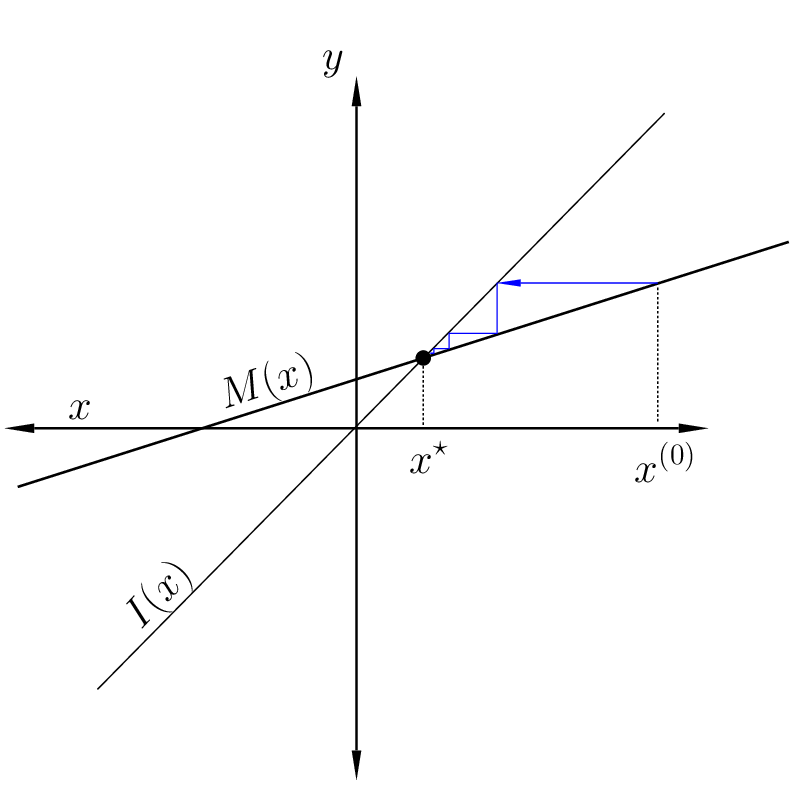}
\includegraphics[width=0.49\textwidth]{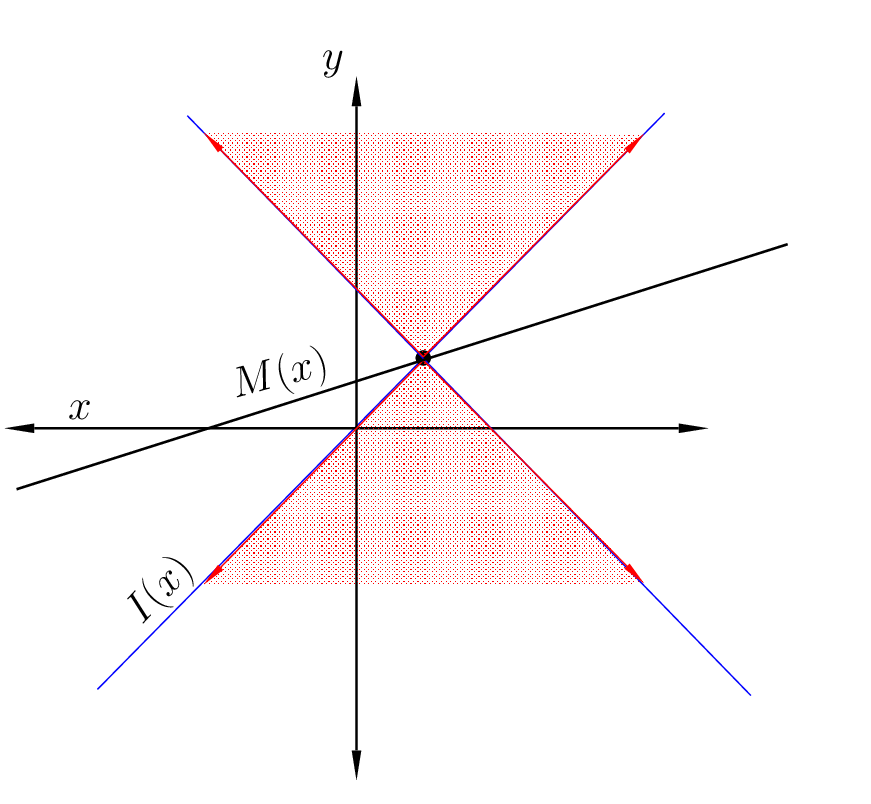}
\caption{\label{fig:fixedpt} Left: Schematic of fixed-point iteration corresponding to
gradient descent of a quadratic function $ f(x) = (1/2)( x - x^\star)^2$.
The fixed-point operator is $M(x) = (I - \alpha \nabla f)(x) = (1-\alpha) x + \alpha x^\star$.
Graphically, fixed-point iteration starts at $M(x_0)$ and alternates between horizontal
projection onto $I(x)$ and vertical projection onto $M(x)$. This iteration converges
to the minimizer $x^\star$, where $M$ and $I$ intersect, if $\alpha$ is chosen appropriately.
Right: Diagram indicating region where fixed-point iteration with $M(x)$ converges to $x^\star$.
Consider the following cases: $\alpha<0$, $M(x)$ is in the red zone because
it has a slope greater than 1
and fixed-point iteration will diverge; $\alpha=0$, $M$ and $I$ coincide and all
iterates equal $x_0$; $0< \alpha< 1$, iterations approach solution from one side; $\alpha=1$,
$M$ is horizontal and convergence to $x^\star$ is achieved in one iteration; $1<\alpha<2$,
iterates approach $x^\star$ but are under-relaxed as successive iterates are on opposite
sides of $x^\star$; $\alpha=2$ iterates oscillate between $x_0$ and $-x_0$; and finally
$\alpha>2$, $M$ is once again in the red zone and fixed-point iteration diverges.
}
\end{figure}
\subsection{Gradient descent as fixed-point iteration}
The gradient descent algorithm to find a minimizer $x^\star$ of a smooth convex function
$f(x)$, involves making an initial guess $x^{(0)}$
and performing the following iteration
\begin{equation}
\label{app:gd}
x^{(k+1)} = x^{(k)} - \alpha \nabla f(x^{(k)}),
\end{equation}
where for each estimate $x^{(k)}$ the gradient descent step involves subracting a
step length parameter $\alpha$ times the gradient at that estimate. In fixed-point operator form,
this iteration is written
\begin{equation*}
x^{(k+1)} = M(x^{(k)}), \; \; \;  M(x) = (I - \alpha \nabla f)(x),
\end{equation*}
where $I$ is the identity operator.
At the minimizer $x^\star$, the gradient of $f$ is zero and hence $x^\star = M(x^\star).$
Whether or not an arbitrary initial guess will converge toward $x^\star$ depends on the step-size $\alpha$.

To illustrate convergence of gradient descent, we take a simple example
of a smooth convex function
\begin{equation*}
f(x) = \frac{1}{2} (x - x^\star)^2.
\end{equation*}
By inspection the minimizer is $x^\star$, and because this is a one-dimensional example
the operator $M$ maps a scalar $x$ to a scalar $y$.
The mechanics and convergence of the fixed-point iteration are illustrated in Fig. \ref{fig:fixedpt}.
From this figure, we can make a couple general observations on conditions for converge
to the minimizer $x^\star$. First, $x^\star$ must be a fixed point of $M$; i.e. $I$ and $M$ intersect at $x^\star$.
Second, if $M$ is either non-decreasing or non-increasing and the magnitude of its slope is less than 1, the
fixed-point iteration will converge to $x^\star$. In stating the second condition, we have appealed to the
fact that our example $M$ maps a scalar to a scalar, so we can talk about $M$ being increasing or decreasing
with $x$ and the slope of $M$ has meaning. For generalization to $n$-dimensional mappings see \ref{sec:mono}.

\begin{figure}[h!]
\centering
\includegraphics[width=0.3\textwidth]{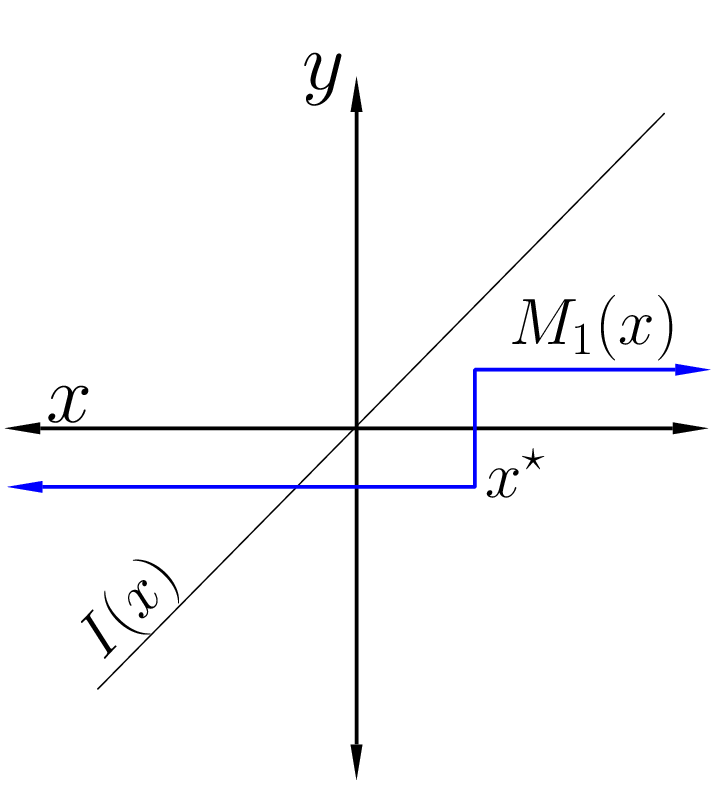}
\includegraphics[width=0.3\textwidth]{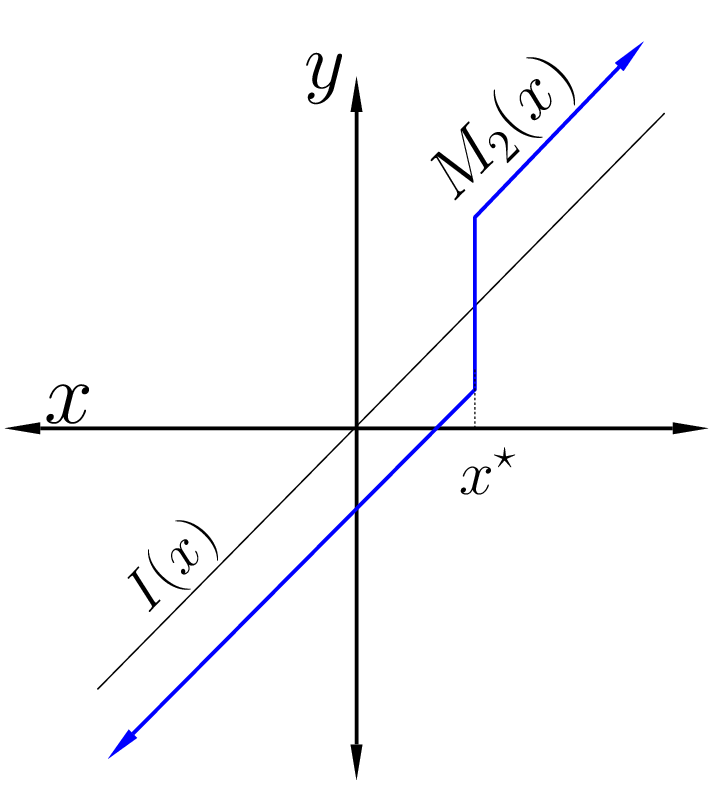}
\includegraphics[width=0.3\textwidth]{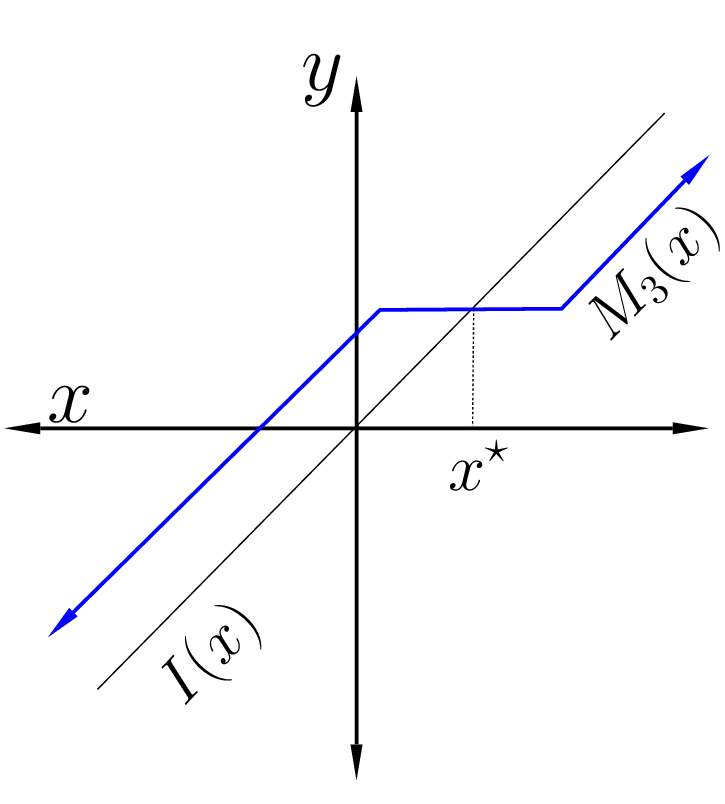}
\caption{\label{fig:proxsteps} The three graphs illustrate construction of the proximal
operator of $f(x) = |x-a|$. Left: Illustration of $M_1(x) = \partial f(x)$. Middle:
Illustration of $M_2(x) = (I + \alpha \partial f)(x)$. Finally, inversion of $M_2$ is simply
a matter of reflecting it about the 45 degree line representing $I(x)$.  Right: Illustration
of $M_3(x) = (I + \alpha \partial f)^{-1}(x) \equiv \prox_{\alpha f}(x)$.}
\end{figure}
\subsection{The proximal point algorithm for non-smooth optimization}
The proximal mapping discussed in Sec. \ref{sec:cppd} arises from the backward Euler step
for gradient descent
\begin{equation}
\label{app:bgd}
x^{(k+1)} = x^{(k)} - \alpha \partial f(x^{(k+1)}),
\end{equation}
which differs from Eq. (\ref{app:gd}) in that the gradient is replaced by the sub-differential
and the sub-differential is evaluated at $x^{(k+1)}$ instead of $x^{(k)}$.
Solving for $x^{(k+1)}$, can be performed with operator algebra
\begin{align*}
& x^{(k+1)} +  \alpha \partial f(x^{(k+1)}) = x^{(k)} \\
& (I + \alpha \partial f) x^{(k+1)} = x^{(k)} \\
& x^{(k+1)}=  (I + \alpha \partial f)^{-1} x^{(k)} = \prox_{\alpha f} ( x^{(k)} ),
\end{align*}
where
\begin{equation}
\label{proxdef}
\prox_{\alpha f} ( x ) \equiv
\argmin_{x^\prime} \left \{ \alpha f(x^\prime) +\frac{1}{2} \|x - x^\prime \|_2^2 \right\}.
\end{equation}
From the above derivation, the fixed-point operator for the proximal point algorithm
is
\begin{equation*}
M(x) = (I + \alpha \partial f)^{-1} x,
\end{equation*}
which can be understood in an intuitive way by breaking down its computation with a simple
graphical example.

The series of graphs in Fig. \ref{fig:proxsteps} shows the steps to deriving the proximal operator
for a simple one-dimensional non-smooth convex function
\begin{equation*}
f(x) = |x - x^\star|,
\end{equation*}
which, again, has a minimum value of zero at $x^\star$.
Going from left to right, the first graph shows $\partial f$ which
has a vertical segment at $x^\star$ corresponding to the multi-valuedness
of the sub-differential. The second graph illustrates $I+ \alpha \partial f$.
For the third and final graph, operator inversion is performed, which simply involves
exchanging $x$ and $y$ or reflection about the 45 degree line. Note that the vertical
line segment becomes horizontal and the proximal operator is single-valued.

The graphical derivation in Fig. \ref{fig:proxsteps} also makes clear a couple of properties
of the proximal mapping of convex functions, which are possibly non-smooth, and
that have a finite minimum. 
First, the proximal mapping converges for any $\alpha$. The argument that this is
the case goes as follows: the subdifferential
of any convex function $f$ is non-decreasing as $x$ increases;
thus $I+ \alpha \partial f$
is strictly increasing with slope greater than or equal to 1; and this in turn
implies that $(I + \alpha \partial f)^{-1}$ is non-decreasing with slope between zero
and one, avoiding the red zones indicated in Fig. \ref{fig:fixedpt}. Second, the proximal mapping
is single-valued: $I+ \alpha \partial f$  cannot have any horizontal segments because
adding the identity to a non-decreasing function yields a strictly increasing function,
thus reflection about the 45 degree line does not allow $(I + \alpha \partial f)^{-1}$
to have any vertical segments, which would indicate multi-valuedness.

\subsection{Monotone operators, the Lipschitz constant, and the resolvent}
\label{sec:mono}

[When this comes up in the saddle point discussion point out that the monotone operator is not a necessarily a gradient]

For general operators $M$ that map $n$-dimensional vectors to
$n$-dimensional vectors, the concept of non-decreasing or non-increasing in one-dimension
is replaced by monotonicity. 
A monotone operator satisfies the condition
\begin{equation}
\label{monotone}
(M(u) -M(v))^\top (u-v) \ge 0 \; \; {\rm for \; all}\; u,v.
\end{equation}
The gradient of a convex function is an example of a monotone operator.
The reason why we need the concept of a montone operator is that it is more
general; not all monotone 
operators can be written as a gradient of a convex function.
The one-dimensional concept of slope is replaced by the Lipschitz constant, which is
the smallest positive real number $L$ such that
\begin{equation}
\label{lipschitz}
\|M(u) - M(v)\|_2 \leq L \|u-v\|_2 \; \; {\rm for \; all}\; u,v.
\end{equation}
Note that we break the convention that capital letters are operators or matrices in the case of $L$.
If $L<1$, $M$ is a contraction, and if $L \le 1$, $M$ is non-expansive.
If $M$ is a contraction, it has a fixed point and fixed point iteration will
arrive at it. The latter is seen easily by letting $v$ be the fixed point: repeatedly
applying $M$ to $v$ will yield $v$, and as a result $M(u)$ is at least a factor $L<1$
closer to $v$ than is $u$.
If $M$ is merely non-expansive, it may or may not have a fixed point and fixed point
iteration might work. Monotone operators can be combined to form other monotone operators
thereby generating new fixed point algorithms. Gradient descent and the proximal point
algorithms, which combine the montone operators $I$ and $\partial f$, are examples of this.

As described in Refs. \cite{he2012convergence,ryu2016primer}
the CPPD algorithm can also be written as fixed point iteration. This requires two generalizations
from the proximal point algorithm. The first generalization is the montone operator called the
resolvent $R$, formed from the operator $F$ by
\begin{equation*}
R = (I + \alpha F)^{-1}.
\end{equation*}
If $F$ is monotone, $R$ is non-expansive.
The proximal operator is the resolvent of $\partial f$, but $R$ is more general than prox
because $F$ may not be the subdifferential of a convex function. Even so, iteration with $R$
is still called the proximal point algorithm.  The second generalization
is the use of a generalized distance metric defined by the positive symmetric matrix $B$
which generalizes Eq. (\ref{lipschitz})
\begin{equation}
\label{gppa}
\|M(u) - M(v)\|_B \leq L \|u-v\|_B, \; \; \; \|u\|_B \equiv \sqrt{u^\top B u}.
\end{equation}
The generalized proximal point algorithm is 
\begin{equation*}
u^{(k+1)} = (B + \alpha F)^{-1} B u^{(k)},
\end{equation*}
which reduces to the proximal point algorithm if $B=I$.

The CPPD algorithm for solving
\begin{equation*}
\min_x, \max_\lambda \left\{ \lambda^\top A x - f^*(\lambda) \right\}
\end{equation*}
is an instance of the generalized proximal point algorithm,
where
\begin{equation*}
B = \left(
\begin{array}{cc}
T^{-1} & -A^\top \\
-A    & \Sigma^{-1}
\end{array}
\right), \; \; \alpha = 1,  \; \; 
F = \left(
\begin{array}{cc}
0 & A^\top \\
- A    & \partial f^*
\end{array}
\right), \; \; u = \left( \begin{array}{c} x\\ \lambda \end{array} \right).
\end{equation*}
Substituting these expressions into Eq. (\ref{gppa}) yields the CPPD algorithm as written
in Eq. (\ref{cppdint1}).

\section{CPPD step parameter computation}
\label{app:steps}

When it comes to realizing various instances of the CPPD algorithm, one of the more
complicated aspects of the implementation is computing the step matrix mappings $\Sigma$
and $T$ from the system matrix $A$. 
This topic is quite broad and problem dependent. We do not attempt to cover all possibilities.
We summarize only our experience in applying CPPD to optimization problems of interest
for X-ray tomographic applications.

Recall from Sec. \ref{sec:matrix-mapping} and \ref{app:fixedpoint} that 
$\Sigma$ and $T$ are positive symmetric matrices and they must satisfy the condition
that $B$ is also positive and symmetric, where
\begin{equation}
\label{cppdmetric}
B = \left(
\begin{array}{cc}
T^{-1} & -A^\top \\
-A    & \Sigma^{-1}
\end{array}
\right).
\end{equation}
We discuss three cases: scalar steps; a diagonal step matrix mapping; and a non-diagonal
step matrix mapping that preconditions the CPPD iteration.

\subsubsection*{Scalar step sizes}
The original CPPD paper \cite{chambolle2011first} employed
scalar steps
\begin{equation*}
\Sigma = \sigma I, \; \; T = \tau I,
\end{equation*}
where $\sigma$ and $\tau$ are positive real numbers, satisfying the condition
\begin{equation*}
\sigma \tau < 1/L^2, \; \; L \equiv \|A\|_2.
\end{equation*}
Note that this condition is equivalent to the condition that $M$ is a positive matrix.
The computation of $L$ can be time consuming, but it can be performed ahead of the CPPD
iteration and stored for the given system matrix $A$.
The CPPD algorithms presented in \cite{chambolle2011first} included $\sigma$ and $\tau$
fixed as a function of iteration number and Nesterov accelerated versions, where
$\sigma$ and $\tau$ vary with iteration number in such a way that can improved convergence rates.
Another reference that we found interesting was an adaptive scheme for balancing primal and dual
progress proposed by Goldstein {\it et al.}  \cite{goldstein2013adaptive}.

We have had some experience in implementing the CPPD Nesterov acceleration schemes
\cite{sidky2013first}; however, when we performed
empirical testing of various implementations of CPPD for
gradient sparsity regularization \cite{jorgensen2015little},
we found that fixing $\sigma$ and $\tau$ as a function of iteration number performed
just as efficiently if not more so than implementations involving various forms of acceleration.
We have since rarely employed CPPD with Nesterov acceleration
or other schemes \cite{goldstein2013adaptive}
for adapting $\sigma$ and $\tau$ as a function of iteration number.
One important piece of information that resulted from testing the adaptive scheme of
Goldstein {\it et al.} \cite{goldstein2013adaptive} is that the ratio of the step sizes
$\sigma$ and $\tau$ is an important tuning parameter as demonstrated in the results
shown in Sec. \ref{sec:results}. When using scalar steps, we set $\sigma$ and $\tau$
according to
\begin{equation*}
\sigma = \rho/L, \; \; \tau = 1/(\rho L),
\end{equation*}
and the step size ratio $\rho$ is varied to find the setting that leads to the most rapid convergence.
Note that these settings violate the strict inequality condition, but in practice we have
never encountered a situation where this setting has led to non-convergence CPPD iteration.

As a practical tip, it is helpful to test the computation of the scalars $\sigma$ and $\tau$,
because they can involve a complicated
power method implementation especially when the optimization problem of interest involves
a system matrix $A$ that is formed by stacking many matrices.
The step scalars $\sigma$ and $\tau$
are set to values that multiply to $1/L^2$. Any error that leads to $L$ significantly
less than $\|A\|_2$ will be discovered immediately, because it will lead to the product
$\sigma \tau$ being too large and divergent iteration.
On the other hand, $L$ larger than $\|A\|_2$ will not be easy to discover because the
CPPD algorithm will still converge to the solution of the optimization problem of interest.
It will, however, do so more slowly than the case where $L$ is computed correctly.
An easy test is to attempt the CPPD iteration with $L =a \|A\|_2$ where $0<a<1$. If convergent
behavior is observed for $a<0.8$, then this is an indication that there is a likely an error
in the computation of $L$. It is not a definitive test because the convergence condition
inequality is not tight. So it is actually
possible, in our experience, to have convergent CPPD iteration for  $0.8<a <1.0$.
A side benefit to performing this test is that a "safe" value $a<1.0$ can be empirically
discovered that leads to slightly faster CPPD convergence than the $a=1.0$ case.

\subsubsection*{Diagonal step matrices}
Pock and Chambolle proposed diagonal step matrices in \cite{Pock2011}
\begin{align*}
\Sigma &= \diag{(|A| 1)^{-1}}\\
T &= \diag{(|A|^\top 1)^{-1}},
\end{align*}
where the $\diag{\cdot}$ operator yields a diagonal matrix with the diagonal elements assigned
to the components of the vector in the argument. The absolute value $|\cdot|$ and inverse
$(\cdot)^{-1}$ operators are applied element-wise to the matrix and vector arguments.
As $A$ is and $m$ by $n$ matrix, the symbol $1$ is interpreted as a $n$-vector and
$m$-vector with all components set to 1 the $\Sigma$ and $T$ equation, respectively.
The diagonal step matrices have two advantages: they are more efficiently computed
than the scalar step sizes because they involve single matrix-vector products
between $|A|$ and 1 instead of repeated matrix-vector products required in the iterative
power method; and they perform preconditioning that may be significant for specific
problems. The efficiency of the diagonal step matrix computation enables application
of the CPPD algorithm to problems where the optimization problem of interest is changing
as a function of iteration number because it is feasible to re-compute the step matrices
at every iteration or every few iterations.  We have taken advantage of this efficiency
property in our mirrored convex/concave (MOCCA) algorithm which address optimization problems
that can be written as a combination of non-convex smooth and convex non-smooth functions \cite{Barber2016mocca}.

In implementing the diagonal step matrices, we have encountered a situation where
they can be more difficult to compute than the scalar step sizes.
If the system matrix $A$ is a product of matrices that may have negative matrix elements,
it can be difficult to obtain $|A|$. Note that this situation is not a problem for the scalar
step sizes because only the matrix $A$ is needed in the power method.
If $A$ is a product of matrices
\begin{equation*}
A = \Pi_i A_i,
\end{equation*}
it can be shown by repeated use of the Cauchy-Schwarz inequality that
\begin{equation*}
|A| \le \Pi_i |A_i|.
\end{equation*}
Taking the product of the absolute value of the matrix factors provides an upperbound on $|A|$,
and using this upperbound for the diagonal step matrices will yield conservative steps
that result in convergent iteration. 
Use of this upperbound is practical when computing $|A_i|$ is efficient and when it does not
result in too much loss of efficiency in the CPPD convergence. We have found this upperbound to be
useful in applying MOCCA to spectral CT image reconstruction \cite{Barber2016}, where the system
matrix involves a product of matrices including a preconditioning matrix that has both negative
and positive matrix elements.

The diagonal step matrices have proven useful, but we do not show examples of their use in this
article. When implementing diagonal step matrices it is important to make use of the the step
size ratio parameter $\rho$, discussed for scalar steps, in order to maximize algorithm efficiency.

\begin{algorithm}
\hrulefill
\begin{algorithmic}[1]
\State k=0
\While{ $k<K$ }
\State Initialize $u_k$
\State j=0
\While{ $j<N_\text{power}$ }
\State $u_k \gets A^\top A u_k$
\State i=0
\While{ $i<k-1$ }
\State $u_k \gets u_k - (u^\top_k u_i) u_i$
\State $i \gets i+1$
\EndWhile
\State $e_k \gets \|u_k\|_2$
\State $u_k \gets u_k/e_k$ \label{finalev}
\State $j \gets j+1$
\EndWhile
\State $k \gets k+1$
\EndWhile
\end{algorithmic}
\hrulefill
\caption{Modified power method for finding $K$ eigenvectors of
$A^\top A$. The parameter $N_\text{power}$ is the number of iterations taken for the
power method. The input to the algorithm is the matrix $A$ and the output are the eigenvectors $u_k$
and corresponding eigenvalues $e_k$.}
\label{alg:modifiedpower}
\end{algorithm}

\subsubsection*{Non-diagonal step matrices for preconditioning}
As discussed in Sec. \ref{sec:ndpc},
a non-diagonal preconditioning step matrix $T$ can be derived by approximating
the inverse of $(A^\top A)$
\begin{equation*}
\Sigma = \sigma I, \: \; T \approx (A^\top A)^{-1},
\end{equation*}
where the scalar $\sigma$ is determined after $T$ is specified
\begin{equation}
\label{TAA}
\sigma = 1/\|T A^\top A\|_2.
\end{equation}
A classic method of obtaining an approximate matrix inverse is to use
a truncated eigenvector expansion.
The matrix $(A^\top A)$ is symmetric and positive semi-definite,
and its eigenvector decomposition is
\begin{equation*}
A^\top A = U E U^\top,
\end{equation*}
where $U$ is an orthogonal matrix and $E$ is diagonal with non-negative
eigenvalues, $e_i$, on the diagonal; without loss of generality, the eigenvalues
are sorted from largest to smallest.
Using a truncated series, the matrix $T$ is expressed
\begin{equation}
\label{nondiagtau}
T = \frac{I}{e_K} + \sum_{i=1}^{K-1} u_i \left(\frac{1}{e_i} - \frac{1}{e_K} \right)
u_i^\top \approx (A^\top A)^{-1},
\end{equation}
where $I$ is the identity matrix; $K$ is the number of eigenvectors used in the expansion,
and the eigenvectors can be found by simple modification of the power method, described
in Algorithm \ref{alg:modifiedpower}.

\begin{figure}[h!]
\centering
\includegraphics[width=0.95\textwidth]{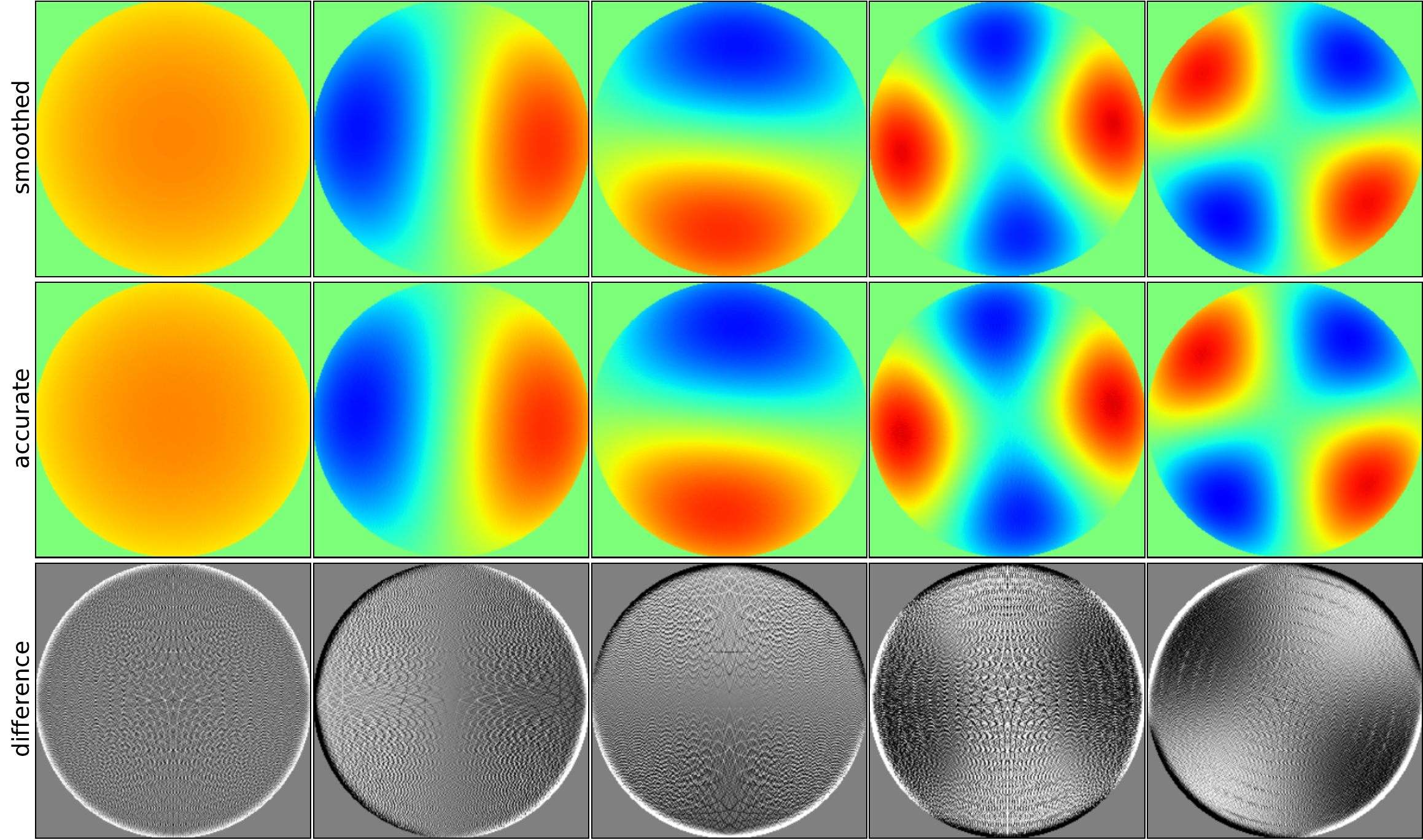}
\caption{\label{fig:evs}
The first five eigenvectors of $X^\top X$ ranked by singular value, going from left to right.
Top row shows the smoothed approximate eigenvectors, obtained by blurring the computed
eigenvector with a gaussian of 4-pixel width. Middle row contains the corresponding
numerically exact eigenvectors. Bottom row shows the difference between approximate
and accurate eigenvectors; these images show the Moire patterns that can enter
the image iterates of CPPD if the accurate eigenvectors are used to form the CPPD preconditioner.
The eigenvector images are shown in a color scale because it displays the eigenvector
structure clearly; the color scale, blue to red, spans [-0.01,0.01]. The gray scale, black to white,
for the difference plots in the bottom row spans [-0.0001,0.0001].}
\end{figure}
For the LSQ problem, presented in Sec. \ref{sec:LSQ},
the matrix $A$ is assigned to $X$, the discrete-to-discrete
approximation of the X-ray transform.
Implementation of non-diagonal preconditioning with $X$ requires some specific
understanding of the structure of $X$. The eigenvalues of $X^\top X$ vary
over orders of magnitude, see for example Ref. \cite{Jakob13}, leading to slow
convergence of first-order iterative algorithms. The eigenvalues of the first few
eigenvectors decreases strongly, so building an approximate inverse of $X^\top X$
using even the first few eigenvectors can increase the efficiency of the CPPD iteration
substantially. Also, the leading eigenvectors tend to dominate the CPPD image iterates at
low iteration number. Careful design of $T$ can improve the image quality of the images
at low iteration number in addition to improving convergence speed.

The non-negative and
relatively smooth nature of the sensing matrix
$X$ is such that, in general, the largest eigenvalue corresponds to an
eigenvector with no spatial nodes. As the eigenvalues decrease, more spatial nodes
are seen in the corresponding eigenvector.
Computing eigenvectors of $X^\top X$ for the configuration specified
in Sec. \ref{sec:LSQ} directly illustrates this trend as seen in Fig. \ref{fig:evs}. 
Due to the discretization of the X-ray transform there is a small
high-spatial frequency component present that can contaminate images
at low iteration number with Moire patterns.
For preconditioning, we
only want to account for the low-spatial frequency dependences and avoid
the high-frequency Moire patterns.
Thus, in designing the $T$ preconditioner, there is benefit to computing
the approximate matrix inverse from spatially smoothed eigenvectors.
To obtain such eigenvectors, the modified power method in Alg. \ref{alg:modifiedpower}
is employed with $A$ replaced by $X$, and the resulting eigenvector
are smoothed
\begin{equation}
u_k \gets Su_k
\end{equation}
where $S$ is a symmetric
matrix that performs discrete convolution with a smooth kernel.
The corresponding leading smoothed eigenvectors are shown in \ref{fig:evs} for
smoothing by a Gaussian with a width of 4 pixels. The smoothed eigenvectors
capture the low-frequency behavior of the accurate eigevectors,
and they can be used in Eq. \ref{nondiagtau} to
construct the matrix $T$. 

For the TV-LSQ problem 
where $A$ consists of the stacked matrix
\begin{equation*}
A = \left[ \begin{array}{c} X \\ \nu \nabla \end{array} \right],
\end{equation*}
the leading eigenvectors
of $SX^\top XS$ are also approximate eigenvectors of $A^\top A$.
The reason for this is that these leading eigenvectors have low spatial
frequency, and as a result they are all nearly in the null space of $\nabla$.
So, once again, we can use the smoothed eigenvectors of $X^\top X$ to form
non-diagonal $T$ using Eq. (\ref{nondiagtau}).
The $k$th eigenvalue is used to determine the parameter $\nu$
\begin{equation*}
\nu^2 = e_K/\|\nabla\|^2_2.
\end{equation*}
For TV-LSQ, smoothed eigenvectors and eigenvalues of $X$ are computed with the modified
power method. The matrix stacking combination parameter $\nu$ is then calculated. Finally, $\sigma$
is obtained from Eq. (\ref{TAA}).

The truncated eigenvector series in Eq. (\ref{nondiagtau}) has the advantage
that it applies to non-standard CT configurations such as sparse-view or
limited angular range sampling. When using standard scanning, an effective
preconditioner can be derived from the fact that $X^\top X$ is equivalent
to convolving with $1/r$ for $X$ representing continuous 2D parallel-beam projection \cite{natterer}.
This observation leads to the cone-filter preconditioner, see for example Ramani
and Fessler \cite{ramani2012splitting}. Even though the cone-filter is derived
for 2D circular, parallel-beam CT, it should still be effective for circular, fan-beam
CT with a discrete-to-discrete system matrix $X$.
Use of the cone-filter as $T$ is also demonstrated in Sec. \ref{sec:LSQ}.
Finally, we point out that the step size ratio parameter, $\rho$, can be used with the matrix mapping
steps. We did not explicitly include it in order to simplify the discussion.

Other examples of preconditioning with splitting methods applied in imaging have
been discussed in the context of CT in \cite{wang2019fast}. 
For magnetic resonance imaging (MRI), Ref. \cite{ong2019accelerating} presents a
preconditioner tailored to application of CPPD to MRI-based optimization problems.

\section{Convergence criteria}
\label{app:convergence}

Recall from Sec. \ref{sec:CPPDback} that the CPPD iteration solves the generic convex minimization
\begin{equation*}
x^\star = \argmin_x F(Ax).
\end{equation*}
Because we are interested in $F(\cdot)$ being possibly nonsmooth, useful convergence criteria
based directly on the objective function $F(Ax)$ can only be done in a problem specific manner.
To illustrate the issue consider $F$ to be one of three common cases: $(i)$ differentiable,
$(ii)$ $\ell_1$-norm, or $(iii)$ an indicator of a convex set.

$(i)$ If $F$ is differentiable, the magnitude of the gradient
\begin{equation*}
g = A^\top \nabla F(Ax)
\end{equation*}
provides a useful convergence criterion because $\|g\|_2$ approaches 0 as $x$ approaches the solution.
Also, if $F$ is locally quadratic about the solution, $\|g\|_2$ is linear in the distance between $x$
and the solution. The objective function value itself is not as useful because it
requires knowledge of the minimum value of $F(Ax)$
and it is not as sensitive as the objective function gradient.

$(ii)$ If $F$ is the $\ell_1$-norm, 
one can check that zero is in the subgradient, i.e. $0 \in \partial \|Ax\|_1$, but this condition is
not useful because zero is in the subgradient only when $Ax$ is exactly zero.
Otherwise, when $Ax \neq 0$ the gradient magnitude 
away from the solution is constant. Thus the gradient
does not yield information on proximity
of $x$ to the minimizer. The objective function itself, however, is informative in this case as long as
the function minimum is known ahead of time.

$(iii)$ If $F$ is an indicator function, i.e. the objective function is $\delta_S(Ax)$, neither the objective
function value nor its subgradient is informative. Instead one can employ a distance function that
yields the distance between $Ax$ and the closest point in the set $S$.

For any given $F$, it may be possible to employ some combination of a gradient, objective function value,
and distance function to provide useful convergence criteria.
What we would like to promote here instead are general convergence conditions that
can be used for any problem where the CPPD iteration applies.

\subsection*{First-order CPPD convergence criteria}
Convenient convergence conditions can be formulated from the intermediate saddle point problem
Eq. (\ref{saddle1}), which we repeat here
\begin{equation*}
\min_{x,y} \max_\lambda \left\{\phi(y) + \lambda^\top(Ax-y)\right\}.
\end{equation*}
The solution of this optimization not only involves finding the desired $x$ but also the 
dual variable $\lambda$ and the splitting variable $y$. Repeating the
conditions from Eqs. (\ref{cond1})-(\ref{cond3}),
the solution to the saddle problem is
\begin{equation*}
 A^\top \lambda = 0, \; \;\;
 \lambda =  \partial \phi(y), \; \;\;
 Ax - y = 0.
\end{equation*}
The middle equation is already used in deriving the CPPD saddle point problem;
recall that Eq. (\ref{saddle2}) only involves $x$ and $\lambda$. So for convergence we only
need to check
\begin{align}
r_\tau^{(k)} &\equiv A^\top \lambda^{(k)} \rightarrow 0, \label{cc1} \\
r_\sigma^{(k)} &\equiv Ax^{(k)} -y^{(k)} \rightarrow 0, \label{cc2} 
\end{align}
as $k \rightarrow \infty$. The equation $r_\tau = 0$ is called the transversality
condition ($\tau$ for $\tau$ransversality),
see for example Proposition 4.4.1 on pg. 336 of Ref. \cite{hiriart1993convex}.
The equation $r_\sigma=0$ is the feasibility condition, but we use a more specific
name for the magnitude of the left-hand side.
We call the quantity $\|r_\sigma\|_2$ the splitting gap ($\sigma$ for $\sigma$plitting)
because it is a meaure of the difference between the splitting variable $y$ and $Ax$.
The first condition can be checked immediately because $\lambda^{(k)}$ is available from
the CPPD iteration, but the second equation requires $y^{(k)}$ to be known. To obtain $y$,
an extra line can be introduced into the CPPD iteration.

[Make sure to clarify the difference between $y \in \partial \phi$ and $y=\partial \phi$,
 especially in the Legendre Transform relations ]

Recall the CPPD iteration Eq. (\ref{finalCPPD})
\begin{align}
x^{(k+1)} =& x^{(k)} - T A^\top \lambda^{(k)}, \notag \\
\bar{x}^{(k+1)} =& 2 x^{(k+1)}  - x^{(k)}, \notag  \\
\lambda^{(k+1)} =& \prox_{\sigma \phi^*} (\lambda^{(k)} + \sigma A \bar{x}^{(k+1)}).
\label{lastline}
\end{align}
The last line is equivalent to
\begin{equation}
\label{y1}
(I + \sigma \partial \phi^*) \lambda^{(k+1)}= \lambda^{(k)}+  \sigma A \bar{x}^{(k+1)},
\end{equation}
and there are two ways to proceed depending on whether it is easier
to compute $(i)$ $\prox_{\phi^*}$ or $(ii)$ $\prox_{\phi}$.

$(i)$ \underline{$\prox_{\phi^*}$ is simpler}:
Using the fact that $y = \partial \phi^*(\lambda)$,
Eq. (\ref{y1}) becomes
\begin{equation}
\label{y2}
\lambda^{(k+1)}  + \sigma  y^{(k+1)}= \lambda^{(k)}+  \sigma A \bar{x}^{(k+1)}.
\end{equation}
Solving for $y^{(k+1)}$ yields
\begin{equation*}
y^{(k+1)}= 
\frac{1}{\sigma}(\lambda^{(k)}-\lambda^{(k+1)}) +  A \bar{x}^{(k+1)}.
\end{equation*}
This update for $y$ can be appended to the CPPD iteration after Eq. (\ref{lastline}).
As a side note, using the $y$-update to directly write an update for the
splitting gap yields,
\begin{align*}
r_\sigma^{(k+1)} &=y^{(k+1)} - Ax^{(k+1)}\\
& = \frac{1}{\sigma}(\lambda^{(k)}-\lambda^{(k+1)}) +  A (2 x^{(k+1)} -x^{(k)}) - Ax^{(k+1)}\\
& = \frac{1}{\sigma}(\lambda^{(k)}-\lambda^{(k+1)}) +  A x^{(k+1)} - A x^{(k)} .
\end{align*}
This expression for the splitting gap is identical to what is called the ``dual residual''
in Ref. \cite{goldstein2013adaptive}.

$(ii)$ \underline{$\prox_{\phi}$ is simpler}:
In this case, the Moreau identity can be exploited to obtain $\prox_{\sigma \phi^*}$ 
\begin{equation}
\label{moreau_iden}
\prox_{\sigma \phi^*}(\lambda) + \sigma \prox_{\phi/\sigma}(\lambda/\sigma) = \lambda.
\end{equation}

Alternatively, and equivalently,
a modified version of the CPPD iteration can be directly derived. 
Using the fact that $y = \partial \phi^*(\lambda)$ and $\lambda = \partial \phi(y)$,
Eq. (\ref{y1}) becomes
\begin{equation*}
(\partial \phi + \sigma ) y^{(k+1)}= \lambda^{(k)}+  \sigma A \bar{x}^{(k+1)},
\end{equation*}
or
\begin{equation}
y^{(k+1)} = \prox_{\phi/\sigma } (\lambda^{(k)}/\sigma + A\bar{x}^{(k+1)}) \label{yline}.
\end{equation}
This line can be inserted before Eq. (\ref{lastline}) in the CPPD iteration, then
the variable $y$ is available to compute the convergence condition Eq. (\ref{cc2}).
But if this line is used to obtain $y$, it is not necessary to perform the additional
prox operation in Eq. (\ref{lastline}). The $\lambda$-update can obtained by solving
for $\lambda^{(k+1)}$ in Eq. (\ref{y2}).
This version of the CPPD algorithm becomes
\begin{align}
x^{(k+1)} =& x^{(k)} - T A^\top \lambda^{(k)}, \label{ycppd1}\\
\bar{x}^{(k+1)} =& 2 x^{(k+1)}  - x^{(k)}, \notag  \\
y^{(k+1)} =& \prox_{\phi/\sigma } (\lambda^{(k)}/\sigma + A\bar{x}^{(k+1)}), \notag \\
\lambda^{(k+1)}=& \lambda^{(k)} + \sigma (A \bar{x}^{(k+1)} -  y^{(k+1)}).
\label{ycppd2}
\end{align}
The CPPD algorithm as written in Eqs.~(\ref{ycppd1})-(\ref{ycppd2})
connects well with the convergence criteria, and, thus, reveals the trade-off between
the step parameters $\sigma$ and $T$. The $x$-update in Eq.~(\ref{ycppd1}) adjusts $x$
in such a way that reduces $|A^\top \lambda|$, and similarly,
the $\lambda$-update in Eq.~(\ref{ycppd2}) adjusts $\lambda$ so that $|Ax-y|$ is reduced. 
Choosing $T$ large increases
progress toward $r_\tau^{(k)} \rightarrow 0$, i.e. Eq.~(\ref{cc1}), and conversely
choosing $\sigma$ large increases progress toward $r_\sigma^{(k)} \rightarrow 0$,
i.e. Eq.~(\ref{cc2}).

\subsection*{The primal-dual gap}
In our first article that demonstrated the use of CPPD for optimization problem
prototyping for CT image reconstruction \cite{SidkyCP:12}, we proposed the conditional primal-dual (cPD)
gap as a convergence check. In that paper, we considered the convex optimization framework
put forth in Ref. \cite{chambolle2011first}, where the general minimization problem was
written
\begin{equation*}
x^\star = \argmin_x \left\{ \phi(Ax) +\eta(x) \right\},
\end{equation*}
and both $\phi$ and $\eta$ are convex functions.

In the simplfied framework discussed here this
gap is derived from the primal problem, involving minimization over $x$,
\begin{equation*}
x^\star = \argmin_x \phi(Ax),
\end{equation*}
and the dual problem, involving maximization over $\lambda$,
\begin{equation*}
\lambda^\star = \argmax_\lambda \left\{ -\delta(A^\top \lambda=0) -\phi^*(\lambda) \right\},
\end{equation*}
where the dual problem is obtained by performing the minimization over $x$ in Eq. (\ref{saddle2}).
For a given $x$ and $\lambda$ the primal-dual gap is the difference between these two objective
functions
\begin{equation*}
\text{PD}(x,\lambda) = \phi(Ax) + \delta(A^\top \lambda=0) + \phi^*(\lambda).
\end{equation*}
This function is positive, when $x$ and $\lambda$ are not solutions of their respective minimization
and maximization problems. When a solution to both is obtained,
\begin{equation*}
 \text{PD}(x^\star,\lambda^*) =0.
\end{equation*}
This condition, however, is not of much practical use for a convergence criterion because it involves
at least one indicator function, which can have infinite value for finite $x$ or $\lambda$.

A practical condition can be obtained by splitting both objective functions into bounded and indicator functions
\begin{equation*}
\phi(Ax) = \phi_\text{bp}(Ax) + \sum_i \delta ( Ax \in \mathcal{P}_i),
\; \; \; \phi^*(\lambda)= \phi^*_\text{bd}(\lambda) + \sum_j \delta ( \lambda \in \mathcal{D}_j).
\end{equation*}
[Make sure to explain all variations in indicator notation]
Progress toward satisfying each of the indicators can be measured by computing distance to each of the primal
convex constraint sets $\mathcal{P}_i$, dual constraint sets $\mathcal{D}_j$, and $A^\top \lambda$ from 0.
The remaining portion of the primal dual gap, which we called the conditional primal-dual gap, is
\begin{equation*}
\text{cPD}(x,\lambda) = \phi_\text{bp}(Ax) + \phi_\text{bd}^*(\lambda).
\end{equation*}
This quantity is bounded for finite $x$ and $\lambda$, and it approaches zero as $x \rightarrow x^\star$
and $\lambda -> \lambda^\star$. Because the indicators are taken away in forming cPD, this function can
have, in general, both positive and negative values.
So for checking convergence the absolute value of this function
should be checked.

Use of $|\text{cPD}(x,\lambda)|$ as a convergence criterion has trade-offs compared
with use of $r_\sigma$ and $r_\tau$.
The conditional primal-dual gap
does not require computation of the splitting variable $y$, but it can yield difficult to interpret
curves because cPD can oscillate about zero and its absolute value is used to check convergence.
In this article, we employ $r_\sigma$ and $r_\tau$ exclusively.

\section{Derivation of CPPD instances}
\label{app:instances}

In this appendix, the derivation for the specific instances of the CPPD algorithm are derived.
The generic convex optimization problem is
\begin{equation}
\label{genericOptK}
x^\star = \argmin_x \phi(Ax).
\end{equation}
The CPPD iteration for solving this equation is
\begin{align*}
x^{(k+1)} =& x^{(k)} - T A^\top \lambda^{(k)}, \\
\bar{x}^{(k+1)} =& 2 x^{(k+1)}  - x^{(k)},  \\
\lambda^{(k+1)} =& \prox_{\sigma \phi^*} (\lambda^{(k)} + \sigma A \bar{x}^{(k+1)}),\\
y^{(k+1)}=& 
\frac{1}{\sigma}(\lambda^{(k)}-\lambda^{(k+1)}) +  A \bar{x}^{(k+1)}.
\end{align*}
As explained in \ref{app:convergence}, the last line that updates the splitting
variable $y$ is included for computing the splitting gap, which is a convergence metric.
It is not strictly needed for obtaining the solution estimate.
In the case that it is easier to compute $\prox_\phi$ than $\prox_{\phi^*}$,
see Eq. (\ref{moreau_iden}).

\subsection{Least-squares (LSQ) minimization}
\label{sec:LSQ}

The first optmization problem is LSQ minimization
\begin{equation}
\label{LSQK}
f^\star = \argmin_f \frac{1}{2} \| Xf-g\|^2_2,
\end{equation} 
where $X$ is the discrete X-ray transform; $f$ represents the image expansion coefficients;
and $g$ is the projection data vector. Derivation of CPPD for this problem
is straightforward because the potential function $\phi$ is differentiable.
Making the following associations,
\begin{equation*}
A = X, \; \; \; x = f, \; \; \; \phi(y) = \frac{1}{2} \|y-g\|^2_2,
\end{equation*}
puts Eq. (\ref{LSQK}) in the form of Eq. (\ref{genericOptK}).
To obtain the CPPD updates, the Legendre transform of $\phi$ and the
$prox$ mappings are needed
\begin{equation}
\label{lsqprox}
\phi^*(y) = \frac{1}{2} \|y \|^2_2 + y^\top g, \; \; \;
\prox_{\sigma \phi^*}(\lambda) = \frac{ \lambda - \sigma g}{1+\sigma}.
\end{equation}
Substituting into the CPPD iteration equations, yields the CPPD-LSQ updates
\begin{align*}
f^{(k+1)} &= f^{(k)} - T X^\top \lambda^{(k)}, \\
\bar{f}^{(k+1)} &= 2f^{(k+1)} -f^{(k)},  \\
\lambda^{(k+1)} &= \left( \lambda^{(k)} + \sigma(X\bar{f}^{(k+1)} - g)\right)/
(1+\sigma).\\
y^{(k+1)} &=
\frac{1}{\sigma}(\lambda^{(k)}-\lambda^{(k+1)}) +  A \bar{f}^{(k+1)}.
\end{align*}

\subsection{Total variation penalized least-squares (TVLSQ)}

For the next algorithm instance, we derive the CPPD update formulas for TVLSQ.
\begin{equation}
\label{TVLSQ}
f^\star = \argmin_f \left\{ \frac{1}{2} \| Xf -g\|^2_2 + \beta
 \|D f \|_1 \right\},
\end{equation}
where $\beta$ is the scalar penalty parameter and $D$ is a finite-differencing
approximation of the image gradient.
Addressing this optimization problem,
while interesting in its own right, serves as a stepping stone to the algorithm
instance for the next problem, which is TV-constrained LSQ (TVCLSQ).
The associations for the TVLSQ optimization problem are

\begin{equation*}
A = \stack{X}{\nu D} , \; \; \;
x = f, \; \; \;
\phi(y) = \frac{1}{2} \|y_s-g\|^2_2 + (\beta/\nu) \| y_g \|_1,
\end{equation*}
where
\begin{equation*}
\nu = \|X\|_2/\|D\|_2;
\end{equation*}
$y_s$ and $y_g$ are the splitting variables for the X-ray transform and image gradient, respectively,
and
\begin{equation*}
y=\stack{y_s}{y_g},
\end{equation*}
where $s$ and $g$ stand for ``sinogram'' and ``gradient'', respectively.
We denote the two terms of the potential $\phi_s$ and $\phi_g$, where
\begin{equation*}
\phi_s(y_s) = \frac{1}{2} \|y_s-g\|^2_2, \; \; \;
\phi_g(y_g) =  (\beta/\nu) \| y_g \|_1.
\end{equation*}
Note that in fitting this optimization problem in the $\phi(Ax)$ form, it is necessary to ``stack'' the linear
transforms, $X$ and $\nu D$, to form the matrix $A$: The matrix $X$ is  $m$ by $n$ 
and in two dimensions $D$ is $2n$ by $n$; the number of columns of both matrices is the same number $n$, the number
of image pixels. The first $m$ rows of $A$ are the same as the $m$ rows of $X$, and rows $m+1$ through $m+2n$ of $A$
are the rows of $\nu D$. The constant $\nu$ is introduced as a factor multiplying $D$ and dividing $\beta$ so
that the matrices $X$ and $\nu D$ will have the same magnitude, i.e. largest singular value. This normalization
factor is useful because the step size parameters $\sigma$ and $T$ depend on $A$, and changing the units of $X$ or
$D$ will accordingly alter the CPPD iteration performance without this normalization.

Because the potential function $\phi$ separates into $\phi_s$ and $\phi_g$,
the proximal mapping separates also. This can be shown by explicitly writing the optimization problem corresponding
to $\prox_{\sigma \phi^*}$
\begin{align*}
\prox_{\sigma \phi^*}(\lambda) &= \argmin_{\lambda^\prime} \left\{
\sigma \phi^*(\lambda^\prime) + \frac{1}{2} \| \lambda^\prime - \lambda \|^2_2 \right\} \\
 &= \argmin_{\lambda_s^\prime,\lambda_g^\prime} \left\{
\sigma \phi^*_s(\lambda_s^\prime) + \frac{1}{2} \| \lambda_s^\prime - \lambda_s \|^2_2 +
\sigma \phi^*_g(\lambda_g^\prime) + \frac{1}{2} \| \lambda_g^\prime - \lambda_g \|^2_2 \right\} \\
 &= \stack{\prox_{\sigma \phi^*_s}(\lambda_s)}{\prox_{\sigma \phi^*_g}(\lambda_g)}.
\end{align*}
Accordingly, the CPPD iteration for TVLSQ takes the form
\begin{align*}
f^{(k+1)} =& f^{(k)} - T( X^\top \lambda_s^{(k)} + \nu D^\top \lambda_g^{(k)}),\\
\bar{f}^{(k+1)} =& 2 f^{(k+1)}  - f^{(k)}, \\
\lambda_s^{(k+1)} =& \prox_{\sigma \phi_s^* } (\lambda_s^{(k)}+\sigma  X\bar{f}^{(k+1)}), \\
\lambda_g^{(k+1)} =& \prox_{\sigma \phi_g^* } (\lambda_g^{(k)}+\sigma \nu D\bar{f}^{(k+1)}), \\
y_s^{(k+1)} =& \frac{1}{\sigma}(\lambda_s^{(k)}-\lambda_s^{(k+1)}) +  X \bar{f}^{(k+1)}, \\
y_g^{(k+1)} =& \frac{1}{\sigma}(\lambda_g^{(k)}-\lambda_g^{(k+1)}) +  \nu D \bar{f}^{(k+1)}.
\end{align*}
where $\prox_{\sigma \phi_s^*}$ and $\prox_{\sigma \phi_g^*}$ need to be evaluated.
The first proximal mapping is the same as Eq. (\ref{lsqprox}).

For the second proximal mapping, it is first necessary to obtain $\phi_g^*$.
Recall, $\phi_g = (\beta/\nu) \| y_g \|_1$, and from \ref{app:legendre} the convex conjugate is
\begin{equation*}
\phi_g^*(\lambda) = \delta ( \|\lambda\|_\infty \le \beta/ \nu ).
\end{equation*}
As $\phi_g^*$ is an indicator function, $\prox_{\sigma \phi_g^*}$ is projection
of the input variable $\lambda$ onto the convex set described by $\|\lambda\|_\infty \le \beta/\nu$,
which is implemented component-wise with the formula for the $i$th component given by
\begin{equation*}
\left[ \prox_{\sigma \phi_g^*}(\lambda_i)\right]_i = \begin{cases}
- \beta/\nu  & \lambda_i  \le -\beta/\nu \\
\lambda_i & -\beta/\nu < \lambda_i <\beta/\nu \\
\beta/\nu & \beta/\nu \le \lambda_i
\end{cases} .
\end{equation*}
Equivalently, this proximal mapping can be written as
\begin{equation}
\label{TVLSQprox}
\prox_{\sigma \phi_g^*}(\lambda) = \frac{(\beta/\nu) \lambda}{\max(\beta/\nu,|\lambda|)},
\end{equation}
where the $\max$ function operates component-wise on $|\lambda|$.
This latter form is slightly more convenient than the previous form, and it
avoids divide-by-zero since the minimum component value of the denominator is $\beta/\nu$.

Installing the specific proximal mappings, the CPPD update equations for TVLSQ become
\begin{align*}
f^{(k+1)} =& f^{(k)} - T( X^\top \lambda_s^{(k)} + \nu D^\top \lambda_g^{(k)}),\\
\bar{f}^{(k+1)} =& 2 f^{(k+1)}  - f^{(k)}, \\
\lambda_s^{(k+1)} =& \frac{\lambda_s^{(k)} + \sigma(X\bar{f}^{(k+1)} - g)}
{1+\sigma}, \\
\lambda^+_g =& \lambda_g^{(k)}+\sigma  \nu D\bar{f}^{(k+1)},\\
\lambda_g^{(k+1)} =& \frac{(\beta/\nu) \lambda^+_g }{\max(\beta/\nu,|\lambda^+_g|)},\\
y_s^{(k+1)} =& \frac{1}{\sigma}(\lambda_s^{(k)}-\lambda_s^{(k+1)}) +  X \bar{f}^{(k+1)}, \\
y_g^{(k+1)} =& \frac{1}{\sigma}(\lambda_g^{(k)}-\lambda_g^{(k+1)}) +  \nu D \bar{f}^{(k+1)}.
\end{align*}
The discussion on setting the step parameters in \ref{app:steps} apply to determining
the specific form of $T$ and $\sigma$ in these update equations.
In particular, we emphasize a few implementation issues. As discussed in \ref{app:steps},
whichever form of $T$ and $\sigma$ is selected, it is important to empirically tune
the step-size ratio parameter $\rho$ because convergence rates can vary by orders
of magnitude as a function of $\rho$. Also, in this presentation of the CPPD updates
for TVLSQ we have selected one particular method for stacking and normalization of the 
linear transforms $X$ and $D$. While other forms of transform normalization could be used,
it is important to do this in one form or another. If this is not done, algorithm performance
will vary with the physical units selected for formulating $X$ and $D$. For optimization problems
involving more than two linear transforms, the stacking and normlization generalizes in a
straight-forward manner.

\subsection{Total variation constrained least-squares (TVCLSQ)}

A related optimization problem to TVLSQ is total variation constrained least-squares (TVCLSQ),
which is formulated as
\begin{equation*}
f^\star = \argmin_f \left\{ \frac{1}{2} \| Xf -g\|^2_2 +\delta \left(
 \| D f \|_1 \le \gamma \right) \right\},
\end{equation*}
where the indicator function encodes a constraint that the image TV is bounded above by $\gamma$.
The TVCLSQ optimization problem is closely related to that of TVLSQ. In fact, if the constraint
parameter $\gamma$ of TVCLSQ and penalty parameter $\beta$ of TVLSQ are chosen appropriately, the
solutions will be identical.
The associations for putting the TVCLSQ optimization problem in the generic form of Eq. (\ref{genericOptK})
are
\begin{equation*}
A = \stack{X}{\nu D} , \; \; \;
x = f, \; \; \;
\phi(y) = \frac{1}{2} \|y_s-g\|^2_2 +\delta \left(
 \|y_g \|_1 \le \nu \gamma \right) ,
\end{equation*}
where $\nu$ is the normalization factor for $D$ as defined in the presentation of the TVLSQ problem;
and $y_s$ and $y_g$ are the splitting variables
for the X-ray transform and image gradient, respectively.
The TVCLSQ objective is split in two terms
\begin{align*}
\phi_s(y_s) &= \frac{1}{2} \|y_s-g\|^2_2, \\
\phi_{\rm gc}(y_g)&=  \delta \left( \|y_g \|_1 \le \nu \gamma \right) 
\end{align*}
Because $\phi_s$ is the same as it is in the TVLSQ example, we need only discuss $\phi_{\rm gc}$
(``gc'' stands for ``gradient constraint''), which
has been altered to reflect the constraint form.

The new proximal mapping needed for TVCLSQ is
$\prox_{\sigma F^*_{\rm gc}}$, which is related to $\prox{\phi_{\rm gc}}$ through the Moreau identity
\begin{align}
\prox_{\sigma \phi^*_{\rm gc}}(\lambda_g) &= \lambda_g - \sigma \prox_{F_{\rm gc}/\sigma}(\lambda_g/\sigma), \nonumber \\
&= \lambda_g - \sigma \proj(\lambda_g/\sigma, \|\lambda_g/\sigma\|_1 \le \nu \gamma), \nonumber  \\
&= \lambda_g - \proj(\lambda_g, \|\lambda_g\|_1 \le \nu \gamma \sigma). \label{moreauTVC}
\end{align}

We write the optimization problem corresonding to the projection (or proximal mapping)
on the right-hand side of Eq. (\ref{moreauTVC})
\begin{equation}
\label{constrainedOptL1}
\proj(\lambda_g, \|\lambda_g\|_1 \le \nu \gamma \sigma)=
\argmin_\lambda \left\{ \frac{1}{2} \| \lambda- \lambda_g\|^2_2 +
\delta( \| \lambda \|_1 \le \nu \gamma \sigma) \right\}.
\end{equation}
In evaluating this optimization problem, there are two cases:\\

$\|\lambda_g\|_1 \le \nu \gamma \sigma:$ When this inequality is satisfied,
\begin{equation*}
\proj(\lambda_g, \|\lambda_g\|_1 \le \nu \gamma \sigma) = \lambda_g.
\end{equation*}
The right-hand side of Eq. (\ref{moreauTVC}) is then $\lambda_g - \lambda_g$, and accordingly,
\begin{equation}
\label{zerocase}
{\rm if}\;\|\lambda_g\|_1<\nu \gamma \sigma, \;\;
\prox_{\sigma \phi^*_{\rm gc}}(\lambda_g) = 0.
\end{equation}

$\|\lambda_g\|_1 > \nu \gamma \sigma:$ In this case, $y_g$ needs to be projected onto the $\ell_1$-ball
of size $\nu \gamma \sigma$.
An efficient
sorting-based algorithm that can perform $\proj(\lambda_g, \|\lambda_g\|_1 \le \nu \gamma \sigma)$
is available in Ref. \cite{Duchi2008}. We also present an alternative that takes advantage of
one of the powerful aspects of splitting, i.e. separation the potential $\phi$ from large linear transforms.
In this non-trivial case where $\|\lambda_g\|_1 > \nu \gamma \sigma$ the constrained optimization
problem of Eq. (\ref{constrainedOptL1}) is equivalent to the following unconstrained optimization
problem
\begin{equation}
\label{unconstrainedOptL1}
\lambda_g^\prime= \argmin_\lambda \left\{ \frac{1}{2} \| \lambda - \lambda_g\|^2_2 +\beta \| \lambda \|_1 \right\},
\end{equation}
for an appropriate choice of penalty parameter $\beta$.
where we use a different solution notation, $\lambda_g^\prime$, because we do not know
{\it a priori} which value of $\beta$ yields
$\lambda_g^\prime = \proj(\lambda_g, \|\lambda_g\|_1 \le \nu \gamma \sigma).$
This optimization problem is the proximal mapping for the $\ell_1$-norm.
i
Equation (\ref{unconstrainedOptL1}) evaluates to a shrinkage operation
\begin{equation*}
\lambda_g^\prime = \prox_{\beta \| \cdot \|_1}(\lambda_g)= \sh(\lambda_g,\beta),
\end{equation*}
where $\beta$ is selected to be the value that
shrinks $\lambda_g$ until it satisfies
$\|\lambda_g\|_1 = \nu \gamma \sigma$.
Thus the equation for $\beta$ becomes
\begin{equation}
\label{impbeta}
\|\sh(\lambda_g,\beta) \|_1 - \nu \gamma\sigma =0.
\end{equation}
This problem can be readily solved numerically by any 1D root-finding algorithm,
e.g. the bisection method. The quantity $\|\sh(\lambda_g,\beta) \|_1$ is a monotonically decreasing
function of $\beta$, and we know that $\beta$ must be somewhere in the finite
interval $[0, \|\lambda_g\|_1]$. ($\beta=0$ is too small because this involves
no shrinking of $\lambda_g$, and $\beta=\|\lambda_g\|_1$ is too big because this value
would shrink $\lambda_g$ to zero magnitude.) This interval serves as input to
the root-finding algorithm.

After solving Eq. (\ref{impbeta}) for $\beta$, the projection operation
in Eq.~(\ref{moreauTVC}) is replaced by shrinkage
\begin{align*}
\prox_{\sigma \phi^*_{\rm gc}}(\lambda_g) &= \lambda_g - \sh(\lambda_g, \beta),\\
                                 &= \frac{\beta \lambda_g}{\max(\beta, |\lambda_g|)},
\end{align*}
where the second line is derived by considering the two cases where components
$[\lambda_g]_i<\beta$ and $[\lambda_g]_i \ge \beta$.

Interestingly, this proximal mapping is related to the proximal mapping derived for TVLSQ
in Eq.~(\ref{TVLSQprox}); the only difference is that $\beta$ is fixed for TVLSQ,
while here for TVCLSQ it depends on $\nu \gamma \sigma$ and the current iteration
of $f$ and $\lambda$.

\begin{algorithm}
\hrulefill
\begin{algorithmic}[1]
\State $f^{(k+1)} \gets f^{(k)} - T \left( X^\top \lambda_s^{(k)} + \nu D^\top \lambda_g^{(k)} \right)$
\label{cppdtvclsq1}
\State $\bar{f} \gets 2f^{(k+1)} -f^{(k)}$
\State $\lambda_s^{(k+1)} \gets \left( \lambda_s^{(k)} + \sigma(X\bar{f} - g)\right)/ (1+\sigma)$
\label{cppdtvclsq2}
\State $\lambda_g^+ \gets \lambda_g^{(k)} + \sigma \nu D \bar{f}$
\If {$\|\lambda_g^+\|_1 > \nu \gamma \sigma $}
\State $\beta^{(k+1)} \gets \solve(\beta,\|\sh(\lambda^+_g,\beta) \|_1 - \nu \gamma\sigma  =0) $
\label{cppdtvclsq3}
\State $\lambda_g^{(k+1)} \gets  \beta^{(k+1)}  \lambda_g^+ / \max \left( \beta^{(k+1)}, \left| \lambda_g^+ \right| \right)$
\label{cppdtvclsq4}
\Else
\State $\beta^{(k+1)} \gets 0$
\State $\lambda_g^{(k+1)} \gets  0$
\label{cppdtvclsq5}
\EndIf
\State $y_s^{(k+1)} \gets \frac{1}{\sigma}(\lambda_s^{(k)}-\lambda_s^{(k+1)}) +  X \bar{f}^{(k+1)}$
\State $y_g^{(k+1)} \gets \frac{1}{\sigma}(\lambda_g^{(k)}-\lambda_g^{(k+1)}) +  D \bar{f}^{(k+1)}$
\end{algorithmic}
\hrulefill
\caption{Pseudocode for the CPPD-TVCLSQ inner loop at iteration
$k$.}
\label{alg:CPPD-TVCLSQ}
\end{algorithm}
Combining all the necessary cases and elements for CPPD-TVCLSQ results in the pseudo-code
presented in Algorithm \ref{alg:CPPD-TVCLSQ}.
The function $\solve(a,b)$ means to solve equation $b$ for variable $a$ and return
the resulting value.


\begin{thebibliography}{10}
\expandafter\ifx\csname url\endcsname\relax
  \def\url#1{{\tt #1}}\fi
\expandafter\ifx\csname urlprefix\endcsname\relax\def\urlprefix{URL }\fi
\providecommand{\eprint}[2][]{\url{#2}}

\bibitem{chambolle2011first}
Chambolle A and Pock T 2011 {\em J. Math. Imag. Vis.\/} {\bf 40} 120--145

\bibitem{rockafellar1970convex}
Rockafellar R~T 1970 {\em Convex analysis\/} (Princeton university press,
  Princeton NJ)

\bibitem{SidkyCP:12}
Sidky E~Y, J{\o}rgensen J~H and Pan X 2012 {\em Phys. Med. Biol.\/} {\bf 57}
  3065--3091

\bibitem{boyd2011distributed}
Boyd S, Parikh N, Chu E, Peleato B and Eckstein J 2011 {\em Found. Trends Mach.
  Learn.\/} {\bf 3} 1--122

\bibitem{barber2024convergence}
Barber R~F and Sidky E~Y 2024 {\em J. Mach. Learn. Res.\/} {\bf 25} 1--46

\bibitem{he2012convergence}
He B and Yuan X 2012 {\em SIAM J. Imaging Sci.\/} {\bf 5} 119--149

\bibitem{Pock2011}
Pock T and Chambolle A 2011 Diagonal preconditioning for first order
  primal-dual algorithms in convex optimization {\em International Conference
  on Computer Vision ({ICCV} 2011)\/} pp 1762--1769

\bibitem{wirgin04}
Wirgin A 2004 {\em arXiv preprint math-ph/0401050\/}

\bibitem{paige1982lsqr}
Paige C and Saunders M~A 1982 {\em ACM Trans. Math. Soft\/} {\bf 8} 43--71

\bibitem{jensen2012implementation}
Jensen T~L, J{\o}rgensen J~H, Hansen P~C and Jensen S~H 2012 {\em BIT Numerical
  Mathematics\/} {\bf 52} 329--356

\bibitem{candes2006robust}
Cand{\`e}s E~J, Romberg J and Tao T 2006 {\em IEEE Trans. Inf. Theory\/} {\bf
  52} 489--509

\bibitem{donoho2006compressed}
Donoho D~L 2006 {\em IEEE Trans. Inf. Theory\/} {\bf 52} 1289--1306

\bibitem{jorgensen2015little}
J{\o}rgensen J~S and Sidky E~Y 2015 {\em Philos. Trans. Royal Soc. {A}\/} {\bf
  373} 20140387

\bibitem{lustig2007sparse}
Lustig M, Donoho D and Pauly J~M 2007 {\em Magn. Reson. Med.\/} {\bf 58}
  1182--1195

\bibitem{sidky2008image}
Sidky E~Y and Pan X 2008 {\em Phys. Med. Biol.\/} {\bf 53} 4777--4807

\bibitem{graff2015}
Graff C~G and Sidky E~Y 2015 {\em Appl. Opt.\/} {\bf 54} C23--C44

\bibitem{Reiser10}
Reiser I and Nishikawa R~M 2010 {\em Med. Phys.\/} {\bf 37} 1591--1600

\bibitem{Jakob13}
J{\o}rgensen J~S, Sidky E~Y and Pan X 2013 {\em IEEE Trans. Med. Imag.\/} {\bf
  32} 460--473

\bibitem{kahan1965pracniques}
Kahan W 1965 {\em Commun. ACM\/} {\bf 8} 40

\bibitem{Higham}
Higham N~J 2002 {\em Accuracy and stability of numerical algorithms, 2nd
  edition\/} (Society for Industrial and Applied Mathematics, Philadelphia, PA)

\bibitem{hiriart1993convex}
Hiriart-Urruty J~B and Lemar{\'e}chal C 1993 {\em Convex analysis and
  minimization algorithms I\/} (Springer, Berlin, Germany)

\bibitem{Beck09}
Beck A and Teboulle M 2009 {\em SIAM J. Imag. Sci.\/} {\bf 2} 183--202

\bibitem{Komodakis15}
Komodakis N and Pesquet J~C 2015 {\em IEEE Sig. Proc. Mag.\/} {\bf 32} 31--54

\bibitem{ryu2016primer}
Ryu E~K and Boyd S 2016 {\em Appl. Comput. Math\/} {\bf 15} 3--43

\bibitem{goldstein2013adaptive}
Goldstein T, Li M, Yuan X, Esser E and Baraniuk R 2013 {\em arXiv preprint
  arXiv:1305.0546\/}

\bibitem{sidky2013first}
Sidky E~Y, J{\o}rgensen J~S and Pan X 2013 {\em Med. Phys.\/} {\bf 40} 031115

\bibitem{Barber2016mocca}
Barber R~F and Sidky E~Y 2016 {\em J. Mach. Learn. Res.\/} {\bf 17} 144:1--51

\bibitem{Barber2016}
Barber R~F, Sidky E~Y, Schmidt T~G and Pan X 2016 {\em Phys. Med. Biol.\/} {\bf
  61} 3784--3818

\bibitem{natterer}
Natterer F 1986 {\em The mathematics of computerized tomography\/} (Society for
  Industrial and Applied Mathematics, Philadelphia, PA)

\bibitem{ramani2012splitting}
Ramani S and Fessler J~A 2012 {\em IEEE Trans. Med. Imag.\/} {\bf 31} 677--688

\bibitem{wang2019fast}
Wang T, Kudo H, Yamazaki F and Liu H 2019 {\em Phys. Med. Biol.\/} {\bf 64}
  art. no. 145006

\bibitem{ong2019accelerating}
Ong F, Uecker M and Lustig M 2019 {\em IEEE Trans. Med. Imaging\/} {\bf 39}
  1646--1654

\bibitem{Duchi2008}
Duchi J, Shalev-Shwartz S, Singer Y and Chandra T 2008 Efficient projections
  onto the $\ell_1$-ball for learning in high dimensions {\em Proceedings of
  the 25th international conference on Machine learning\/} (ACM) pp 272--279

\end{thebibliography}

\providecommand{\newblock}{}

\end{document}